\documentclass[12pt, oneside,reqno]{amsart}
\usepackage{amsmath, amsfonts, amssymb}
\usepackage{eucal}
\usepackage{mathrsfs}
\usepackage{latexsym}
\usepackage{cite}
\usepackage{bbm}        
\usepackage{bm}         
\usepackage{mathtools}
\usepackage{relsize}
\usepackage{etoolbox}
\usepackage{hyperref}
\usepackage{graphicx}
\usepackage{fancyhdr, float}
\usepackage{array}
\usepackage{esint}
\usepackage{tikz}
\usetikzlibrary{quotes,angles}
\usetikzlibrary{calc}
\usetikzlibrary{decorations.pathreplacing}
\usetikzlibrary{3d}
\usetikzlibrary{patterns,patterns.meta}
\usetikzlibrary {shapes.geometric}

\usepackage{etex,caption}
\usepackage{extarrows,pgfkeys}

\makeatletter
\def\@secnumfont{\bfseries}
\makeatletter

\makeatletter
\renewcommand\section{\@startsection{section}{1}{\z@}{-3.5ex \@plus -1ex \@minus -.2ex}{2ex \@plus .2ex}
{\centering\normalfont\large\scshape}}
\makeatother

\patchcmd{\subsection}{-.5em}{.5em}{}{}

\makeatletter
\renewcommand\subsubsection{\@startsection{subsubsection}{3}{\z@}{1ex \@plus .4ex \@minus .1ex}{-1em}
{\normalfont\footnotesize\bfseries}}
\makeatother

\usepackage[utf8]{inputenc}

\newcommand{\R}{\mathbb R}
\newcommand{\rn}{\mathbb R^n}
\newcommand{\sn}{S^{n-1}}
\newcommand{\bn}{B^n}

\newcommand{\kon}{\mathcal K_o^n}
\newcommand{\HH}{\mathcal{H}}

\newcommand{\bnu}{\pmb{\nu}}
\newcommand{\balpha}{\pmb{\alpha}}

\newcommand{\Hn}{\mathcal H^{n-1}}
\newcommand{\chara}[1]{{\mathbbm{1}_{#1}}}
\newcommand{\lkappa}{\ds\mathlarger{\kappa}}

\newcommand{\widebar}[3]{\mathrlap{\hspace{#2}\overline{\scalebox{#1}[1]{\phantom{\ensuremath{#3}}}}}\ensuremath{#3}}
\newcommand{\barV}{\widebar{.7}{2pt}{V}}

\newcommand{\ds}{\displaystyle}
\newcommand{\st}{\scriptstyle}
\newcommand{\sst}{\scriptscriptstyle}
\newcommand{\blb}{\raise.3ex\hbox{$\st \pmb \lbrack$}}
\newcommand{\sblb}{\raise.1ex\hbox{$\sst \pmb \lbrack$}}
\newcommand{\brb}{\raise.3ex\hbox{$\st \pmb \rbrack$}}
\newcommand{\sbrb}{\raise.1ex\hbox{$\sst \pmb \rbrack$}}

\newcommand{\bla}{\raise.2ex\hbox{$\st\pmb \langle$}}
\newcommand{\sbla}{\raise.1ex\hbox{$\sst\pmb \langle$}}
\newcommand{\bra}{\raise.2ex\hbox{$\st\pmb \rangle$}}
\newcommand{\sbra}{\raise.1ex\hbox{$\sst\pmb \rangle$}}
\newcommand{\blrb}{\raise.3ex\hbox{$\st \pmb | $}}
\newcommand{\sblrb}{\raise.1ex\hbox{$\sst \pmb | $}}
\newcommand{\wul}[1]{{\blb #1 \brb}}
\newcommand{\swul}[1]{{\sblb #1 \sbrb}}

\newcommand{\wt}{\widetilde}
\newcommand{\wts}{\,\wt{+}\,}
\newcommand{\psum}{\,{+\raisebox{-.5pt}{$_{\negthinspace\kern-1.5pt p}$}}\,}
\newcommand{\spsum}{\,{+\raisebox{-.5pt}{$_{\negthinspace\kern-1pt \sst p}$}}\,}
\newcommand{\qsum}[1]{\,{+_{\negthinspace\kern-2pt \lower -2pt \hbox{$_{_{#1}}$}}}\,}
\newcommand{\osum}{{+_{\negthinspace\kern-2pt {\rm{o}}}}\,}
\newcommand{\dpsum}{\,{\tilde+_{\negthinspace\kern-1pt p}}\,}
\newcommand{\dqsum}[1]{{\,\wt+_{\negthinspace\kern-1pt #1}}\,}

\newcommand{\lsub}[1]{\hskip -1.2pt\lower.4ex\hbox{$_{#1}$}}

\newcommand{\rhok}{\rho\hskip -1pt\lower.4ex\hbox{$_{K}$}}
\newcommand{\nuk}{\nu{\hskip -2pt\lower.2ex\hbox{$_{K}$}}}

\newcommand{\rk}{r{\hskip -2pt\lower.2ex\hbox{$_{K}$}}}
\newcommand{\srk}{r{\hskip -2pt\lower.2ex\hbox{$_{\sst K}$}}}
\newcommand{\thin}[1]{\negthinspace #1 \negthinspace}

\newcommand{\sbr}[1]{\item[\scriptsize\rm ({\bf #1})]}
\newcommand{\dr}[1]{{\rm ({\small \bf #1})}}
\newcommand{\fr}[1]{\noindent{\footnotesize {\bf #1.}}}

\newcommand{\vs}[1]{\vspace{#1 pt}}

\newcommand{\conv}{\operatorname{conv}}
\newcommand{\interior}{\operatorname{int}}

\makeatletter
\def\l@subsection{\@tocline{2}{0pt}{2.8pc}{5pc}{}}
\makeatother

\let\oldtocsection=\tocsection
\let\oldtocsubsection=\tocsubsection
\renewcommand{\tocsection}[2]{\hspace{0em}\oldtocsection{#1}{#2}}
\renewcommand{\tocsubsection}[2]{\small \oldtocsubsection{#1}{#2}}

\makeatletter
\def\@tocline#1#2#3#4#5#6#7{\relax
  \ifnum #1>\c@tocdepth 
  \else
    \par \addpenalty\@secpenalty\addvspace{#2}%
    \begingroup \hyphenpenalty\@M
    \@ifempty{#4}{%
      \@tempdima\csname r@tocindent\number#1\endcsname\relax
    }{%
      \@tempdima#4\relax
    }%
    \parindent\z@ \leftskip#3\relax \advance\leftskip\@tempdima\relax
    \rightskip\@pnumwidth plus4em \parfillskip-\@pnumwidth
    #5\leavevmode\hskip-\@tempdima #6\nobreak\relax
    \ifnum#1<0\hfill\else\hfill\fi\hbox to\@pnumwidth{\small\@tocpagenum{#7}}\par
    \nobreak
    \endgroup
  \fi}
\makeatother

\begin{document}

\title[Minkowski Problems]
{Minkowski Problems for Geometric Measures}

\author[Y. Huang]{Yong Huang}
\address{Institute of Mathematics,
Hunan University,
Changsha, 410082, China}
\email{huangyong@hnu.edu.cn}

\author[D. Yang]{Deane Yang}
\address{Department of Mathematics,
Courant Institute of Mathematical Sciences,
New York University,
251 Mercer Street,
New York, NY 10012, USA}
\email{deane.yang@courant.nyu.edu}

\author[G. Zhang]{Gaoyong Zhang}
\address{Department of Mathematics,
Courant Institute of Mathematical Sciences,
New York University,
251 Mercer Street,
New York, NY 10012, USA}
\email{gaoyong.zhang@courant.nyu.edu}

\subjclass{52A38}
\keywords{Convex body, polytope, Minkowski problem, logarithmic Minkowski problem,
dual Minkowski problem, $L_p$ Minkowski problem, quermassintegral, dual quermassintegral,
chord integral, surface area measure,
cone-volume measure, area measure, curvature measure, dual curvature measure, chord measure.}

\thanks{Research of Huang supported, in part, by NSFC Grants (12171144,12231006).
Research of Yang and Zhang supported, in part, by NSF Grant  DMS--2005875.}

\maketitle

\numberwithin{equation}{section}

\newtheorem{theo}{Theorem}[section]
\newtheorem{coro}[theo]{Corollary}
\newtheorem{lemm}[theo]{Lemma}
\newtheorem{defi}[theo]{Definition}
\newtheorem{prop}[theo]{Proposition}
\newtheorem{conj}[theo]{Conjecture}
\newtheorem{exam}[theo]{Example}
\newtheorem{prob}[theo]{Problem}
\newtheorem{rema}[theo]{Remark}
\newtheorem*{MET}{\footnotesize Minkowski Existence Theorem}

\theoremstyle{definition}
\newtheorem*{definition-non}{Definition}

\begin{abstract}
This paper describes the theory of Minkowski problems for geometric measures
in convex geometric analysis. The theory goes back to Minkowski and Aleksandrov
and has been developed extensively in recent years.
The paper surveys classical and new Minkowski problems studied in
convex geometry, PDEs, and harmonic analysis, and structured in a conceptual
framework of the Brunn-Minkowski theory, its extensions, and related subjects.
\end{abstract}

\tableofcontents

\section{Introduction}\

Convex geometric analysis studies --- from an analytic perspective ---
global geometric invariants and locally defined geometric measures of convex bodies in Euclidean space \cite{S14}.
Global geometric invariants, such as volume and surface area, quantify the size of a convex body,
while geometric measures characterize its shape.
The geometric measures extend the concept of pointwise curvature invariants of smooth hypersurfaces
to local invariants of arbitrary convex bodies. They are essential for studying nonsmooth
convex bodies such as polytopes. The most important geometric measures  arise as the differentials
of global geometric invariants viewed as
functionals on the space of convex bodies.

A geometric measure of a convex body is usually defined on or pushed to the unit sphere via either the radial or Gauss map.
A basic example is the surface area measure which is the differential of the volume of a convex body
and is a measure on the unit sphere of outer unit normals of the convex body.
A fundamental problem in convex geometric analysis is the characterization of geometric measures.
In other words, given a measure on the unit sphere, does it arise as a specific geometric measure of a convex body?
The characterization of surface area measure is known as the classical Minkowski problem,
 first studied by Minkowski for polytopes in 1897. For this reason,
the characterization of a geometric measure is called the Minkowski problem for that measure.

In this article we survey old and new results on Minkowski problems for a range of geometric measures.
Perhaps more importantly, we show that everything can be organized in a cohesive conceptual framework
that arises from the Brunn-Minkowski theory, its extensions, and related subjects.

\vspace{3pt}

During the period 1890--1909, Minkowski \cite{M1897, M1903}, building on the work of Brunn, created
the Brunn-Minkowski theory of convex bodies in the Euclidean space, which lies at
the core of convex geometry.  A fundamental result of the theory is the Brunn-Minkowski inequality,
which extends the classical isoperimetric inequality
and demonstrates the concavity of volume with respect to the vector addition in the Euclidean space.
In the 1930s, Aleksandrov and Fenchel \cite{Al2, FJ} studied more
general geometric invariants and proved the Aleksandrov-Fenchel inequalities for mixed volumes.
Since then, the study of sharp inequalities satisfied by geometric invariants of convex bodies
has been a major focus in convex geometry. Beyond that, the geometric inequalities are closely
linked to analytic inequalities for functions and measures on Euclidean space. The best known
example of this is the equivalence between the isoperimetric and the $L^1$ Sobolev inequality.
Gardner's Bulletin survey \cite{G02bull} explains the relationships between the Brunn-Minkowski
inequality and many other geometric and analytic inequalities. Lutwak's survey paper \cite{L93hcg}
describes a number of affine geometric inequalities of affine geometric invariants.
Schneider's book \cite{S14} is the definitive reference
for the foundations of the Brunn-Minkowski theory.

Lacking the notion of a surface area measure, Minkowski formulated the Minkowski problem in terms of
facet areas for polytopes and the Gauss curvature for smooth convex bodies \cite{M1897, M1903}.
He solved the problem for these classes of bodies.
Aleksandrov and Fenchel-Jessen unified Minkowski's results by defining surface
area measure for all convex bodies. They resolved the Minkowski problem for arbitrary convex bodies.

However, the study of geometric measures of convex bodies and the corresponding Minkowski
problems became a major focus of research in convex geometry only after Lutwak
introduced $L_p$ surface area measure and its $L_p$ Minkowski problem in 1993 \cite{L93jdg}.
In 2016, dual curvature measures were constructed in the paper \cite{HLYZ16acta}.
The study of geometric measures and Minkowski
problems has greatly enriched the field of convex geometry.

The Minkowski problem for a geometric measure is simply the following: Given a measure
on the unit sphere, what are the necessary and sufficient conditions that there exists
a convex body whose geometric measure is equal to the given measure? If such a convex
body exists, to what extent is it unique? In the setting of smooth convex bodies,
this can be stated in terms of solving a fully nonlinear partial differential equation
on the unit sphere.
Significant advances in solving Minkowski problems have been made in recent years for
not only geometric measures that arise in the classical Brunn-Minkowski theory,
but also those in the dual Brunn-Minkowski theory, the $L_p$ Brunn-Minkowski theory,
integral geometry, and other subjects.

\section{Elementary Minkowski Problems}\

Minkowski problems are natural geometric problems, related to the construction
of convex bodies
in the Euclidean space. To illustrate their geometric meaning, we describe three
Minkowski problems for polygons in the plane. This helps us understand better the
general Minkowski problems in any dimension.

\subsection{Minkowski Problem for Polygons} \label{mp-polygon}

More than a century ago, Minkowski \cite{M1897} studied the problem of finding a
convex polytope whose facets have prescribed surface areas
and outer unit normals. In the plane, the problem can be stated as follows:
\vs{4}

{\it Given positive numbers $l_1, \ldots, l_m$ and unit vectors $v_1, \ldots, v_m$  in $\R^2$,
construct a convex polygon $P$ in $\R^2$ so that its sides have lengths
$l_1, \ldots, l_m$ and corresponding outer normals $v_1,\ldots, v_m$.}
\vs{5}


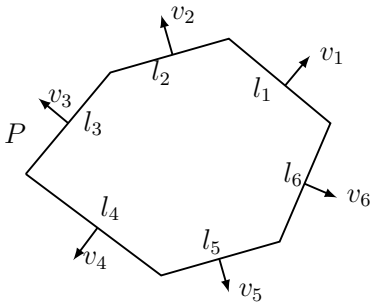
\begin{center}
\scalebox{0.9}{
\begin{tikzpicture}[scale=0.5]
	\coordinate (A) at (0,0) ;
	\coordinate (B) at (2.5,3);
	\coordinate (C) at (6,4);
	\coordinate (D) at (9,1.5);
	\coordinate (E) at (7.5,-2);
	\coordinate (F) at (4,-3);
	\coordinate (v3) at (0.37,2.24);
	\coordinate (v2) at (4,4.7);
	\coordinate (v1) at (8.4,3.5);
	\coordinate (v4) at (1.4,-2.55);
	\coordinate (v5) at (6,-3.5);
	\coordinate (v6) at (9.2,-0.7);
	
	\draw[thick] (A)--(B)--(C)--(D)--(E)--(F)--(A)node at (2.5,-1){$l_4$} node at (5.5,-2){$l_5$} node at (7.9,0){$l_6$} node at (2,1.5){$l_3$} node at (4,3){$l_2$} node at (7,2.5){$l_1$};	
	
	\coordinate[label={left:$P$}] (p) at (0.3,1.2);
	
    \draw[latex-,thick](v3) -- ($(A)!(v3)!(B)$) ;
    \coordinate[label={right:$v_3$}] (a) at (v3);
    \draw[latex-,thick](v2) -- ($(B)!(v2)!(C)$) ;
	\coordinate[label={right:$v_2$}] (a) at (v2);
    \draw[latex-,thick](v1) -- ($(C)!(v1)!(D)$) ;
	\coordinate[label={right:$v_1$}] (a) at (v1);
	\draw[latex-,thick](v4) -- ($(A)!(v4)!(F)$) ;
	\coordinate[label={right:$v_4$}] (a) at (v4);
	\draw[latex-,thick](v5) -- ($(F)!(v5)!(E)$) ;
	\coordinate[label={right:$v_5$}] (a) at (v5);
	\draw[latex-,thick](v6) -- ($(E)!(v6)!(D)$) ;
	\coordinate[label={right:$v_6$}] (a) at (v6);
\end{tikzpicture}
}
\end{center}

There is a discrete geometric measure on the unit circle associated with a convex polygon $P$:
\[
S(P, \cdot) = \sum_{i=1}^m l_i \, \delta_{v_i},
\]
where $\delta_{v_i}$ is the unit point mass (delta measure) at $v_i$.
This measure, concentrated at the normal vectors and with weights of the side lengths of a polygon,
has an analogue for an arbitrary convex body
in higher dimension, called the {\it surface area measure} of the convex body.

The {\it Minkowski problem} for convex polygons reads: Given a discrete measure $\mu$
on the unit circle, is there a convex polygon $P$ in $\R^2$ such that $S(P,\cdot) = \mu$?

This can be solved by a straightforward construction. One can assume that the unit
vectors $v_1, \ldots, v_m$ are listed in counterclockwise order. Obviously, the outer
unit normals of a convex polygon cannot lie on a half-circle. It is therefore necessary
to assume this for the unit vectors $v_1, \ldots, v_m$. We can now try to construct a polygon
as follows: Start with an oriented line segment $N_1$ of length $l_1$ that lies at a
counterclockwise right angle to $v_1$. Next, draw an oriented line segment $N_2$ of length $l_2$
that starts at the endpoint of $N_1$, and is counterclockwise orthogonal to $v_2$.
Repeat this for oriented line segments $N_3, \ldots, N_m$.  These oriented line segments will bound a convex
polygon if and only if the endpoint of $N_m$ is the starting point of $N_1$. It is
not hard to show that this holds if and only if $l_1v_1 + \cdots + l_mv_m = 0$.
Thus, the solution to the Minkowski problem for polygons is the following:

\vs{3}

{\it Let $l_1, \ldots, l_m$ be positive numbers  and
$v_1, \ldots, v_m$ be unit vectors in $\R^2$.
Then there exists a convex polygon $P$ in $\R^2$  with side lengths $l_1, \dots, l_m$
and corresponding outer unit normals $v_1,\ldots, v_m$ if and only if
$l_1v_1 + \cdots + l_mv_m = 0$
and the unit vectors
$v_1,\ldots, v_m$ do not all lie in a unit half-circle.
Moreover, the polygon $P$ is unique up to translation.}

\vs{3}

For a detailed elementary proof of this result, see Klain \cite{Klain}.

\subsection{Aleksandrov Problem for Polygons}\label{ap2}

Aleksandrov \cite{A42}
went on to define a new geometric measure, called {\it integral curvature}, one of the
curvature measures constructed later by Federer \cite{F59}. Aleksandrov proposed and solved
the Minkowski problem for integral curvature, now called the {\it Aleksandrov problem}.

For a convex polygon in $\R^2$ that contains the origin in its interior, the unit vectors from the origin to
the vertices are called the radial vectors of the polygon, and the angle at a vertex formed by
the normals of two adjacent sides is called a normal angle of the polygon.
The Aleksandrov problem can be stated as follows:
\vs{4}

{\it Given positive numbers $\alpha_1,\ldots, \alpha_m$ and unit vectors $u_1, \ldots, u_m$ in $\R^2$,
construct a convex polygon $P$ in $\R^2$ that contains the origin in its interior
with normal angles $\alpha_1, \ldots, \alpha_m$ and radial vectors
$u_1, \ldots, u_m$.
}


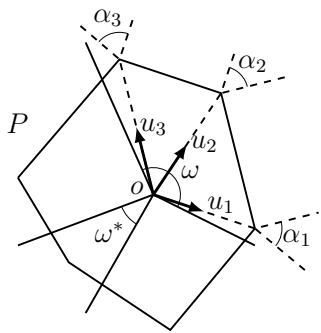
\vs{4}

\begin{center}
\scalebox{0.9}{
\begin{tikzpicture}[scale=0.5]
	\coordinate (A) at (0,0) ;
	\coordinate (B) at (3,3.5) ;
	\coordinate (C) at (6,2.5) ;
	\coordinate (D) at (7,-1.5) ;
	\coordinate (E) at (4.5,-4.5) ;
	\coordinate (F) at (1.5,-2.5) ;
	\coordinate (G) at (2,4) ;
	\coordinate (H) at (7,-2) ;
	\coordinate (I) at (0,-2) ;
	\coordinate (J) at (2,-4) ;
	\coordinate (O) at (4,-0.5) ;
	\coordinate (A1) at (-2,0) ;
	\coordinate (B1) at (3,35/6) ;
	\coordinate (B2) at (4,14/3) ;
	\coordinate (C2) at (8,10/3) ;
	\coordinate (C3) at (6,10.5) ;
	\coordinate (D3) at (9,-1.5) ;	
	\coordinate (D4) at (9.5,-1.5) ;
	\coordinate (E4) at (4.5,-7.5) ;	
	\coordinate (a1) at ($(B2)!(B)!(C2)$) ;
	\coordinate (a2) at ($(A1)!(B)!(B1)$) ;	
	\coordinate (c1) at ($(B2)!(C)!(C2)$) ;
	\coordinate (c2) at ($(C3)!(C)!(D3)$) ;	
	\coordinate (d1) at ($(C3)!(D)!(D3)$) ;
	\coordinate (d2) at ($(D4)!(D)!(E4)$) ;	
	\coordinate (u1) at (5.5,-1) ;	
	\coordinate (u2) at (5,1) ;	
	\coordinate (u3) at (3.5,1.5) ;	
	
	\draw[thick] (A)--(B)--(C)--(D)--(E)--(F)--(A);	
	\draw[thick,dashed] (O)--(B);
	\draw[thick,dashed] (O)--(C);
	\draw[thick,dashed] (O)--(D);
	
	\draw[-latex,very thick] (O)--(u1)node at (6,-0.8){$u_1$} ;
	\draw[-latex,very thick] (O)--(u2)node at (5.5,1){$u_2$} ;
	\draw[-latex,very thick] (O)--(u3)node at (4,1.5){$u_3$} ;

	\draw[thick] (H)--(O)--(G);
	\draw pic["$\ \ \ \omega$", draw=black,  angle eccentricity=1.2, angle radius=0.4cm]{angle=H--O--G};
	\draw[thick] (I)--(O)--(J);
	\draw pic["$\omega^*$", draw=black,  angle eccentricity=1.7, angle radius=0.5cm]{angle=I--O--J};
	
	\coordinate[label={left:$P$}] (p) at (0.6,1.6);
	\coordinate[label={left:$o$}] (p) at (4,-0.3);
	
	\draw[dashed,thick] (a2)--(B) -- (a1) ;
	\draw pic["$\alpha_3$", draw=black,  angle eccentricity=1.6, angle radius=0.4cm]{angle=a1--B--a2};	
	\draw[dashed,thick] (c2)--(C) -- (c1) ;
	\draw pic["$\alpha_2$", draw=black,  angle eccentricity=1.6, angle radius=0.4cm]{angle=c2--C--c1};	
	\draw[dashed,thick]  (d2)--(D) -- (d1) ;
	\draw pic["$\alpha_1$", draw=black,  angle eccentricity=1.6, angle radius=0.4cm]{angle=d2--D--d1};	
\end{tikzpicture}
}
\end{center}

There is a discrete geometric measure on the unit circle associated with this convex polygon,
\[
J(P, \cdot) = \sum_{i=1}^m \alpha_i \, \delta_{u_i}.
\]
This measure, concentrated at the radial vectors and with weights of the normal angles
of the polygon, was extended by Aleksandrov to an arbitrary convex body
 in any dimension, called the integral curvature of the convex body.

 The {\it Aleksandrov problem} for convex polygons reads: Given a discrete measure $\mu$
on the unit circle, is there a convex polygon $P$ in $\R^2$ that contains the origin in its interior
such that $J(P,\cdot) = \mu$?

Since the sum of normal angles of a convex polygon is equal to $2\pi$,
there is the necessary condition, $\alpha_1+\cdots+\alpha_m=2\pi$.
Also, the given unit vectors
 cannot all lie in a half circle. But these conditions are not sufficient.

Assume that the radial vectors $u_1, \ldots, u_m$ are given on the unit circle in a counterclockwise order.
Let $\omega \subset S^1$ be an angular sector that contains only the radial vectors $u_i, u_{i+1}, \ldots, u_j$,
where $1 \le i \le j \le m$.
Then for every normal vector $v$ at a boundary point of $P$
whose radial direction is inside $\omega$, there exists $u\in \omega$ such that $u\cdot v >0$, that is, $v$
is outside of the polar angle,
\[
\omega^* = \{v \in S^1 : v\cdot u \le 0 \text{ for all } u\in \omega\}.
\]
This implies that
\[
\alpha_i + \cdots + \alpha_j + |\omega^*| < 2\pi, \ \ \ \text{i.e.} \ \ \
\alpha_i + \cdots + \alpha_j < \pi + |\omega|.
\]
This condition has also to be satisfied by the normal angles. Aleksandrov showed that these necessary conditions are
also sufficient. The solution to the Aleksandrov problem can be stated as follows:

{\it If $\alpha_1,\ldots, \alpha_m$ are positive numbers and $u_1, \ldots, u_m$ are unit vectors in $\R^2$
in a counterclockwise order, then there exists a convex polygon $P$ in $\R^2$ containing the origin in its interior
so that it has normal angles
$\alpha_1, \ldots, \alpha_m$ and radial vectors $u_1,\ldots, u_m$ if and only if the unit vectors $u_1,\ldots, u_m$
do not all lie in a half circle,
\[
\alpha_1+\cdots+\alpha_m=2\pi, \ \ \text{and } \ \  \alpha_i+ \cdots + \alpha_j < \pi + |\omega|, \ \ 1\le i\le j \le m,
\]
for any angle $\omega$ that contains the radial vectors $u_i, \ldots, u_j$.
}

One could try to give a construction to the Aleksandrov problem in $\R^2$ similar to that of the Minkowski
problem described in the previous section. But it does not seem easy to have one.

\subsection{Logarithmic Minkowski Problem for Polygons}\label{stancu}

Given a convex polygon $P$ and an interior point $o \in P$, each side $\ell$ of $P$ uniquely determines
a triangle that has $\ell$ as a side and $o$ as a vertex. In addition to the side lengths and the normal
angles of a convex polygon, the areas of these triangles are also natural measurements of the polygon.
Prescribing the triangle areas of a polygon is called the {\it logarithmic Minkowski problem} for polygons,
which can be stated as follows:
\vs{4}

{\it Given positive numbers $a_1, \ldots, a_m$ and unit vectors $v_1, \ldots, v_m$  in $\R^2$,
construct a convex polygon $P$ in $\R^2$ with triangle areas
$a_1, \ldots, a_m$ and corresponding outer unit normals $v_1,\ldots, v_m$.
}

\begin{center}
\scalebox{0.9}{
\begin{tikzpicture}[scale=0.5]
	\coordinate (A) at (0,0) ;
	\coordinate (B) at (2.5,3);
	\coordinate (C) at (5,3.4);
	\coordinate (D) at (9.1,1.5);
	\coordinate (E) at (7,-2);
	\coordinate (F) at (3.5,-3);
	\coordinate (O) at (5,0.8);
	\coordinate (v4) at (0.4,2.3);
	\coordinate (v3) at (3.7,4.4);
	\coordinate (v2) at (7.6,3.6);
	\coordinate (v5) at (1,-2.5);
	\coordinate (v6) at (5.7,-3.6);
	\coordinate (v1) at (9,-1);
	
	\draw[thick] (A)--(B)--(C)--(D)--(E)--(F)--(A)		
	node at (6.5,0){$a_1$}
	node at (6,2){$a_2$}
	node at (4.4,2.2){$a_3$}
	node at (3,1.4){$a_4$}
	node at (3.4,-0.6){$a_5$}
	node at (5.3,-1){$a_6$}
	node at (4.5,0.4){$o$};	
	\draw[thick,dashed] (C)--(O)--(B);
	\draw[thick,dashed] (E)--(O)--(D);
	\draw[thick,dashed] (A)--(O)--(F);
	
	\coordinate[label={left:$P$}] (p) at (0.3,1.2);
	
    \draw[latex-,thick](v4) -- ($(A)!(v4)!(B)$) ;
	\coordinate[label={right:$v_4$}] (a) at (v4);
	\draw[latex-,thick](v3) -- ($(B)!(v3)!(C)$) ;
	\coordinate[label={right:$v_3$}] (a) at (v3);
	\draw[latex-,thick](v2) -- ($(C)!(v2)!(D)$) ;
	\coordinate[label={right:$v_2$}] (a) at (v2);
	\draw[latex-,thick](v5) -- ($(A)!(v5)!(F)$) ;
	\coordinate[label={right:$v_5$}] (a) at (v5);
	\draw[latex-,thick](v6) -- ($(F)!(v6)!(E)$) ;
	\coordinate[label={right:$v_6$}] (a) at (v6);
	\draw[latex-,thick](v1) -- ($(E)!(v1)!(D)$) ;
	\coordinate[label={right:$v_1$}] (a) at (v1);
	
\end{tikzpicture}
}
\end{center}

There is another discrete geometric measure on the unit circle associated with a polygon $P$,
\[
V(P, \cdot) = \sum_{i=1}^m a_i \, \delta_{v_i},
\]
which concentrates at the unit normal vectors and has weights of triangle areas.
The total measure of $V(P,\cdot)$ is obviously the area of $P$, denoted by $V(P)$.
This measure can be defined
for an arbitrary convex body in higher dimensions,
called the {\it cone-volume measure} of the convex body.

The logarithmic Minkowski problem
 for convex polygons reads: Given a discrete measure $\mu$
on the unit circle, is there a convex polygon $P$ in $\R^2$ such that $V(P,\cdot) = \mu$?

The logarithmic Minkowski problem is much more difficult to solve than the Minkowski problem
and the Aleksandrov problem. In fact, even the planar discrete case, as stated above,
has not been solved completely.

If the prescribed measure is symmetric, then there can be an obstruction. The simplest example
is if the prescribed unit normals are $v_1, v_2, v_3=-v_1, v_4=-v_2$ and the corresponding
prescribed areas are $a_1, a_2, a_3=a_1, a_4=a_2$. It is clear that for such a measure,
the polygon must be a parallelogram centered at the origin. This implies that $a_1 = a_2$.
Therefore, if the prescribed measure has $a_1 \ne a_2$, then there is no solution. On the
other hand, if $a_1 = a_2$ and $P$ is a solution, then uniqueness does not hold, because
the polygon $\phi P$, where $\phi$ is a linear transformation such that $\phi v_1 = bv_1$ and
 $\phi v_2 = b^{-1}v_2$, $b>0$, is also a solution.

More generally, if $P$ is an $m$-sided polygon such that, for some $3 \le q \le m-1$,
$v_q = -v_1$ and $a_q = a_1$ and $Q$ is the convex hull of the edges normal to $v_1, v_q$,
then it is not difficult to show that
\[
  a_1 + a_q = \frac{1}{2}V(Q) \le \frac{1}{2}V(P),
\]
where equality holds if and only if $P$ is a parallelogram. This gives a necessary condition
for the prescribed measure to be the cone-volume measure of a polygon.

For origin-symmetric polygons, Stancu \cite{St02} showed, using a heat flow argument,
that this is also a sufficient condition and proved the following:
\vs{2}

{\it  Let $a_1, \ldots, a_k$ be positive numbers and $v_1, \ldots, v_k, -v_1, \ldots, -v_k$
be unit vectors ordered counterclockwise in $\R^2$. Then there exists an origin-symmetric convex
polygon $P$ with triangle areas $a_1, \ldots, a_k, a_1, \ldots, a_k$ and corresponding unit normal vectors
$v_1, \ldots, v_k, -v_1, \ldots, -v_k$
if and only if
\begin{align*}
a_i &< \frac12 (a_1+\cdots+ a_k), \ \ \ i=1,\ldots, k,  \ \ \ \text{when } k>2, \ \ \text{or} \\
a_1 &=a_2,  \ \ \ \text{ when } k=2.
\end{align*}
}
Partial results were also obtained in the general case, which remains open in $\R^2$.

\section{Basic Concepts of Convex Bodies} \label{BMt}\

To make the article readable for people from different fields, we try to provide basic concepts and
facts in this and later sections. Without affecting the reading, one can jump to the next section
and easily look back for a concept if necessary.

\subsubsection{Basic Notations}
Let $\rn$ be the $n$-dimensional Euclidean space, and $\bn$ the $n$-dimensional unit ball
in $\rn$ centered at the origin. Sometimes, we also write $\bn$ as $B$.
The boundary of $\bn$ is the $(n-1)$-dimensional unit sphere
denoted by $\sn$. Let $G_{n,i}$ denote the Grassmann manifold of $i$-dimensional subspaces in $\rn$.
For $x,y\in \rn$, let $x\cdot y$ be the Euclidean inner product
and $|x|$ the Euclidean norm.

The integration with respect to the Lebesgue measure in $\rn$ or on a hypersurface in $\rn$  is
denoted by $dx$.
The integration with respect to the normalized Haar measure on $G_{n,i}$ is denoted by $d\xi$.
Let $\mathcal H^{k}$ be the $k$-dimensional Hausdorff measure in $\rn$.

Denote by $C^+(\sn)$ the set of positive continuous functions on $\sn$, and by $C^+_e(\sn)$
the set of positive continuous even functions on $\sn$.
For a Borel measure $\mu$ on $\sn$, denote its total mass by $|\mu|$.

\subsubsection{Convex bodies}
A {\it convex body} $K$ in $\rn$ is a compact convex subset with nonempty interior. Denote by
$\interior K$ the set of interior points of $K$, and $\partial K$ the boundary of $K$.
Let $\mathcal K_o^n$ be the set of all convex bodies that contain the origin in their interiors.
We often call an origin-symmetric convex body a symmetric convex body.

The {\it polar body} of $K\in \kon$ is defined to be
\[
K^* = \{ x \in \rn : x\cdot y \le 1 \text{ for all } y\in K\}.
\]

Denote by $V(K)$ the volume of $K$ and by $S(K)$
the surface area of $\partial K$ that is the Hausdorff measure $\mathcal H^{n-1}(\partial K)$.
Let $\omega_n=\frac{\pi^{n/2}}{\Gamma(1+n/2)}$ be the volume of $\bn$.  The surface area of $\sn$
is $\frac{2\pi^{n/2}}{\Gamma(n/2)}=n\omega_n$.

A convex body $P$ is an $n$-dimensional {\it polytope} if it is the convex hull of
finitely many points $x_1, \ldots, x_m $ in $\rn$.

\subsubsection{Homogeneous Functions and Monge-Amp\`ere Type Equations on $\sn$}
For a twice-differentiable function $h$ in $\rn$, let $\nabla h$
and $\nabla^2 h$ be the gradient and Hessian in $\rn$,
respectively. If $h$ is homogeneous of degree $1$, then $\nabla h$ is homogeneous of degree $0$, and
\[
y\cdot \nabla h(y) = h(y), \quad \nabla^2 h(y)y=0, \ \ \ y\in \rn.
\]
Thus, any $v\in \sn$ is a unit eigenvector of $\nabla^2 h(v)$ corresponding to its zero eigenvalue.
Then the restriction $\nabla^2 h(v)|_{v^\perp}$ is a linear transformation in $v^\perp$ whose matrix
representation under an orthonormal basis of $v^\perp$ is $(\nabla_{ij} h + h\delta_{ij})$,
where $(\nabla_{ij} h)$ is the spherical Hessian with respect to the orthonormal basis and $\delta_{ij}$ is the
Kronecker delta. If $h$ is also a convex function, then $(\nabla_{ij} h + h\delta_{ij})$ is positive
semi-definite.

Most of the partial differential
equations associated with Minkowski problems are Monge-Amp\`ere type equations on $\sn$ which can
be stated as follows: Given a nonnegative continuous function $f$ on $\sn$,
find a function $h$ on the unit sphere $\sn$ such that
\begin{equation}\label{monge}
F(h, \nabla h)\,\det \big(\nabla_{ij} h + h\delta_{ij}\big) =  f.
\end{equation}

\subsubsection{Support Function}
A convex body has two basic functional descriptions, namely, its support
and radial functions.

The {\it support function} $h_K : \rn \to \R$  of a convex body $K$ in $\rn$ is defined by
\[
h_K(y) = \max\{x\cdot y : x\in K\}, \ \ \ y\in \rn.
\]
The support function is convex and homogeneous of degree 1.

For each $v\in \sn$, the hyperplane
\begin{equation*}
H_K(v) = \{x \in \rn : x\cdot v = h_K(v)\}
\end{equation*}
is called the {\it supporting hyperplane to $K$ with unit normal $v$}, which extends
the concept of a tangent hyperplane.
Geometrically, the value $h_K(v)$, $v\in\sn$, is the signed distance from the origin to the support hyperplane
$H_K(v)$.

\subsubsection{Radial Function}

The {\it radial function} $\rhok : \rn\setminus\{0\} \to (0,\infty)$ of a convex body $K$
in $\mathcal K_o^n$ is defined by
\[
\rhok(x) = \max\{t > 0 : tx\in K\}, \ \ \ x\in \rn.
\]
It is homogeneous of degree $-1$.
If $u \in \sn$, then $\rhok(u)$ is the radial distance from the origin to the boundary
$\partial K$ in direction $u$. When $x\in \partial K$, $\rhok(x)=1$.

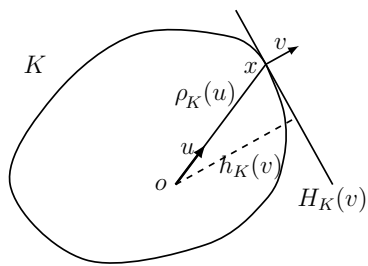
\vs{3}


\begin{tabular}{cl}
  \begin{tabular}{c} \hspace{20pt}
\scalebox{0.8}{
\begin{tikzpicture}[scale=0.5]
\coordinate (A1) at (-1.3,5);
\coordinate (A2) at (3,4);
\coordinate (A3) at (3.6,0.9);
\coordinate (A4) at (3,-0.5);
\coordinate (A5) at (1,-2);
\coordinate (A6) at (-3,-2.5);
\coordinate (A7) at (-5.5,0);
\coordinate (A8) at (-4,3);

\draw[thick]  plot [smooth,tension=0.6]
coordinates {(A2)(3.6,2.4) (A3) (A4) (A5) (A6) (A7) (A8) (A1) (1.9,4.9)(A2)};
\draw[thick](2,5.8)--(A2)--(47/9,0) node at (47/9,-0.5){$H_K(v)$};
\draw[thick](0,0)--(A2) node at (2.5,3.9){$x$};
\draw[-latex,thick]  (A2) -- ($(0,12)!(A2)!(20/3,0)$) ;
\draw[-latex,very thick]  (0,0) -- (1,4/3) node at (0.4,1.2){$u$} node at (1,3){$\rhok(u)$};
\draw[dashed,thick]  (0,0) -- ($(2,5.8)!(0,0)!(47/9,0)$) node at (2.5,0.6){$h_K(v)$};

\coordinate[label={left:$K$}] (k) at (-4,4);
\coordinate[label={left:$o$}] (o) at (0,0);
\coordinate[label={left:$v$}] (v) at (4,4.7);
\end{tikzpicture}
}
\end{tabular}
 &\begin{tabular}{ll}
    support function  &$h_K(v) = x\cdot v$ \\
    radial function   &$\rhok(u)= |x|$ \\
    radial map        &$r_K(u) =x$ \\
    Gauss map         &$\nu_K(x) = v$ \\
    radial Gauss map  &$\alpha_K(u) = v$
  \end{tabular} \\
 \end{tabular}


The convexity of $h_K$ implies that
$h_K$ is twice-differentiable almost everywhere in $\rn$, that is,
$h_K$ has the Taylor expansion at almost all $y\in \rn$,
\[
h_K(y+z) =h_K(y) + z\cdot\nabla h_K(y) +\frac12 z\cdot \big(\nabla^2 h_K(y)z\big) + o(|z|^2).
\]
The Hessian $\nabla^2h_K(y)$ is positive semi-definite.


\subsubsection{Curvature Function}
The convexity of $h_K$ implies that the restriction $\nabla^2 h_K(y)|_{y^\perp}$
is a positive semi-definite  linear transformation on  $y^\perp$. If $v \in \sn$
and $h_K$ is twice differentiable at $v$, then
\begin{equation}\label{0.18}
f_K(v) = \det\big(\nabla^2 h_K(v)|_{v^\perp}\big)
\end{equation}
is called the {\it curvature function of $K$}. If $h_K$ is $C^2$, then
\[
f_K(v) = \det \big(\nabla_{ij} h_K(v) + h_K(v)\delta_{ij}\big).
\]
If $f_K(v) > 0$, then $f_K(v)$ is the reciprocal of the Gauss curvature at the boundary point of $K$
with outer unit normal $v$.

\subsubsection{Wulff Shape and Convex Hull} A basic construction of a convex body from a function is the Wulff shape.
Given a positive continuous function $h$ on $\sn$, a convex body $\wul{h}$, called the {\it Wulff shape}
of $h$, is defined by
\[
\wul{h} = \{x\in \rn : x\cdot v \le h(v) \text{ for all } v\in \sn\}.
\]
It contains the origin in its interior. Intuitively, it is the largest convex body whose support
function is nowhere greater than $h$. Obviously,
\[
K=\wul{h_K}, \ \ \ \text{for  } K\in\kon.
\]

Another basic construction of a convex body from a function is the convex hull.
Given a positive continuous function $\rho$ on $\sn$, a convex body $\bla\rho\bra$,
called the {\it convex hull}
of $\rho$, is defined as the convex hull of the points $\rho(u)u$, $u\in \sn$, i.e.,
\begin{equation}\label{0.19}
\bla\rho\bra = \conv\{\rho(u)u : u\in \sn\}.
\end{equation}
It is also a convex body that contains the origin in its interior. It is the smallest
convex body whose radial function is nowhere less than $\rho$, and
\[
K=\bla\rhok\bra, \ \ \ \text{for }  K\in \kon.
\]

Wulff shape, convex hull, and polar are related by the following formula (see, for example, \cite{HLYZ16acta})
\[
\wul{1/\rho} = \bla\rho\bra^*.
\]

\subsubsection{Hausdorff Metric}

Convergence of sequences of convex bodies can be defined using the {\it Hausdorff distance},
\[
  d_H(K,L) = \|h_K - h_L\|_\infty.
\]

A basic fact of convergence of sequences of convex bodies is the Blaschke selection theorem:
{\it If a sequence $\{K_i\}$ of convex bodies in $\rn$ is uniformly bounded, then it has a subsequence
that converges to a compact convex set in $\rn$.}
It is a version of the classical Bolzano-Weierstrass theorem in $\rn$.

\vs{3}

\subsubsection{Radial and Gauss Maps} \label{rGm}
There are two fundamental maps between the boundary of a convex body and the unit sphere,
known as the radial and Gauss maps.

The {\it radial map} $\rk : \sn \to \partial K$ of a convex body $K$ in $\mathcal K_o^n$ is defined by
\[
\rk(u) = \rhok(u)u, \ \ \ u\in \sn.
\]
Geometrically, $r_K$ maps a radial direction $u$ to the boundary point $\rhok(u)u$ in that direction.
The radial map is a homeomorphism between the radial sphere and the boundary of $K$.
The boundary $\partial K$ is differentiable at $\rk(u)$ if and only if $\rhok$ is
differentiable at $u$. The gradient $\nabla \rhok(u)$ in $\rn$ is an inward pointing
normal vector of $\partial K$ at $\rk(u)$.

A point $x \in \partial K$ is called a regular or smooth point if $\partial K$ is differentiable at
$x$. Denote by $\partial'K$ the set of regular points of $\partial K$.
The boundary $\partial K$ has a unique outer unit normal $\nuk(x)$ at a regular point $x\in \partial'K$.
This defines the {\it Gauss map},
$
\nuk : \partial' K \to \sn.
$
That is, the Gauss map $\nuk$ maps a regular point to the unique outer unit normal at the point.
The Gauss map can be expressed using the gradient of the radial function,
\[
\nuk(x) = -\frac{\nabla \rhok(x)}{|\nabla \rhok(x)|}, \quad x\in \partial'K.
\]
Since the boundary $\partial K$ is differentiable almost everywhere
with respect to the Hausdorff measure $\Hn$, the Gauss map is defined almost everywhere on
$\partial K$.

The {\it radial Gauss map} $\alpha_K$, which maps a unit radial vector to the corresponding
outer unit normal vector, is defined almost everywhere on the radial sphere $\sn$ by
\[
\alpha_K = \nuk \circ r_K,
\]
that is,
\[
\alpha_K(u) = -\frac{\nabla \rhok(u)}{|\nabla \rhok(u)|}\text{ if } u\in r_K^{-1}(\partial'K).
\]

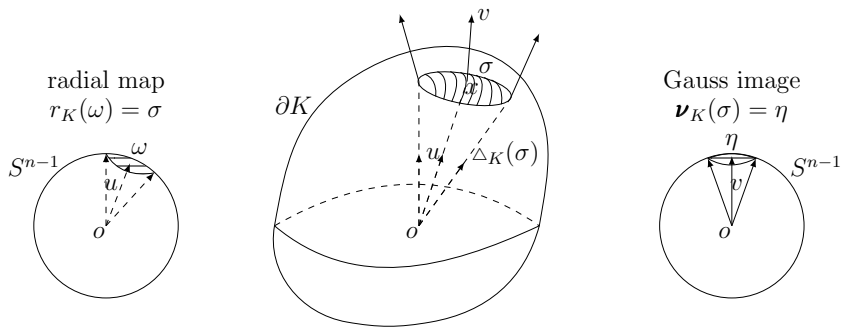
\begin{center}
\scalebox{0.8}{
		\begin{tikzpicture}[scale=0.4]
		\draw (-13,0) circle [radius=3cm]
		node at (-13.3,-0.3){$o$}
		node at (-16,2.5) {$\sn$}
		node at (-13,6) {radial\ map}
		node at (-13,4.8){$r_K(\omega)=\sigma$}
		node at (-11.6,3.2){$\omega$};
		\draw [-latex,dashed](-13,0)--(-13,3) node at (-12.8,1.7) {$u$};
		\draw [-latex,dashed](-13,0)--(-11,2.2);
		\draw [-latex,dashed](-13,0)--(-12,2.6);
		\draw (-13,3)[out=-60,in=190] to (-11,2.2) ;
		\fill[pattern={Lines[angle=30,yshift=.5pt]},pattern color=black]
		(-13,3)[out=-60,in=190] to (-11,2.2) [out=140,in=0] to (-13,3);

		\draw (13,0) circle [radius=3cm]
		node at (12.7,-0.3){$o$}
		node at (16.5,2.5) {$\sn$}
		node at (13,6) {Gauss\ image}
		node at (13,4.8){$\bnu_K(\sigma)=\eta$}
		node at (13,3.5){$\eta$}
		node at (-5.2,5) {$\partial K$};
		\draw [-latex](13,0)--(13,2.9)
		node at (13.2,1.7) {$v$} ;
		\draw [-latex](13,0)--(14,2.8);
		\draw [-latex](13,0)--(12,2.8);
		\filldraw[pattern={Lines[angle=40,yshift=.5pt]},pattern color=black]
		(12,2.8) [out=-30,in=-155] to (14,2.8) [out=160,in=20] to (12,2.8);
		
		\draw (-6,0)[out=-40,in=-155] to (5,0)
		node at (-0.3,-0.3){$o$};
		\draw[dashed] (-6,0)[out=30,in=145] to (5,0);
		\draw (-6,0)[out=80,in=-145] to (-3,6)
		[out=35,in=130] to (4.4,6)
		[out=-55,in=80] to (5,0)
		[out=-105,in=10] to (0.5,-4)
		[out=-170,in=-10] to (-3,-4)
		[out=170,in=-100] to (-6,0);
		\node (start) at (1.9,5.7)[draw,ellipse,minimum height=0.5cm,minimum width=1.57cm,rotate=-10] {};
		\foreach \angle in {155,143,125,100,70,45,25,13}  	
		\draw (node cs:name=start,angle=\angle+15) .. controls +(-20:0.6cm) and +(0,0) .. (node cs:name=start,angle=-\angle);
		\draw[dashed](0,0)--(0,6);
		\draw[dashed](0,0)--(2,6) node at (2.2,5.7) {$x$} node at (2.8,6.6) {$\sigma$};
		\draw[dashed](0,0)--(3.8,5.3);
		\draw[dashed,-latex](0,0)--(0,3);
		\draw[dashed,-latex](0,0)--(1,3) node at (0.6,3) {$u$};
		\draw[dashed,-latex](0,0)--(1.9,2.65) node at (3.6,3) {${\st\triangle}_K(\sigma)$};
		\draw[-latex](0,6)--(-0.8,8.8);
		\draw[-latex](2,6)--(2.2,8.8) node at (2.7,8.8) {$v$};
		\draw[-latex](3.8,5.3)--(5,8);		
		\end{tikzpicture}
}
\end{center}

\vs{5}

\subsubsection{Normal Cone and Gauss Image}

If $K\in\kon$,
then the {\it normal cone} of $K$ at $x \in \partial K$ is the set of all outer
unit normals of supporting hyperplanes that contain $x$, i.e.,
\[
N(K,x)=\{v\in\sn: \text{$(y-x)\cdot v \le 0$ for all $y\in K$}\}=\{v\in\sn : x\in H_K(v)\}.
\]
The normal cone at $x$ contains only one unit vector $v$ if and only if $\rho$ is differentiable at $x$.

For $\sigma\subset\partial K$, the {\it Gauss image of $\sigma$}, also called the
{\it spherical image of $\sigma$}, is the union of the normal cones at each $x \in \sigma$, i.e.,
\begin{equation*}
\bnu\lsub K(\sigma) = \{v\in\sn : x\in H_K(v)\ \text{for some}\ x\in\sigma \}= \bigcup_{x\in \sigma} N(K,x) \subset\sn.
\end{equation*}
If $\partial K$ is differentiable at $x \in \partial K$, then there is a unique supporting hyperplane,
which implies that the normal cone has only a single vector and therefore
\[
  \bnu\lsub K(\{x\})= N(K,x) = \{\nuk(x)\}.
\]

For a set $\eta\subset\sn$ of unit normal vectors, the {\it reverse Gauss image of $\eta$} is defined by
\begin{equation}\label{reverseGaussImage}
\bnu_K^*(\eta)=\{ x\in \partial K : x\in H_K(v) \text{ for some } v\in \eta\} \subset \partial K.
\end{equation}

When $K$ contains the origin in its interior, it is more convenient to use
the {\it radial Gauss image} of $\omega \subset \sn$, defined to be
\begin{equation}\label{radial-gauss-image}
\balpha\lsub K(\omega) = \bnu\lsub K(\rk(\omega)),
\end{equation}
or equivalently,
\[
\balpha\lsub K(\omega)=\big\{v\in\sn : \rk(u)\in H_K(v) \text{ for some } u\in\omega \big\}
=\bigcup_{u\in\omega} N(K,r_K(u)).
\]
The radial Gauss image maps subsets of the radial sphere to subsets of the normal sphere.

The {\it reverse radial Gauss image} $\balpha_K^*$ is defined by
\begin{equation}\label{reverse-radial-gauss-image}
\balpha_K^*(\eta)   =\big\{u\in\sn : \rk(u)\in H_K(v) \text{ for some } v\in\eta \big\}, \quad \eta\subset\sn,
\end{equation}
that is,
\[
\balpha_K^*(\eta) = r_K^{-1}(\bnu_K^*(\eta)).
\]

The reverse radial Gauss image maps subsets of the normal sphere to subsets of the radial sphere.
 For an arbitrary convex body $K$
(that contains the origin in its interior), the reverse radial Gauss image of $K$
and the radial Gauss image of $K^*$ satisfy (see \cite{HLYZ16acta})
\[
\balpha_{K^*}(\eta) = \balpha_K^*(\eta), \quad \eta \subset \sn.
\]

\section{Descriptions of Minkowski Problems and Variational Methods}\

\subsubsection{Geometric Measures of Convex Bodies}
The essential view is that an interesting  geometric
measure should be the differential of an important geometric invariant.
Well-known geometric measures arise in such a way.
For example,
surface area measure is the differential of volume,
area measures are differentials of quermassintegrals,
and dual curvature measures are differentials
of dual quermassintegrals. We give a description of geometric measures.

Let $\mathcal K$ be a class of convex bodies in $\rn$ and
$\mathcal B$ the set of finite Borel measures on $\sn$.
A geometric measure of convex bodies over $\mathcal K$ is a measure-valued map,
\[
 M : \mathcal K \to \mathcal B, \quad
K \mapsto  M(K,\cdot),
\]
which is the differential of a global geometric invariant $W$.
That is, for each convex body $K$,
\begin{equation}\label{vf}
\frac{dW(K_t)}{dt}\Big|_{t=0^+}
= \int_{\sn} f \, d{ M(K,\cdot)},
\end{equation}
where $K_t$, $0\le t<\delta$, is a perturbation family of convex bodies generated by $K$, a geometric
operation, such as the Wulff shape or convex hull, and some continuous
function $f$ on $\sn$, such as a support function, a radial function, or any continuous function. 
The class $\mathcal K$ is assumed to contain every such perturbation family $K_t$.

The essential difficulty in proving such a variational formula
is that convex bodies are not smooth
and the Wulff shape and convex hull are also not smooth as operations.

\subsubsection{General Minkowski Problem}
{\it Given a finite Borel measure $\mu$ on $\sn$,
find necessary and sufficient conditions on the measure to guarantee the existence of a
convex body $K$ in $\rn$ that is a solution to the geometric measure equation (called the Minkowski equation)
\begin{equation}\label{mequ}
 M(K,\cdot)=\mu.
\end{equation}
}

In the smooth case and when $\mu$ has a density $f$, the measure
equation \eqref{mequ} often becomes a Monge-Amp\`ere type equation \eqref{monge}.
The associated uniqueness problem reads as follows.

\subsubsection{Uniqueness Problem of Geometric Measure}
{\it For two convex bodies $K$ and $L$, determine how $K$ and $L$ are related if
$ M(K,\cdot) =  M(L,\cdot)$.
}

\vspace{2pt}

The main method to solve the measure equation \eqref{mequ} is a variational method.
For different measure equations, different variational methods are needed which
have similar steps but require different techniques and details to carry them out.
We give a rough description of the steps.

\subsubsection{Variational Method}

\dr{1} Construct an optimization problem of some geometric invariant that is related to the geometric measure.
Note that the space of
convex bodies is not an appropriate setting for a variational problem, since arbitrary
deformations of a convex body are not necessarily convex. Instead, it is necessary to use
a function space and the Wulff shape $\wul{h}$ or convex hull $\bla h \bra$
of a positive continuous function $h$. In view of the formula \eqref{vf}, a possible example is
the maximization problem
\begin{equation}\label{op}
\sup\Big\{W(\wul{h}) :  \int_{\sn} \varphi(h)\, d\mu=1, \ h\in C^+(\sn)\Big\},
\end{equation}
where $\varphi : (0, \infty) \to \R$ is some continuous function, for instance,
$\varphi(t) = t^p$ or $\varphi(t) = \log t$.
The maximum may be achieved when $h$ is the support function of a convex body.

\dr{2} Prove the existence of solutions to the optimization problem. One needs to show that
an optimizing sequence is convergent and its limit gives a non-degenerate convex body.
To do this, one shows that
it suffices to consider only optimizing sequences of support functions of convex bodies
and to show that an optimizing sequence is uniformly bounded. Then the Blaschke selection theorem
implies that the optimizing sequence converges to the support function of a compact convex set.
Then one needs to show that the compact convex set has nonempty interior, and also often needs to
show that the origin is inside the interior.

\dr{3} Prove that a solution of the optimization problem gives a solution of the Minkowski problem,
that is, the Euler-Lagrange equation of the optimization problem \eqref{op} is the measure
equation \eqref{mequ}. The key is to show that the geometric invariant $W$ is differentiable
and its differential
is the geometric measure $M$. That is, one needs to show that the variational formula \eqref{vf}
also holds true for $t= 0^-$.

\section{Minkowski Problems in the Brunn-Minkowski Theory}\

In this section, we describe Minkowski problems of geometric measures associated with mixed volumes.

\subsection{Basics of the Brunn-Minkowski Theory}

\subsubsection{Mixed Volumes}
The Brunn-Minkowski theory studies mixed volumes that arise from the interaction of
volume, the Minkowski addition, and differentiation.
If $K, L$ are convex bodies in $\rn$, the Minkowski sum $K+L$ and dilation $tK$, $t>0$, are
convex bodies defined by
\begin{equation}\label{Madd}
K+L=\{x+y : x\in K, \ y\in L\}, \ \ \ tK = \{ tx : x \in K \}.
\end{equation}
There are identities,
\[
  h_{K+L} = h_K + h_L, \ \ \ h_{tK} = th_K.
\]

The first mixed volume $V_1(K,L)$ of two convex bodies $K$ and $L$ is, up to a constant factor,
the one-sided derivative of the volume with respect to the Minkowski addition,
\[
  V_1(K,L) = \frac{1}{n}\lim_{t\rightarrow 0^+} \frac{V(K+tL)-V(K)}{t}.
\]
Higher order mixed volumes can be defined as follows. It can be shown (see, for example,
Schneider \cite[Theorem 5.1.7]{S14}) that $V(K+tL)$ is a polynomial with respect to $t$.
The mixed volumes $V_i(K,L)$, $0 \le i \le n$, of $K$ and $L$ are defined in terms of
the coefficients of the polynomial,
\[
V(K + tL) = \sum_{i=0}^n \binom ni V_i(K, L) t^{i}, \ \ \ t> 0.
\]
More generally, one can define the mixed volume of $m$ convex bodies, $K_1, \ldots, K_m$,
to be the coefficients of the homogeneous polynomial,
\[
  V(t_1K_1 + \cdots + t_mK_m) = \sum_{i_1, \cdots, i_n = 1}^m
 V(K_{i_1}, \dots, K_{i_n})\, t_{i_1}\cdots t_{i_n}, \ \ \  t_1, \ldots, t_m > 0,
\]
where the coefficients
$V(K_{i_1},\cdots,K_{i_n})$ are nonnegative, symmetric in the indices, and dependent only on the bodies
$K_{i_1},\cdots,K_{i_n}$. The invariant $V(K_{i_1},\cdots,K_{i_n})$ is called the {\it mixed volume} of convex bodies
$K_{i_1},\cdots,K_{i_n}$,  with the property
$V(K_{i_1},\cdots,K_{i_n})= V(K_{i_1})$ when $i_1=\cdots=i_n$.

The $i$th mixed volume of $K$ and $L$ is given by
\[
  V_i(K,L) = V(K, \dots, K, L, \dots, L),
\]
where $K$ appears $n-i$ times, and  $L$ appears $i$ times.

Mixed volumes are multi-linear with respect to
the Minkowski addition. By this and the Riesz representation theorem,
it is easily shown that for convex bodies
$K_{2}, \ldots, K_{n}$ in $\rn$ there exists a unique
finite Borel measure $S(K_{2}, \ldots, K_{n}, \cdot)$ on $\sn$ such that
\[
V(K_{1}, K_2, \cdots,K_{n})= \frac1n \int_{\sn} h_{K_{1}}(v)\, dS(K_{2}, \ldots, K_{n}, v),
\]
for any convex body $K_1$ in $\rn$.
The measure $S(K_{2}, \ldots, K_{n}, \cdot)$ is called the {\it mixed area measure} of
$K_{2}, \ldots, K_{n}$. It is translation invariant. The centroid of a mixed area measure is at the origin,
\[
\int_{\sn} v\, dS(K_2,\ldots, K_2, v)=0.
\]

Mixed volumes satisfy the Aleksandrov-Fenchel inequality,
\[
V(K_1, K_2, K_3, \ldots, K_n)^2 \ge V(K_1, K_1,  K_3, \ldots, K_n)V(K_2, K_2,  K_3, \ldots, K_n).
\]
One corollary is the following:
\[
V(K_1, K_2, \ldots, K_n)^n \ge V(K_1) V(K_2) \cdots V(K_n),
\]
which implies Minkowski's first inequality for mixed volumes,
\begin{equation}\label{b1.7}
V_1(K, L)^n \ge V(K)^{n-1} V(L).
\end{equation}
Equality holds if and only if $K$ and $L$ are homothetic.
When $L$ is the unit ball, \eqref{b1.7} becomes the classical isoperimetric inequality.

The Minkowski inequality is equivalent to the Brunn-Minkowski inequality,
which states that volume is a log-concave function with respect to Minkowski sum. Specifically,
given convex bodies $K$ and $L$ in $\rn$,
\begin{equation}\label{BMI}
V((1-t)K+tL) \ge V(K)^{1-t} V(L)^t, \ \ \ t>0,
\end{equation}
with equality if and only if $K$ and $L$ are homothetic. An equivalent form of the inequality is
\begin{equation}\label{BMI2}
V(K+L)^{\frac{1}{n}} \ge V(K)^{\frac{1}{n}} + V(L)^{\frac{1}{n}}.
\end{equation}
There are by now many proofs of the Brunn-Minkowski inequality. See, for example,
the book by Schneider \cite{S14} and the survey by Gardner \cite{G02bull}.

\subsubsection{Quermassintegrals}
The most important mixed volumes are quermassintegrals, which include
volume, surface area and mean width.

The first mixed volume $V_1(K,L)$ can be viewed as the surface area of $K$
with respect to the (possibly asymmetric) norm on $\rn$ whose unit ball is $L$.
In particular, if $B$ is the standard Euclidean unit ball,
$$
nV_1(K,B) = \Hn(\partial K) = S(K)=\text{ Euclidean surface area of }\partial K.
$$

More generally, the {\it $i$th quermassintegral of $K$} is defined to be
\[
  W_i(K) = V_i(K,B).
\]
It turns out that
\[
  W_0(K)=V(K),\ W_1(K)= \frac1n S(K), \ \text{ and }W_n(K)=V(B)=\omega_n.
\]

The surface area of a polytope is the sum of areas of the $(n-1)$-dimensional faces.
More generally, the quermassintegral $W_{n-i}(K)$ is a weighted
sum of the areas of the $i$-dimensional faces of a polytope.

When $K$ has a smooth boundary, quermassintegrals are
the integrals of intermediate curvatures of the hypersurface $\partial K$,
\[
W_{i+1}(K)=\frac1n \int_{\partial K} H_i(x)\, dx, \quad i=0, 1, \ldots, n-1,
\]
where $H_0=1$ and $H_i$ is the $i$th intermediate curvature of $\partial K$, that is, the $i$th
normalized elementary symmetric function of the principal curvatures of $\partial K$.

For a general convex body $K$,
the quermassintegrals can be defined as integrals of the areas of projections of $K$ onto subspaces,
\[
W_{n-i}(K) = \frac{\omega_n}{\omega_i} \int_{G_{n,i}} V_i(K|\xi)\, d\xi,
\]
where $K|\xi$ denotes the orthogonal projection
of $K$ onto the subspace $\xi \in G_{n,i}$.

\subsubsection{Area Measures}

Geometric measures associated with mixed volumes are mixed area measures. The most important case is about
two convex bodies. For convex bodies $K$ and $L$ in $\rn$, the {\it $i$th mixed area measure} $S_i(K,L, \cdot)$ is
defined by
\[
S_i(K,L, \cdot) = S(K, \ldots, K, L, \ldots, L, \cdot), \ \ \ i=0,\ldots, n-1,
\]
where $K$ appears $i$ times and $L$ appears $n-i-1$ times. When $L$ is the unit ball, the area measure
\[
S_i(K, \cdot) = S_i(K,B, \cdot), \ \ \ i=0,\ldots, n-1,
\]
is called the {\it $i$th area measure of $K$}. The $(n-1)$th area measure $S_{n-1}(K,\cdot)$ is called
the surface area measure of $K$, usually denoted by $S(K,\cdot)$ or $S_K$.

The total mass of the $i$th area measure is
the $(n-i)$th quermassintegral with factor $n$,
\begin{equation}\label{b3.2}
S_i(K,S^{n-1}) = nW_{n-i}(K), \quad i=0,\ldots, n-1.
\end{equation}

When the boundary $\partial K$ is $C^2$ and has positive curvature, then
\[
dS_i(K, v) = \sigma_i\, dv,  \quad i=0,\ldots, n-1,
\]
where
\[
  \sigma_i=\frac{H_{n-i-1}}{H_{n-1}}
\]
is called the $i$th curvature function and is equal to the $i$th normalized elementary
symmetric function of the principal radii (i.e., reciprocals of the principal curvatures).
When $i=0$, $S_0(K, \cdot)$ is spherical Lebesgue
measure and $\sigma_0=1$. When $i=n-1$, $S_{n-1}(K,\cdot)=S(K,\cdot)$ is the surface area measure of $K$ and
$\sigma_{n-1}= f_K$ is the curvature function of $K$.

Mixed area measures are one-sided derivatives of mixed volumes. Given convex bodies $K, Q, L$,
\begin{equation}\label{b3.3}
\frac{d}{dt}V_{n-i-1}(K+tQ, L)\Big|_{t=0^+} =\frac {i+1}n \int_{\sn} h_Q(v)\, dS_{i}(K,L,v), \quad i=0,\ldots, n-1.
\end{equation}

\subsubsection{Uniqueness Problem of Area Measures}
A convex body is uniquely determined, up to a translation, by
its $i$th area measure. More generally, for convex bodies $K$, $Q$, and a smooth convex body $L$, $i=1,\ldots, n-1$,
\[
S_i(K, L,\cdot) = S_i(Q, L, \cdot) \ \ \ \Longrightarrow \ \ \ K = Q + x_0, \ \ \ \text{for some } x_0 \in \rn.
\]
When $L$ is not smooth, the uniqueness problem of mixed area measures is still unsolved.
It follows from
mixed volume inequalities and their equality conditions (see Schneider \cite{S14}, p. 448).

\subsubsection{Surface Area Measure}

The surface area measure  $S(K, \cdot)$ is a fundamental geometric measure of convex bodies.
It generalizes and interpolates between the extreme cases of a smooth positively curved convex body
and a convex polytope. It has a more direct geometric definition. The idea is
to use the reverse Gauss image \eqref{reverseGaussImage} and define the surface area measure
to be the pullback of $\Hn(\partial K,\cdot)$
by the reverse Gauss image,
\begin{equation}\label{b1.1}
S(K, \eta) = \Hn\big(\bnu_K^*(\eta)\big), \ \ \  \eta\subset \sn,
\end{equation}
where $\bnu_K^*(\eta)$ is the reverse Gauss image of a Borel set $\eta$ defined in \eqref{reverseGaussImage}. That is,
$S(K,\eta)$ is equal to the surface area of the subset $\bnu_K^*(\eta) \subset \partial K$ whose unit normals
belong to $\eta$.
Equivalently, it can be defined by
\begin{equation}\label{b1.2}
\int_{\sn} f(v)\, dS(K, v) = \int_{\partial'K} f(\nuk(x))\, dx,
\end{equation}
for any bounded Borel function $f$ over $\sn$. The total mass of surface area measure is the surface area of $K$,
i.e., $S(K,\sn) = S(K)$.

By the Lebesgue-Radon-Nikodym decomposition theorem, the surface area measure
can be decomposed into the sum of two mutually singular measures $\lambda_K$ and $\zeta_K$,
\begin{equation}\label{b1.4}
S(K,\cdot) = \lambda_K + \zeta_K, \ \ \ \ d\lambda_K = f_K(v)\, dv,
\end{equation}
where $\lambda_K$ is absolutely continuous with respect to the spherical Lebesgue measure
with density $f_K$, the curvature function of $K$, see \eqref{0.18}.

 When $K$ is a polytope with outer unit normals $v_1, \ldots, v_m$ to its facets and corresponding facet areas
$a_1, \ldots, a_m$, then $S(K,\cdot)$ is the discrete measure,
\begin{equation}\label{b1.3}
S(K,\cdot) =\zeta_K= \sum_{i=1}^m a_i \, \delta_{v_i},
\end{equation}
where $\delta_{v_i}$ is the delta measure concentrated at $v_i$.

\subsection{Classical Minkowski Problem for Surface Area Measure} \label{classical}

In this section we present in detail the classical Minkowski problem first introduced by Minkowski.
It asks what are the necessary and sufficient conditions for a Borel measure on the unit sphere
to be the surface area measure of a convex body. The classical Minkowski problem is the fundamental
one in the Brunn-Minkowski theory, as formulated by Fenchel and Jessen \cite{FJ} and by Aleksandrov \cite{Al3}:

\vs{3}

\fr{The Minkowski Problem} {\it Given a finite Borel measure $\mu$ on $\sn$, find the necessary and sufficient conditions
such that there exists
a convex body $K$ in $\rn$ whose surface area measure is $\mu$, that is, the measure equation
$S(K,\cdot) = \mu$ has a solution $K$.}

\vs{3}

When the surface area measure of a convex body is discrete, then the body is a polytope
and the measure is given by \eqref{b1.3}. The Minkowski problem can therefore be viewed as
that of finding a convex polytope with prescribed facet areas and corresponding outer unit normals.
The discrete Minkowski problem in $\R^2$ is described in \S\ref{mp-polygon}. In $\rn$, it states:
\vs{3}

\fr{Discrete Minkowski Problem}
{\it Given positive numbers $f_1, \ldots, f_m$, and unit vectors $v_1, \ldots, v_m$ in $\rn$,
find the necessary and sufficient conditions such that there exists a convex polytope $K$ in $\rn$
whose facets have areas $f_1, \ldots, f_m$ and corresponding outer unit normals $v_1, \ldots, v_m$.}

\vs{3}

When $\mu$ is absolutely continuous with respect to spherical Lebesgue measure, then,
by \eqref{b1.4}, the Minkowski problem is equivalent to that of finding a convex body
with prescribed curvature function \eqref{0.18}.

Since the surface area measure of a convex body is invariant under translations,
there are infinitely many solutions to the Minkowski problem for a given measure.
Moreover, as we have already seen in \S\ref{mp-polygon}, the translation invariance
also implies that there is an obstruction to a measure being the surface area measure
of a convex body. Below, we show that the Minkowski problem for a given finite measure has a solution
if and only if the measure is not concentrated on a closed hemisphere and
has its centroid at the origin, that is, it satisfies \eqref{balanced}.

\subsubsection{Solution to the Minkowski Problem}\label{4.1}

Minkowski solved this problem and also proved the uniqueness of the solutions for
the discrete and continuous cases \cite{M1897, M1903}.
Aleksandrov \cite{Al2, Al3} and Fenchel \& Jessen \cite{FJ} independently solved the problem
for arbitrary convex bodies by introducing surface area measure and using a variational method
that transforms the Minkowski problem into an optimization problem. They proved the following:

\begin{MET} If $\mu$ is a finite Borel measure on $\sn$,
 then the measure equation $S(K,\cdot)=\mu$ has a solution if and only if
$\mu$ is not concentrated on a closed hemisphere and
\begin{equation}\label{balanced}
\int_{\sn} v\, d\mu(v)=0.
\end{equation}
\end{MET}

\vs{3}

\subsubsection{Outline of Proof for the Minkowski Existence Theorem}

We give a sketch of Aleksandrov's variational approach.  The main steps
are as described in the previous section.

The optimization problem associated with the classical Minkowski problem is the following
maximization problem of volume,
\begin{equation}\label{b2.1}
\sup\Big\{V(\wul{h}) :  \int_{\sn} h\, d\mu=1, \ h\in C^+(\sn)\Big\}.
\end{equation}
This maximization problem is motivated
by the mixed volume inequality \eqref{b1.7}. If $\mu$ is the surface area of $K$,
the maximum is achieved only when $h$ is, up to a positive scalar factor, the support function of $K$.

To prove the existence of solutions to the maximization problem, one needs to show
that an optimizing sequence is uniformly bounded. The uniform bound
comes from the equation $\int_{\sn} h\, d\mu=1$ quite easily.
Then Blaschke's selection theorem
implies that the optimizing sequence converges to the support function of a compact convex set.
The non-degeneracy of the limiting compact convex set
follows from the maximality of its volume.

To show the Euler-Lagrange equation of the maximization problem is the equation of the Minkowski problem,
one needs to show that $V(\wul{h})$ is differentiable with respect to $h \in C^+(\sn)$.
This is a key nontrivial step in Aleksandrov's proof. Given $h \in C^+(\sn)$, it suffices
to consider variations of the form
\[
h_t(v) =h(v) + t f(v),
\]
where $f \in C(\sn)$ and $t\in \R$ is small. Aleksandrov's formula is
\begin{equation}\label{b2.2}
\frac{dV(\wul{h_t})}{dt}\Big|_{t=0} = \int_{\sn} f(v)\, dS(\wul{h},v), \ \ \ f\in C(\sn).
\end{equation}
The proof of this formula uses properties of the surface area measure of the Wulff shape and
Minkowski's mixed volume inequality \eqref{b1.7} (see, for example, Lemma 7.5.3 in Schneider \cite{S14}, or
Lemma 1 in Haberl-LYZ \cite{HLYZ10}).

\vs{5}

\subsubsection{Continuous Minkowski Problem}

Recall that if the surface area measure $S_K$ is absolutely continuous with respect to
the spherical Lebesgue measure,
\[
  dS_K = f_K\,dv,
\]
then the density $f_K$ is the curvature function of $K$. If the boundary of $K$ is $C^2$,
then $f_K$ is continuous and satisfies
\[
  f_K(\nu_K(x)) = \frac{1}{\lkappa_K(x)},
\]
where $\nu_K: \partial K \rightarrow \sn$ is the Gauss map and $\lkappa_K: \partial K \rightarrow (0,\infty)$
is the Gauss curvature.

If we assume that
\[
  d\mu = f\,dv,
\]
where $f \in C^+(\sn)$, then we can ask if there is a solution to the Minkowski problem
that is a convex body $K$ with $C^2$ boundary. If so, then $f$ is the curvature function of $K$, i.e.,
\[
  f_K(\nu_K(x)) = f(\nu_K(x)), \ \text{ for all } x \in \partial K.
\]

Since the curvature function is the reciprocal of Gauss curvature, this is equivalent to the following:

\vspace{3pt}

\fr{Continuous Minkowski Problem} {\it
Given a positive continuous function $g$ on $\sn$, find the necessary and sufficient conditions
such that there exists a closed convex hypersurface $\Sigma$ in $\rn$ whose Gauss curvature
$\lkappa$ and Gauss map $\nu$ satisfy
$\lkappa(\nu^{-1}) = g.$
}

\vspace{3pt}

On one hand, any solution $\Sigma$ to the continuous Minkowski problem gives a solution to the Minkowski problem
\begin{equation}\label{cMP}
  dS_K = \frac1{g(v)} \,dv,
\end{equation}
where $K$ is the convex body with $\partial K =\Sigma$.
On the other hand, the Minkowski existence theorem shows that there is a solution $K$ to
the Minkowski problem \eqref{cMP} if and only if
\[\int_{\sn} \frac{v}{g(v)} \, dv = 0.\]

Therefore, if we can show that the body $K$ has $C^2$ boundary, then it is a solution to the
continuous Minkowski problem. This is equivalent to showing that the support function $h_K$ is $C^2$
because the positivity of the curvature function $f_K=1/g$.

A generalized continuous Minkowski problem in the differential geometry is to find an embedding
 $G: \sn \to \Sigma \subset \rn$ such that  $\lkappa(G) = g$ for a given positive continuous function $g$ on $\sn$.
A solution was obtained by Gluck \cite{Gl1, Gl2}.

\subsubsection{Continuous Minkowski Problem as a PDE}

The continuous Minkowski problem is equivalent to the Monge-Amp\`ere equation
\begin{equation}\label{Mp1}
\det(\nabla_{ij} h + h\delta_{ij} ) = f,
\end{equation}
where $f$ is a positive continuous function on $\sn$.  The solution $h$ satisfies that
\begin{equation}\label{positive-definite}
  (\nabla_{ij}h + h\delta_{ij}) \ \ \text{is positive definite}
\end{equation}
if and only if it is the support function of a convex body $K$ whose curvature function is $f$.
This is also equivalent to \eqref{Mp1} being an elliptic partial differential equation.

Therefore, if $h$ is a solution to \eqref{Mp1} that satisfies \eqref{positive-definite},
then it is the support function of a solution to the continuous Minkowski problem. In particular,
the continuous Minkowski problem becomes a regularity question for solutions to elliptic
Monge-Amp\`ere equations on a sphere. This is a difficult question studied by many scholars.
Major contributions were made by Lewy \cite{Lew38}, Nirenberg \cite{N53},
Cheng-Yau \cite{ChengYau}, Pogorolov \cite{P78}, Caffarelli \cite{Caf90Ann},
culminating in the following result:

\vs{3}

\fr{Regularity of the Minkowski Problem} {\it If $g \in C^\alpha(\sn)$, $g>0$, and
\[
\int_{\sn} \frac v{g(v)}\, dv =0,
\]
then there exists a convex body $K$ in $\rn$ with $C^{2,\alpha}$ boundary
such that its Gauss curvature  $\lkappa$
as a function of its outer unit normal is equal to $g$.
}

\vs{3}

It is of interest to generalize the continuous
Minkowski problem of prescribing Gauss curvature to intermediate curvatures, that is, replacing
the Gauss curvature $\lkappa$ by the $i$th
intermediate curvature $H_i$ of $\partial K$.
Guan-Guan \cite{GG02} solved the symmetric case by proving:

\vs{5}

{\it Let $g \in C^{1,1}(\sn)$ be positive and even and $i\in \{1, \ldots, n-2\}$.
Then there exists an origin-symmetric convex body $K$ in $\rn$ with $C^{3,\alpha}$
boundary and positive curvature
so that its $i$th intermediate curvature $H_i$
as a function of its outer unit normal is equal to $g$.
}
\vs{5}

They also obtained results for the nonsymmetric case.
Aleksandrov \cite{A56} and Chern \cite{Chern1, Chern2} posed the uniqueness problem
of prescribing intermediate curvatures.

\subsection{Minkowski Problems of Area Measures}

Which finite Borel measure on the unit sphere can be a mixed area measure of a convex body?
This is the problem of prescribing mixed area measures, which states:

\subsubsection{Minkowski Problem of Mixed Area Measures}
{\it Given a finite Borel measure $\mu$ on $\sn$, a convex body $L$ in $\rn$, and $i\in \{1, \ldots, n-1\}$,
find the necessary and sufficient conditions such that there exists
a convex body $K$ in $\rn$ satisfying $S_i(K,L,\cdot) = \mu$.}

\vs{5}

The case $i=n-1$ is the Minkowski problem. The general case has not been studied except
when $L=B^n$ which is the following:

\subsubsection{Christoffel-Minkowski Problem} {\it Given a finite Borel measure
$\mu$ on $\sn$ and $i\in \{1, \ldots, n-1\}$,
find the necessary and sufficient conditions such that there exists
a convex body $K$ in $\rn$ satisfying $S_i(K, \cdot) = \mu$.}

\vs{5}

When $\mu$ has a density $f$, this problem is equivalent
to the partial differential equation
\[
\sigma_i\big(\nabla_{kj}h + h\delta_{kj}\big) =f,
\]
where $\sigma_i$ denotes the normalized $i$th elementary symmetric function of eigenvalues
of the matrix, that is, the sum of $i$th principal minors of the matrix divided by $\binom {n-1}i$.

\vs{4}

The case $i=1$ is the Christoffel problem with corresponding partial differential equation
\[
\Delta_{\sn} h + (n-1)h = f,
\]
where $\Delta_{\sn}$ is the spherical Laplacian.
The smooth case in $\R^3$ was first treated by Christoffel \cite{C1865}
as a problem of prescribing the sum of radii of curvature.
The general case in all dimensions was solved independently by Firey \cite{F68, F70}
and Berg \cite{B69}. Schneider \cite{S77} gave more explicit necessary and sufficient
conditions for the existence of solutions in the polytope case.
The case of even measures and symmetric convex bodies
was studied in connection with the spherical Radon transform
 by Goodey-Weil \cite{GW92}.
A short approach to the Christoffel problem was found by Grinberg and Zhang \cite{GZ99}.
A Fourier transform approach was given by Goodey-Yaskin-Yaskina \cite{GYY11}.
When $f$ is a positive Lipschitz function, Li-Wan-Wang \cite{LWW20} gave simpler
necessary and sufficient conditions for the existence of solutions.

For $1<i<n-1$, the Christoffel-Minkowski problem has remained unsolved for measures.
For sufficiently smooth
convex bodies of revolution, Firey \cite{F70.1} gave necessary and sufficient conditions for the existence
of solution. For centrally symmetric convex bodies with $C^{2,\alpha}$ boundary and positive curvature,
Zhang \cite{Z94} showed that the Christoffel-Minkowski problem is solvable in a neighborhood of $f$
when it is solvable for a smooth $f$. A breakthrough
was achieved by Guan-Ma \cite{GM03} who found sufficient conditions for both existence and
regularity of solutions to the Christoffel-Minkowski problem when $f$ is assumed to be sufficiently smooth:

{\it Let $i\in \{1, \ldots, n-2\}$ and $\mu$ a finite Borel measure on $\sn$ with positive density $f\in C^{1,1}(\sn)$.
If the centroid of $\mu$ is at the origin, i.e.,
\[
\int_{\sn} v f(v) \, dv=0,
\]
and $f^{-\frac1i}$ is convex when it is extended to $\rn\setminus\{0\}$ homogeneously of degree $1$,
then there exists a convex body $K$ of class $C^{3,\alpha}$ and positive curvature so that $S_i(K, \cdot) = \mu$.
}

\vs{4}

The proof in \cite{GM03} required an additional technical assumption on the density $f$.
This was shown to be unnecessary by Sheng-Trudinger-Wang \cite{STW} and also
by Andrews and Ma in an unpublished note. Bryan-Ivaki-Scheuer\cite{BIS20} recently used
a fully nonlinear curvature flow to derive the existence result under those conditions of $f$.

The continuity method was used to prove the results above. A natural question is whether
Christoffel-Minkowski problems could be solved using a variational method
similar to that used by Minkowski, Aleksandrov, and Fenchel for
 the classical Minkowski problem. The main obstacle
for doing this is the lack of a differential formula for quermassintegrals like
the one for volume \eqref{b2.2}.
The one-sided differential formula \eqref{b3.3} is not sufficient for transforming the
Christoffel-Minkowski problem into an optimization problem.

\subsection{Minkowski Problems of Curvature Measures}

\subsubsection{Curvature Measures of Convex Bodies}

Given any $x \in \rn$, let $d(K,x)$ be the distance from $x$ to $K\in \kon$, and, given $t > 0$,
denote the annulus with thickness $t$ by
\[
  K_t = \{ x\notin K\ : d(K,x) \le t \} = (K+tB)\setminus K.
\]
Consider the metric projection map
\[
  p\lsub{K} : \rn \setminus K \to \partial K,
\]
where $p\lsub{K}(x) \in \partial K$ is the unique point such that
\[
|p\lsub{K}(x)-x| = d(K,x).
\]
Given a Borel set $\omega\subset \sn$, let
\begin{equation}\label{cmp}
A_t(K, \omega) = K_t \cap p_K^{-1}(r_K(\omega)).
\end{equation}
There is the Steiner-type formula,
\begin{equation}\label{cme}
V(A_t(K,\omega)) = \frac1n \sum_{i=0}^{n-1} \binom ni t^{n-i} C_i(K,\omega),
\end{equation}
where $C_i(K,\cdot)$ is a Borel measure on $\sn$, called the {\it $i$th curvature measure} of $K$.
\vs{5}

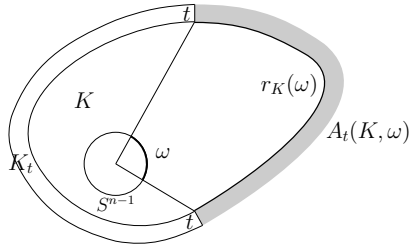
\begin{center}
\scalebox{0.6}{
\begin{tikzpicture} [scale=0.35]
	\draw (0,0) circle [radius=2cm] node at (0,-2.5) {$\sn$};
	
	\draw (-5,0) to [out=-60,in=-150] (5,-3);
    \draw [line width =1.6pt] (5,-3) to [out=30,in=-30](12,7) to [out=150,in=0](5,9);
    \draw (5,9) to [out=180,in=120](-5,0)
     node at (-6,0){\large $K_t$}
     node at (-2,4){\large $K$}
	node at (11,5) {\large $r_K(\omega)$}
	node at (16,2) {\large $A_t(K,\omega)$};
	
	\fill[fill=white!80!black]
	(5,10.1)to [out=0,in=145](13.2,7.6) to [out=-40,in=30] (5.5,-3.9)to (5,-3) to [out=30,in=-30](12,7) to [out=150,in=0](5,9) to (5,10.1);
	\draw [line width=0.6pt](5.5,-3.9) to (5,-3);
    \draw [line width=0.6pt](5,10.1) to (5,9)
    node at (4.7,-3.7){\large$t$}node at (4.5,9.5){\large$t$};
	\draw (0,0) to (5,9);
	\draw (0,0) to (5,-3);

	\draw [line width=1.2pt] (1.7,-1) arc (-30:60:2) node at (3,0.7){\large$\omega$};
	\draw (5,10.1) to  [out=183,in=20] (0,9.1) to [out=-157,in=70](-6.6,3)to [out=-105,in=120] (-6,-0.5) to [out=-60,in=165] (0,-4.9) to [out=-10,in=-155] (5.5,-3.9);
\end{tikzpicture}
}
\end{center}

Federer \cite{F59} defined curvature measures as measures in $\rn$
for sets of positive reach that are the sets such that the metric
projection map can be defined. Convex bodies are sets of positive reach.
Schneider \cite{S78} gave a unified treatment of curvature measures and
area measures for convex bodies. As shown in \cite{S78}, curvature measures
of a convex body concentrate on its boundary.
The curvature measures defined above for a convex body $K\in \kon$
are pushforwards of Federer's curvature measures to the unit sphere $\sn$
by the inverse radial map $r_K^{-1}$. Since the radial map $\rk$ is
bi-Lipschitz, one can define the curvature measures equivalently on
either $\partial K$ or $\sn$, although they are different measures.

The 0th curvature measure $C_0(K,\cdot)$ was defined
and studied earlier by Aleksandrov, who called it the {\it integral curvature of $K$}.
It can be defined, for Borel $\omega \subset \sn$,  by
\begin{equation}\label{4.2}
C_0(K, \omega) = \mathcal H^{n-1}(\balpha_K(\omega)),
\end{equation}
where $\balpha_K$ is the radial Gauss image map defined by \eqref{radial-gauss-image}.
In other words, $C_0(K,\omega)$ is the spherical Lebesgue measure of $\balpha_K(\omega)$,
which has total measure
$C_0(K,\sn)$ equal to the surface area of $\sn$. Its geometric
meaning is more apparent when $K$ is a polytope. In this case, $C_0(K,\cdot)$ is a discrete measure
concentrated on the radial directions of the vertices, where the mass of the measure at a vertex
is the solid angle formed by the unit normals.

On the other hand, the $(n - 1)$th curvature measure $C_{n-1}(K, \cdot)$ on $\sn$ can be defined,
for each Borel $\omega\subset\sn$, by
\begin{equation*}
C_{n-1}(K, \omega) = \mathcal H^{n-1}(r_K(\omega)),
\end{equation*}
that is, $C_{n-1}(K,\omega)$ is the surface area of the subset of $\partial K$ whose radial directions
lie in $\omega$.
Recall that the surface area measure $S(K,\eta)=S_{n-1}(K,\eta)$ is the surface area of the subset
of $\partial K$ corresponding to the set $\eta$ of
normal directions. Thus,
\begin{equation}\label{4.3}
C_{n-1}(K,\balpha^*_K(\eta)) = S_{n-1}(K, \eta),
\end{equation}
for each Borel $\eta\subset\sn$, where $\balpha^*_K$ is the reverse radial Gauss image
map \eqref{reverse-radial-gauss-image}.

The total mass of the $i$th curvature measure gives
the $(n-i)$th quermassintegral,
\begin{equation*}
C_i(K,S^{n-1}) = nW_{n-i}(K), \quad i=0,\ldots, n-1.
\end{equation*}

If $\partial K$ is $C^2$ and has positive curvature, the curvature measure $C_i(K,\cdot)$ can be written as
\[
C_{i}(K,\omega) = S_{i}(K, \balpha_K(\omega))
= \int_{r_K(\omega)} H_{n-i-1} \, d\Hn,  \ \ \  \omega \subset \sn,
\]
where $r_K$ is the radial map and $H_{n-i-1}$ is the $(n-i-1)$th mean curvature.
Since the surface area element of $\partial K$ has the formula,
\[
d\Hn = \rho_K^{n-2}(u) |\nabla \rho_K(u)|\, du,
\]
we have
\[
dC_{i}(K, u)
= H_{n-i-1}(r_K(u)) \rho_K^{n-2}(u) |\nabla \rho_K(u)|\, du,  \ \ \  u\in \sn, \ \ i=0, \ldots, n-1,
\]
where $du$ is the surface area element of $\sn$.

\subsubsection{Uniqueness Problem of Curvature Measures}
The uniqueness of a convex body with given curvature measure $C_i(K,\cdot)$ was shown by
Aleksandrov \cite{A42} for $i=0$, and
by Schneider \cite{S78} for $i=1, \ldots, n-1$. It states:
{\it Let $K$ and $L$ be convex bodies in $\kon$ and
$i\in \{0, 1, \ldots, n-1\}$. If
\[
C_i(K, \cdot)=C_i(L,\cdot),
\]
then $K=L$ when $i\in \{1,\ldots, n-1\}$ and $K$ is a dilation of $L$ when $i=0$.}

\subsubsection{Minkowski Problems of Curvature Measures}\label{AP}

This can be stated as follows:

\vs{5}

{\it Given a finite Borel measure $\mu$ on $\sn$ and $i\in \{0, 1, \ldots, n-1\}$.
What are the necessary and sufficient conditions such that there exists
a convex body $K$ in $\kon$ satisfying $C_i(K,\cdot) = \mu$?}

\vs{5}

The
cases of $i=1, \ldots, n-1$ are called the {\it Aleksandrov-Minkowski} problems.
The case $i=0$ is the {\it Aleksandrov problem} of prescribing integral curvature.
The discrete Aleksandrov problem has a clear geometric meaning. The case in $\R^2$ is described in \S\ref{ap2}.
In $\rn$, it states:

\vs{3}
{\it Given positive numbers $f_1,\ldots, f_m$ and unit vectors $u_1, \ldots, u_m$ in $\rn$, find the necessary
and sufficient conditions so that there exists a convex polytope $P$ whose vertices have outer normal
angles $f_1, \ldots, f_m$ and radial unit vectors $u_1, \ldots, u_m$.
}
\vs{3}

Aleksandrov \cite{A42} gave a complete solution to the Aleksandrov problem.

\subsubsection{Solution to the Aleksandrov Problem} Aleksandrov \cite{A42} proved the following:

\vs{5}

{\it
  If $\mu$ is a finite Borel measure on $\sn$, then
there exists a convex body $K$ in $\mathcal K_o^n$ such that
$C_0(K, \cdot) =\mu$
if and only if $\mu$ satisfies the Aleksandrov condition,
\begin{equation}\label{aleksandrov-condition}
|\mu| = n\omega_n > \Hn(\omega^*) + \mu(\omega),
\end{equation}
for every compact and spherically convex set $\omega\in \sn$, where $\omega^*$ is the polar set of $\omega$,
\[\omega^*=\{v\in \sn : v\cdot u \le 0 \ \text{for all } u \in \omega\}.\]
}

 Aleksandrov proved the discrete case for polytopes first by a topological argument using what is now known as
 the Aleksandrov  mapping lemma. He extended the solution to the general case by approximation.
 A mass transport approach to the Aleksandrov problem was given by Oliker \cite{Ol} and later differently by Bertrand \cite{B16}.

\subsubsection{Solution to the Symmetric Aleksandrov Problem}

 A longstanding open problem was whether there exists a direct variational proof,
similar to that of the Minkowski problem, demonstrating the existence of a solution
to the Aleksandrov problem. Prior to the work \cite{HLYZ16acta}, it was unknown whether
there is a geometric invariant whose differential is the integral curvature.
It was shown in \cite{HLYZ16acta} that the integral curvature is related to the differential of
the {\it entropy of convex body $K$} in $\mathcal K_o^n$,
\begin{equation}\label{ent}
E(K) = \int_{\sn} \log \rhok (u) \, du.
\end{equation}
In particular, if $\rho_t$ is a logarithmic variation of the radial function of a convex body $K$, i.e.,
\[
\log \rho\lsub{t}(v) = \log \rhok(v) + t g(v) + o(t, v)\text{ and }\lim_{t\to 0} \frac{o(t,v)}t = 0 \ \ \text{uniformly,}
\]
and $\bla{\rho\lsub{t}}\bra$ is the convex hull of $\rho_t$, see \eqref{0.19}, then
\begin{equation}\label{cm1}
\frac d{dt} E(\bla{\rho\lsub{t}}\bra^*)\Big|_{t=0} = -\int_{\sn} g(v)\, dC_0(K, v).
\end{equation}
Using the logarithmic variation of the radial function is to avoid an extra factor of radial function so that
one gets the exact differential formula above.
Using this formula, a direct variational proof of the solution to the symmetric case of
the Aleksandrov problem was given in \cite{HLYZ18jdg}.
It was also observed in \cite{HLYZ18jdg} that the Aleksandrov condition \eqref{aleksandrov-condition}
holds for any even measure on $\sn$.
Thus, for the symmetric case, we have

\vs{5}

 {\it If $\mu$ is a finite even Borel measure on $\sn$, then
there exists an origin-symmetric convex body $K$ in $\rn$ such that
$C_0(K, \cdot) =\mu$
if and only if $\mu$ is not concentrated on a great sub-sphere of $\sn$ and $|\mu|=n\omega_n$. }
\vs{3}

For the general case, a direct variational proof of Aleksandrov's theorem was given in \cite{BLYZZ20cpam}.
The associated maximization problem is the following,
\[
\sup_{\rho \in C^+(\sn)} \Big\{\frac1{|\mu|}\int_{\sn} \log \rho \, d\mu
+ \frac1{n\omega_n} E(\bla\rho\bra^*) \Big\}.
\]

When the
measure $\mu$ has a density function $f:\sn \to \R$, the partial differential
equation of the Aleksandrov problem can be written as a Monge-Amp\`ere type equation on $\sn$:
\begin{equation}\label{APp}
 \det(\nabla_{ij} h+ h\delta_{ij}) = \frac{|\nabla h|^{n}}h f,
\end{equation}
where $h$ is the unknown function on $\sn$ to be found.
A solution $h$ to this equation is the support function of $K^*$, where $K$ is a solution
to the Aleksandrov problem.

Regularity of solutions to the Aleksandrov problem was investigated
by Pogorelov \cite{P69}, Guan-Li \cite{GL97com} and Oliker \cite{Ol2}.

\subsubsection{Existence and Regularity of the Aleksandrov-Minkowski Problem}
Guan-Lin-Ma \cite{GLM09} and Guan-Li-Li \cite{GLL12duke} studied the Aleksandrov-Minkowski problem
for the cases $i=1, \ldots, n-2$ when the given measure has an everywhere positive density.
Guan-Li-Li \cite{GLL12duke} proved the following existence and regularity theorem:

\vs{5}

{\it Let $i\in \{1, \ldots, n-2\}$ and $\mu$ be a finite Borel measure on $\sn$ with
density $f\in C^2(\sn)$
and $f>0$. If $f^{-\frac1i}$ is convex when $f$ is extended to $\rn\setminus\{0\}$
homogeneously of degree $-n$,
then there exists a convex body $K$ of $C^{3,\alpha}$ and positive curvature
so that $C_i(K, \cdot) = \mu$.
}

\vs{5}

Guan-Lin-Ma \cite{GLM09} had proved the same theorem earlier but with a stronger assumption,
namely that  $f^{-\frac1i}$ is strictly convex when $f$ is extended to $\rn\setminus\{0\}$
homogeneously of degree $-n$.

\section{Dual Minkowski Problems in the Dual Brunn-Minkowski Theory}\label{dBMt}\

The dual Brunn-Minkowski theory of star bodies was introduced by Lutwak in 1970s \cite{L75paci}.
It is based on a remarkable
duality in convex geometry, related to but distinct from the standard duality between
a convex body and its polar. The duality is a heuristic concept that is used as a guiding principle.
 An excellent explanation of
this conceptual duality is given in Schneider \cite{S14} (p. 507).
The dual Brunn-Minkowski theory plays a central role in the study of intersection bodies and,
in the late 1980s and 1990s, was a crucial element in the solution to
 the Busemann-Petty problem (\cite{L88adv, G94annals, Z99annals, GKS99, K03, K05}).
For convex bodies,
 problems in the dual Brunn-Minkowski theory are often deeper and significantly more
 challenging than their analogues in the classical Brunn-Minkowski theory.
 This duality between projections of convex bodies and central intersections of star bodies
 is thoroughly discussed in Gardner \cite{G06book}.

Central to the dual Brunn-Minkowski theory are dual mixed volumes and the newly discovered
dual curvature measures.  Until recently, it was not known how to derive the
dual curvature measures using variations of dual mixed volumes.

The characterization of dual curvature measures is called the dual Minkowski problem.
The fundamental case, namely the logarithmic Minkowski problem, is the characterization
of the cone-volume measure. The general dual Minkowski problem exhibits more subtleties
and is more challenging than the classical Minkowski problem in the Brunn-Minkowski theory.

\subsection{Basics of the Dual Brunn-Minkowski Theory}

\subsubsection{Dual Mixed Volumes}

A {\it star set} in $\rn$ is a compact set $K$, where if $x \in K$, then the line segment
from $0$ to $x$ also lies in $K$. Given $u \in S^{n-1}$, let $\R^+u$
denote the ray from $o$ in the direction $u$.

If $M, N$ are star sets, then $M\cap (\R^{+}u)$ and $N\cap (\R^{+}u)$ are overlapping
line segments and their Minkowski sum is a longer line segment in the same direction.
The {\it radial sum} of $M$ and $N$ is defined to be the star set $M\wts N$,
where for each $u \in \sn$,
\begin{equation*}
(M\wts N)\cap \R^+u =M\cap (\R^{+}u) + N\cap (\R^{+}u).
\end{equation*}
There is the formula,
\begin{equation*}
\rho\lsub{M\wt + N}(x) = \rho\lsub{M}(x) + \rho\lsub{N}(x), \ \ \ x\in \rn\setminus\{0\},
\end{equation*}
where the radial function of a star set $M$,
$\rho\lsub{M} : \rn\setminus\{0\} \to (0,\infty)$, is defined to be
\[
\rho\lsub{M}(x) = \max\{t > 0 : tx\in M\}, \ \ \ x\in \rn\setminus\{0\}.
\]
If a star set $M$ contains the origin in its interior and the radial function
$\rho\lsub{M}$ is continuous, $M$ is called a {\it star body}.

Let $M_1, \ldots, M_m$ be star bodies in $\rn$.
The volume of $t_1 M_1 \wts \cdots \wts t_m M_m$, $t_1, \ldots, t_m >0$,
\begin{equation}\label{dual-combination}
V(t_1 M_1 \wts \cdots \wts t_m M_m) = \frac1n \int_{\sn}
\big(t_1\rho\lsub{M_1}(u) + \cdots +t_m\rho\lsub{M_m}(u)\big)^n \, du,
\end{equation}
is a homogeneous $n$th-order polynomial in $t_1, \dots, t_m$. It can therefore be written as
\[
V(t_1M_1\wts \cdots\wts t_mM_m)=\sum_{i_1,\ldots, i_n=1}^m \wt V(M_{i_1},\cdots,M_{i_n})\,
t_{i_1}\cdots t_{i_n},
\]
where the coefficients
$\wt V(M_{i_1},\cdots,M_{i_n})$ are nonnegative, symmetric functions of the bodies
$M_1,\cdots,M_m$. They are called {\it dual mixed volumes}.
By \eqref{dual-combination}, they have the integral representation
\begin{equation*}
\wt V(M_{i_1},\cdots,M_{i_n})
=\frac1n\int_{S^{n-1}}\rho\lsub{M_{i_1}}(u)\cdots\rho\lsub{M_{i_n}}(u)\,du.
\end{equation*}

For $q \in \R$, the {\it $q$th dual mixed volume} of star bodies $M$ and $N$ is defined to be
\[
\wt V_q(M,N)=\frac1n\int_{S^{n-1}}\rho_M^{n-q}(u) \rho_N^q(u)\, du,
\]
and the {\it $q$th dual quermassintegral}
 $\wt W_{q}(M)$ of a star body $M$ is
\begin{equation*}
\wt W_{q}(M) = \wt V_q(M,B) =\frac1n \int_{S^{n-1}}\rho_M^{n-q}(u)\,du,
\end{equation*}
where $B$ is the unit ball in $\rn$.

Let
\begin{equation}\label{dqmi}
\wt V_{q}(K) =\wt W_{n-q}(K)= \frac1n \int_{\sn} \rho_K^q (u)\, du.
\end{equation}

When $q\neq0$, the dual quermassintegral $\wt V_{q}(K)$ can be written as
\[
\wt V_{q}(K) = \begin{cases}\ds \frac qn \int_K |x|^{q-n} dx &q>0 \\
                            \ds \frac {|q|}n \int_{\rn \setminus K} |x|^{q-n} dx &q<0,  \end{cases}
\]
which is related to the Riesz potential of the characteristic functions $\chara K$ and $\chara{\rn\setminus K}$.

The integer cases $\wt V_1(K), \ldots, \wt V_n(K)$ have a clear geometric meaning.
They are integrals of lower dimensional volumes of intersections of $K$ with subspaces,
\begin{equation}\label{dqxi}
\wt V_i(K) = \frac{\omega_n}{\omega_i} \int_{G_{n,i}} V_i(K\cap\xi)\, d\xi, \ \ \ i=1,\ldots, n.
\end{equation}

The normalized dual quermassintegral $\barV_q(K)$ is defined to be
\[
  \barV_q(K) = \begin{cases}
    \Big(\dfrac{\wt V_q(K)}{\omega_n}\Big)^\frac1q &\text{ if }q\ne 0\\[10pt]
    \exp\Big(\displaystyle\frac1{n\omega_n}\int_{\sn} \log \rhok(u)\, du\Big) &\text{ if }q = 0.
  \end{cases}
\]

\subsubsection{Dual Curvature Measures}
Huang-LYZ \cite{HLYZ16acta} discovered the geometric measures associated with dual quermassintegrals, called
{\it dual curvature measures}.
For a convex body $K$ that contains the origin in its interior and $q\in \R$,
the {\it $q$th dual curvature measure} $\wt C_q(K,\cdot)$ is defined by
\begin{equation}\label{5.4}
\wt C_q(K,\eta) = \frac1n \int_{\balpha^*_K(\eta)} \rho_K^q(u)\, du
= \frac1n \int_{\sn} \chara{\balpha^*_K(\eta)}(u) \rho_K^q(u)\, du,
\end{equation}
for each Borel set $\eta\subset\sn$.
It can also be defined by an integral over the boundary $\partial K$,
\[
\wt C_q(K,\eta)=\frac1n\int_{\nu_K^{-1}(\eta)} (x\cdot \nuk(x)) |x|^{q-n} \, dx,
\]
or as a linear functional,
\[
\int_{\sn} f(v)\, d\wt C_q(K,v) = \frac1n \int_{\sn} f(\alpha\lsub K(u)) \rho_K^q(u)\, du.
\]

For a polytope $P$ with outer unit normals $v_1,\ldots,v_m$, the dual curvature measure of $P$ is a
discrete measure given by
\[
\wt C_q(P, \cdot) = \sum_{i=1}^m  c_i \delta_{v_i}\text{, where } c_i=\int_{\balpha_P^*(v_i)} \rho_K^q(u)\, du.
\]

If $K$ is strictly convex, then the dual curvature measure can be expressed in terms of
the surface area measure and support
function,
\[
d\wt C_q(K,\cdot) = \frac1n h_K |\nabla h_K|^{q-n} dS(K,\cdot).
\]

When $q=0$, $\wt C_0(K,\cdot)$ is called the {\it radial angle measure}. For a polytope $P$,
\[
\wt C_0(P,\cdot) = \sum_{i=1}^m  c_i \delta_{v_i}, \ \ \ c_i=\Hn(\balpha_P^*(v_i)),
\]
where the weights $c_i$ are the radial angles of the cones formed by the origin and facets of $P$.

Up to a constant factor, the radial angle measure of a convex body is
the integral curvature of the polar body,
\[
\wt C_0(K,\cdot)=\frac1n C_0(K^*,\cdot).
\]

When $q=n$, $\wt C_n(K,\cdot)$ is the cone-volume measure of $K$, denoted by $V(K, \cdot)$ or $V_K$,
\[
d\wt C_n(K,\cdot) = \frac1n h_K dS(K,\cdot) = dV(K,\cdot).
\]

The fact that cone-volume measure and integral curvature are related to each other was not known until
it was discovered that both measures are dual curvature measures.

The concept of dual curvature measures in the dual Brunn-Minkowski theory corresponds to the concept of
curvature measures in the Brunn-Minkowski theory, see \cite{HLYZ16acta} for detailed explanations.
Here, we mention the duality of geometric constructions of curvature measures and dual curvature measures.

Let $K$ be a convex body that contains the origin in its interior in $\rn$.
Define the {\it radial projection map} $\wt p\lsub{K} : \rn \setminus K \to \partial K$ by
\[
\wt p\lsub{K}(x) = \rhok(x)x,
\]
for $x \in \rn \setminus K$. For $x \in \rn$, the {\it radial distance} $\wt d(K,x)$ of $x$ to $K$, is defined by
\[
\wt d(K,x) = \begin{cases} |x-\wt p\lsub{K}(x)|  & x\notin K \\
0 &x\in K.
\end{cases}
\]

For $t>0$, consider the radial annulus with thickness $t$,
\[
\wt K_t = \{x\notin K : \wt d(K, x) \le t\} = (K \wts tB) \setminus K.
\]

For $t > 0$ and a Borel set $\eta \subset \sn$, define the local radial annulus,
\begin{equation}
\wt A_t(K,\eta) =\wt K_t \cap \wt p^{-1}_{K}(\bnu_K^*(\eta)).
\label{dcs}
\end{equation}

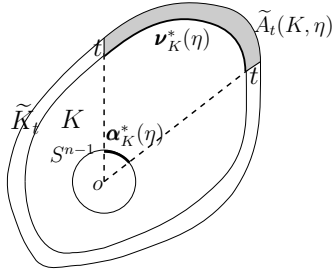
\begin{center}
\scalebox{0.7}{
\begin{tikzpicture}
	 [scale=0.3]
	\draw (0,0) circle [radius=2cm]
	node at (-2,1.8) {$ \sn$}
	node at (-0.4,-0.2) {$o$}
	node at (2,3) {${\balpha^*_K(\eta)}$}
	node at (-5,4) {\large{$\wt K_t$}};
	\draw [line width=1.5pt] (1.6,1.2) arc (40:92:2) ;
	
    \draw (-5,0) to [out=90,in=-140] (0,8);
    \draw [line width=2pt](0,8) to [out=40,in=160] (7,10) to [out=-20,in=90] (9,7)
    node at (5, 9) {$\bnu^*_K(\eta)$};
    \draw (9,7) to [out=-90,in=45] (7.3,0) to [out=-135,in=20] (0,-4.2) to [out=-170,in=0] (-2,-4.5) to [out=190,in=-90](-5,0)
	node at (12,10) {$\large \widetilde{A}_t (K,\eta)$}
	node at (-2,4){\large $K$};
	
	\fill[fill=white!80!black]
	(9.9,7.6) to [out=90,in=-17] (7.8,11)to [out=165,in=40] (0,9.1) to (0,8) to [out=40,in=160] (7,10) to [out=-20,in=90] (9,7) to (9.9,7.6) ;
	
	\draw (9.9,7.6) to [out=90,in=-17] (7.8,11)to [out=165,in=40] (0,9.1)to [out=-140,in=50] (-4,5.2) to [out=-130,in=80] (-6.1,0.5) to [out=-95,in=110] (-5.4,-3) to [out=-70,in=125] (-4.3,-5) to [out=-45,in=180] (-1.5,-5.3) to [out=10,in=-160](1.6,-4.6) to [out=20,in=-130](8.1,-0.8) to [out=55,in=-90](9.9,7.6) ;

	\draw[dashed,thick](0,0)--(0,8);
	\draw[dashed,thick](0,0)--(9,7);
	
	\draw [line width=0.6pt](0,8) to (0,9.1) node at (-0.4,8.5) {\large$t$};
	\draw [line width=0.6pt](9,7) to (9.9,7.6) node at (9.5,6.7) {\large$t$};
	
	
\end{tikzpicture}
}
\end{center}

Then there is the dual Steiner-type formula,
\begin{equation}
V(\wt A_t(K,\eta))= \sum_{i=0}^{n} \binom ni t^{n-i} \wt C_i(K,\eta). \label{dcmf}
\end{equation}

Formulas \eqref{dcs} and \eqref{dcmf} are dual to formulas \eqref{cmp} and \eqref{cme}.

Dual curvature measures are differentials of dual quermassintegrals, as shown by
the following differential formulas,
\begin{equation}\label{df1}
\frac{d}{dt} \log \barV_q(\bla \rho_t \bra^*)\Big|_{t=0}
= - \frac{1}{\wt V_q(\bla \rho_{0} \bra^*)}
\int_{\sn} g(u) \, d\wt C_q(\bla \rho_{0} \bra^*, u),
\end{equation}
where $\bla\rho_t\bra$ is the convex hull of the logarithmic family of variations $\rho_t$ defined by
\[
\log\rho_t(u) = \log \rho_0(u) + t g(u) + o(t,u),  \ \ \ \lim_{t\to 0}\frac{o(t,u)}t=0 \ \text{ uniformly.}
\]
The term $o(t,u)$ is for generality which makes it easy to reformulate in the dual case. That is,
it can also be formulated as differential formulas of dual quermassintegrals of Wulff shapes,
\begin{equation}\label{df2}
\frac{d}{dt} \log \barV_q(\blb h_t \brb)\Big|_{t=0}
= \frac{1}{\wt V_q(\blb h_0 \brb)}
\int_{\sn} f(v) \, d\wt C_q(\blb h_0 \brb, v),
\end{equation}
where $\blb h_t \brb$ is the Wulff shape of the logarithmic family of variations $h_t$ defined by
\[
\log h_t(v) = \log h_0(v) + t f(v) + o(t,v), \ \ \ \lim_{t\to 0}\frac{o(t,v)}t=0 \ \text{ uniformly.}
\]
When $q=n$, \eqref{df2} is equivalent to \eqref{b2.2}. Similar to the role of Aleksandrov's differential
formula of volume in the classical Minkowski problem, differential formulas \eqref{df1} and \eqref{df2}
are critical for solving the Minkowski problems of dual curvature measures.
The proofs of \eqref{df1} and \eqref{df2}
were given by Huang-LYZ \cite[Thereoms 4.4 and 4.5]{HLYZ16acta}.

\subsubsection{Cone-Volume Measure}
Cone-volume measure $V(K,\cdot)$ of a convex body $K\in \kon$ can be defined more geometrically.
For a Borel set
$\sigma \subset \partial K$, let the {\it radial cone of $\sigma$} be the convex hull of $\sigma$ and the origin
(see the figure in \S\ref{rGm}),
\begin{equation}\label{radial-cone}
{\scriptstyle\triangle}_{K}(\sigma)=\{sx: 0\le s \le 1,\  x\in \sigma\}.
\end{equation}
If $K$ is a polytope, we will call the radial cone formed by the origin and
a facet of $K$ a {\it facet radial cone} of the polytope.

Recall from \eqref{reverseGaussImage} that the reverse Gauss image of a subset $\eta \subset \sn$,
denoted $\bnu_K^*(\eta)$ , is the set of points in $\partial K$ whose normal cones intersect $\eta$.
The cone-volume measure of a Borel subset $\eta \subset \sn$ is defined to be
\[
V(K, \eta) = \HH^n({\st \triangle}_{K}(\bnu_K^*(\eta)),\quad \eta\subset \sn.
\]

For example, the cone-volume measure of a polytope $P$, whose facets have
radial cones ${\st\triangle}_1,\ldots,{\st\triangle}_m$ and corresponding
outer unit normals $v_1,\ldots,v_m$,
is the discrete measure,
\begin{equation}\label{dcvm}
V(P, \cdot) = \sum_{i=1}^m  \HH^n({\st\triangle}_i) \, \delta_{v_i},
\end{equation}
where $\HH^n({\st\triangle}_i)$ is the volume of ${\st\triangle}_i$ and $\delta_{v_i}$
is the delta measure at $v_i\in\sn$.

The volume of the radial cone can be written as an integral either on the unit sphere or on the boundary
of the convex body using the radial function or support function,
\begin{equation*}
V(K, \eta) = \int_{{\st\triangle}_K(\bnu_K^*(\eta))} dx =\frac1n \int_{\balpha_K^*(\eta)} \rho_K^n(u)\, du
= \frac1n \int_{\nu_K^{-1}(\eta)} (x\cdot \nuk(x))\, dx.
\end{equation*}
The last equality above shows that cone-volume measure and surface area measure satisfy the formula,
\begin{equation}\label{5.1}
dV(K, \cdot) = \frac1n h_K \, dS(K,\cdot).
\end{equation}

Given an affine transformation $\phi \in SL(n)$, the cone-volume measure of $\phi K$ satisfies
\[
V(\phi K,\eta) = V(K, \overline{\phi^{-t}\eta}), \quad \eta\subset \sn, \ \ \phi\in SL(n),
\]
where
\[
  \overline{\phi^{-t}\eta}=\left\{\frac{\phi^{-t}v}{|\phi^{-t}v|} : v\in \eta\right\}
\]
and $\phi^{-t}$ is the inverse transpose of $\phi$.

When the boundary of $K$ is $C^2$ and has positive Gauss curvature $\lkappa$,
the cone-volume measure of $K$ has a positive density,
\[
dV(K,v) =  \frac{h_K}{n\lkappa}\, dv.
\]

\subsection{Logarithmic Minkowski Problem of Cone-Volume Measure}\label{5.2}

The Minkowski problem for cone-volume measure is the following.
\vs{3}

\fr{Logarithmic Minkowski Problem} {\it Given a finite Borel measure $\mu$ on $\sn$,
find the necessary and sufficient conditions such that there exists a convex body $K$ in $\rn$
whose cone-volume measure is $\mu$, that is, the measure equation $V(K,\cdot) = \mu$ has a solution $K$.}

\vs{3}

When the cone-volume measure of a convex body is discrete, then the body is a polytope
and the measure is given by \eqref{dcvm}. The logarithmic Minkowski problem can therefore be viewed
as a problem of finding a convex polytope with prescribed outer unit normals and cone volumes
of the corresponding facet radial cones.
The discrete logarithmic Minkowski problem in $\R^2$ is described in \S\ref{stancu}. In $\rn$, it states:
\vs{3}

\fr{Discrete Logarithmic Minkowski Problem}
{\it Given positive numbers $f_1, \ldots, f_m$, and unit vectors $v_1, \ldots, v_m$ in $\rn$,
find the necessary and sufficient conditions such that there exists a convex polytope $K$ in $\rn$
whose facet radial cones have volumes $f_1, \ldots, f_m$ with corresponding facet outer unit normals $v_1, \ldots, v_m$.}

\vs{3}

When the measure $\mu$ has a positive density $f$, the question can be
formulated as an elliptic Monge-Amp\`ere equation,
\begin{equation}\label{logp}
\det(\nabla_{ij} h + h\delta_{ij}) = \frac fh.
\end{equation}

The first results in $2$-dimension were given
by Gage and Li \cite{GL94} for measures with smooth densities and by Stancu \cite{St02} for discrete measures.
Chou and Wang \cite{CW06} extended Gage and Li's result to higher dimensions.
B\"or\"oczky-LYZ \cite{BLYZ13jams} then solved
the problem for arbitrary even measures in all dimensions.
There are only partial results for the logarithmic Minkowski problem for arbitrary measures,
because the obstructions to a solution are not fully understood.

The necessary and sufficient conditions for an even measure
to be the cone-volume measure of a symmetric convex body are known as the subspace concentration condition.

\subsubsection{Subspace Concentration Condition}\label{sec:ssc}

For an origin-symmetric convex body $K \subset \rn$, we first look at how large
the cone-volume measure $V_K(\xi\cap \sn)$ can be for a 1-dimensional subspace $\xi$.
If $\xi$ is spanned by unit vector $v$, then
\begin{equation}\label{cone}
V_K(\xi\cap \sn)= V_K(\{v,-v\}) = \frac{2}{n}h F,
\end{equation}
where $h = h_K(v)$, $F= \mathcal{H}^{n-1}(K_v)$, and $K_v=\partial K \cap H_K(v)=\bnu_K^*(\{v\})$.

Note that the volume of the convex hull of $K_v$ and $K_{-v}$ is at least $2hF$.
This is because the intersection of a hyperplane (parallel to $H_K(v)$) with the convex hull is
$t K_v + (1-t) K_{-v}$ for some $0\le t\le 1$. By the Brunn-Minkowski inequality,
\[
\Hn(t K_v + (1-t) K_{-v}) \ge (tF^\frac1{n-1} +(1-t)F^\frac1{n-1})^{n-1} = F.
\]
Since the convex hull is contained in $K$, we have
$2hF \le V(K)$.
By this and \eqref{cone}, it follows that
\begin{equation}\label{1-concentration}
  V_K(\xi\cap \sn) \le \frac{1}{n}V(K).
\end{equation}
This motivates a general condition for cone-volume measure:
\vs{3}

{\it A finite Borel measure $\mu$ on $\sn$ is said to satisfy the
{\it subspace concentration condition} if, for any subspace $\xi$ of
$\rn$,
\begin{equation}\label{ssc}
 \mu(\xi\cap \sn)\leq \frac {\dim \xi}n \, \mu(\sn),
\end{equation}
and equality holds if and only if there exists a linear subspace
$\xi' \subset \rn$ complementary to $\xi$ such that the support of $\mu$ lies in $\sn\cap(\xi\cup \xi')$.}

\vs{3}

The measure $\mu$ on $S^{n-1}$ is said to satisfy the {\it
strict subspace concentration inequality}
if the inequality in \eqref{ssc} is strict for each subspace $\xi
\subset \rn$ such that $0<\dim \xi<n$.

He-Leng-Li \cite{HLL} and Henk-Sch\"urmann-Wills \cite{HSW} independently showed that
the cone-volume measure of a symmetric polytope satisfies the subspace concentration condition,
with an alternate proof given by Xiong \cite{X1}.
B\"or\"oczky-LYZ \cite{BLYZ13jams} showed that the cone-volume measure of any symmetric
convex body satisfies the subspace concentration condition.

If $K$ is not necessarily symmetric but its centroid is at the origin, then its
cone-volume measure satisfies the subspace
concentration condition. This was proved for polytopes in $\R^3$ by Xiong \cite{X1},
for polytopes in $\rn$
by Henk-Linke \cite{HL14}, and for all convex bodies by  B\"or\"oczky-Henk \cite{BH16}.

Applications of the subspace concentration condition to reverse affine isoperimetric
inequalities were given in \cite{HLL, X1, HL14, BH16}. Connections of it to the so-called
$K$-stability property in algebraic geometry \cite{CDS15I, CDS15II, CDS15III} was observed by
Hering-Nill-S\"uss  \cite{HNS19}.

\subsubsection{Solution to the Symmetric Logarithmic Minkowski Problem}\label{5.3}

 If the logarithmic Minkowski problem is restricted to even measures and symmetric convex bodies,
 it is called the symmetric logarithmic Minkowski problem.

The symmetric logarithmic Minkowski problem was solved for polygons in $\R^2$
by Stancu \cite{St02}. For arbitrary convex bodies in all dimensions, a complete solution to
the symmetric logarithmic Minkowski problem was established by B\"or\"oczky-LYZ \cite{BLYZ13jams}.
They proved the following:

\vs{3}

{\it A non-zero finite even Borel measure on the unit sphere $\sn$ is
the cone-volume measure of an origin-symmetric convex body
in $\rn$ if and only if it satisfies the subspace concentration
condition.}

\vs{3}

The proof in \cite{BLYZ13jams} uses a variational method.  A Wulff shape construction is needed,
as in the proof of the Minkowski problem, to define a differentiable functional whose minimum
gives a solution to the logarithmic Minkowski problem. In particular, given a measure $\mu$ on $\sn$,
if there is a function $h \in C_e^+(\sn)$ that minimizes the functional
\begin{equation}\label{lmp}
\Phi_\mu[h] = V(\blb h\brb)^{-\frac1n} \exp\int_{\sn} \log h \, d\mu,
\end{equation}
then there exists a convex body $K$ whose support function $h_K$ also minimizes this functional.
That $h_K$ satisfies the Euler-Lagrange equation for the functional implies that the cone-volume measure
of $K$ is a scalar multiple of $\mu$. Rescaling $K$ then proves the theorem.

We now need to show that a minimizer exists.
Since the functional of the minimization problem is monotone, it is easily shown that
there exists a minimizing sequence, where the functions are all support functions of convex bodies.
Since the functional is invariant under rescaling of the bodies, one can assume that the convex bodies
have bounded diameter. It follows that there exists a subsequence that converges to a convex set.
The difficult part is proving that the limit is in fact a convex body, i.e, it has a nonempty
interior. The proof of this uses sharp estimates for both the volume $V(\blb h\brb)$ of the Wulff shape
and the integral $\int_{\sn} \log h \, d\mu$. The main obstacle is to deal with the concentration of
$\mu$ in subspaces.
To overcome it, the crucial idea is to construct a spherical partition that can be used to express
the concentrations of $\mu$ in subspaces of all dimensions.
Then the subspace concentration condition is invoked to obtain sharp estimates
for the integral in the minimization functional. This technique is also important
in solving other Minkowski problems.

\subsubsection{Solution to the Non-Symmetric Logarithmic Minkowski Problem}

When a convex body is not origin-symmetric, its cone-volume measure may not satisfy
the subspace concentration condition.
However, Zhu \cite{Zhu1}, B\"or\"oczky-Heged\H{u}s-Zhu \cite{BHZ16}, and Chen-Li-Zhu \cite{CLZ19}
showed that the subspace concentration condition is sufficient for prescribing
cone-volume measure.

A linear subspace $\xi$ in $\rn$, $0< \dim\xi <n$, is called an {\em essential subspace}
with respect to a Borel measure
$\mu$ on $\sn$ if $\xi \cap \text{supp}\mu$ is not contained in a closed hemisphere
of the subsphere $\xi \cap \sn$.
If $\mu$ is discrete and the vectors in its support  are in general position, i.e.,
any $k$ vectors in supp$\xi$,
$1\le k\le n$, are linearly independent, then it is easily seen there are no
essential subspaces with respect to $\mu$.
\vs{3}

For arbitrary convex polytopes, the following solution to the logarithmic Minkowski problem was proved:
\vs{3}

{\it If $\mu$ is a discrete measure on $\sn$ that is not concentrated on any closed
hemisphere and satisfies the subspace concentration condition for all essential subspaces,
then $\mu$ is the cone-volume measure of a polytope in $\rn$ containing the origin in its interior.}
\vs{3}

For polytopes whose outer unit normals are in general position, this was first proved
by Zhu \cite{Zhu1}. B\"or\"oczky-Heged\H{u}s-Zhu \cite{BHZ16} proved it for general convex polytopes.
Chen-Li-Zhu \cite{CLZ19} then used an approximation argument to prove:
\vs{3}

{\it  If $\mu$ is a non-zero finite Borel measure on $\sn$ that satisfies
the subspace concentration condition, then there exists a convex body $K$ containing
the origin in $\rn$ such that $V(K,\cdot)=\mu$.}
\vs{1}

The necessary and sufficient condition for the non-symmetric case is still unknown.
Liu-Lu-Sun-Xiong \cite{LLSX24} found a necessary condition for a convex body in $\R^2$
with centroid at the origin and solved the logarithmic Minkowski problem for quadrilaterals.

\subsubsection{Minkowski Problems and Self-Similar Solutions of Gauss Curvature Flows}

Another approach to solving Minkowski problems for smooth solutions is using
anisotropic Gauss curvature flow equations. One considers an anisotropic Gauss curvature
flow equation related to the partial differential equation of a Minkowski problem,
and normalize the flow equation by rescaling so that the self-similar solution of
the flow equation satisfies the PDE.

Let $X(\cdot,t):  \sn \to \rn$ be a parametrization
of the boundaries of a family of smooth convex bodies $K_t \in \kon$ with positive
Gauss curvature $\lkappa(\cdot, t)$ and outer unit normal $\nu(\cdot, t)$ at $X(\cdot, t)$,
and let $f$ be a positive smooth function on $\sn$. For the logarithmic Minkowski problems,
the anisotropic Gauss curvature flow equation is
\begin{equation}\label{aGf}
\frac{\partial X}{\partial t} (u,t) = -  \lkappa(u,t) f(\nu(u,t)) \nu(u,t).
\end{equation}

When $f=1$, this is the Gauss curvature flow first studied by Firey \cite{F74},
subsequently by Tso \cite{T85}, Chow \cite{Chow87},
Andrews \cite{And99}, Andrews-Chen \cite{AC12}, and Huang-Liu-Xu \cite{HLX15}.
The existence and regularity aspects of the Gauss curvature flow were resolved by Tso \cite{T85}.
The smooth convergence to a ball in $\R^3$ was first proved by Andrews \cite{And99}.
For higher dimensions, Guan-Ni \cite{GN17} proved the smooth convergence to a self-similar solution,
namely a soliton. Brendle-Choi-Daskalopoulos \cite{BCD17} proved that the soliton is a ball,
also proved in Choi \cite{Choi17}.
The general case was studied by Chou-Zhu \cite{CZ01}, Andrews \cite{And00}, and others.

The main technical part of the approach is to solve the normalized anisotropic Gauss curvature flow equation
for all $t>0$ and prove that the solution converges to a smooth self-similar solution as $t\to\infty$.
This requires finding a priori estimates of solutions to the flow equation.

This approach has been applied to finding smooth solutions of many other Minkowski problems
of geometric measures. For example,
Chou-Wang\cite{CW00} introduced the following logarithmic Gauss curvature flow
\begin{equation}\label{lGf}
\frac{\partial X}{\partial t}(u,t)=-{\rm log}\frac{\lkappa(u,t)}{f(\nu(u,t))}\nu(u,t),
\end{equation}
and showed  that its solution converges to the smooth solution of the classical Minkowski problem.

\subsection{Logarithmic Brunn-Minkowski Conjecture}

\subsubsection{Uniqueness Conjecture of Cone-Volume Measure}
Solutions to the logarithmic Minkowski problem
need not be unique.
In particular, if the cone-volume measure of a symmetric convex body $K$ is supported
on a union of linear subspaces $\xi_1, \dots, \xi_m \subset \rn$ of positive dimension such that
$\rn = \xi_1 \oplus \cdots \oplus \xi_m$, then
\[
  K = (K\cap \xi_1) + \cdots + (K\cap\xi_m).
\]
In that case, given any $a_1, \dots, a_m > 0$ such that $a_1\cdots a_m = 1$, the convex body
\begin{equation}\label{decomposition}
  K' = (a_1K\cap \xi_1) + \cdots + (a_mK\cap\xi_m)
\end{equation}
has the same cone-volume measure as $K$. This leads to the following conjecture:

\vspace{3pt}

\fr{Uniqueness Conjecture of Cone-Volume Measure}  {\it If $K, L \subset \rn$ are symmetric convex bodies such that
  $V(K,\cdot) = V(L,\cdot)$, then there exists a direct sum decomposition
  $\rn = \xi_1 \oplus \cdots \oplus \xi_m$
  such that
  \[
    K = (K\cap \xi_1) + \cdots + (K\cap\xi_m), \ \ \ L = (a_1K\cap \xi_1) + \cdots + (a_mK\cap\xi_m),
  \]
  where $a_1\cdots a_m = 1$, $a_1, \ldots, a_m >0$.
}

\vspace{3pt}

Note that a direct sum decomposition with $m>1$ is not a smooth convex body. Thus,
if $K$ and $L$ are smooth, then the decomposition is trivial with $m=1$ and $\xi_1=\rn$, and therefore
the question becomes whether the following holds:
\begin{equation}\label{logBMP}
V(K,\cdot) = V(L,\cdot) \ \ \ \Longrightarrow \ \ \ K=L.
\end{equation}

When $K$ and $L$ are symmetric convex bodies in $\R^2$, the conjecture was proved
by B\"or\"oczky-LYZ \cite{BLYZ12adv}. Ma \cite{Ma15} gave a different proof, and
Xi-Leng \cite{XL16} extended this to a larger class of not necessarily symmetric
$2$-dimensional bodies. In earlier work, Gage \cite{Ga93} proved the conjecture
for smooth symmetric bodies with positive curvature in $\R^2$,
and Stancu \cite{St03} did it for symmetric polygons.
If $K$ and $L$ are convex bodies that are assumed to be symmetric about
coordinates hyperplanes, the conjecture was established by Saroglou \cite{Saro1}.

When the bodies are assumed to be smooth, there are many partial results.
In 1974 Firey \cite{F74} proved \eqref{logBMP} when $K$ is smooth and symmetric
and $L$ is a ball centered at the origin. In other words, a smooth symmetric
convex body $K$ with constant cone-volume measure is a ball.  Firey conjectured
that this holds even if $K$ is not assumed to be symmetric:
\vs{3}

{\it If the ratio of the
support function and Gauss curvature of a convex body $K$ in $\rn$ is constant,
${h_K}/{\lkappa} = constant$,
then $K$ is a ball.}
\vs{3}

Firey's conjecture
for $3$-dimensional convex bodies was proved by Andrews \cite{And99}.
Andrews-Guan-Ni \cite{AGN16} obtained partial results in $\rn$.
The conjecture was fully solved by Choi \cite{Choi17}. Further results
were proved by Brendle-Choi-Daskalopoulos \cite{BCD17}.

Kolesnikov and Milman \cite{KM18} proved the uniqueness conjecture of cone-volume measure
 in a range of cases that includes
when $K$ and $L$ are positively curved $C^2$ symmetric convex bodies lying in a sufficiently
small $C^2$-neighborhood  of the unit ball in the space $\ell_q^n$, $q\ge 2$,
when the dimension $n$ is sufficiently large.
Chen-Huang-Li-Liu \cite{CHLL20} proved it when $K$ is a $C^{2,\alpha}$ symmetric body
that lies in a sufficiently small neighborhood of the ball $B^n$ and $L$ is any
symmetric convex body in $\rn$.

\subsubsection{Logarithmic Minkowski Inequality}

B\"or\"oczky-LYZ \cite{BLYZ12adv, BLYZ13jams} observed that an affirmative answer to
the uniqueness conjecture of cone-volume measure implies the following logarithmic Minkowski inequality:

\vspace{3pt}

{\em
Given symmetric convex bodies $K, L \subset \rn$,
\begin{equation}\label{lMi}
\int_{\sn} \log\frac{h_L}{h_K}\, dV(K,\cdot) \ge \frac{V(K)}n \log\frac{V(L)}{V(K)}.
\end{equation}
}

The implication can be shown as follows. Given a symmetric convex body $K$, note that
\eqref{lMi} is equivalent to
\begin{equation}\label{Phi-inequality}
  \Phi_\mu[h_K] \le \Phi_\mu[h_L],
\end{equation}
where $\Phi_\mu$ is the functional defined by \eqref{lmp} and $\mu = V_K$. On the other hand,
the proof of B\"or\"oczky-LYZ \cite{BLYZ13jams} shows that there exists a symmetric convex body $K'$
such that $V_{K'} = V_K$ and, for any symmetric convex body $L$,
\begin{equation}\label{Phi-inequality2}
  \Phi_\mu[h_{K'}] \le \Phi_\mu[h_L].
\end{equation}
If the uniqueness conjecture holds, then $K$ and $K'$ satisfy \eqref{decomposition},
which in turn implies that $\Phi_\mu[h_K] = \Phi_\mu[h_{K'}]$. Therefore, \eqref{Phi-inequality}
holds for any symmetric convex body $L$.

Thus, the uniqueness results cited in the previous section all imply the logarithmic Minkowski inequality
under the assumptions described. Precisely, inequality \eqref{lMi} is true if $K$ is an
origin-symmetric convex body in a $C^{2,\alpha}$ neighborhood of the Euclidean unit ball \cite{CHLL20},
$K$ is an origin-symmetric convex body in $\R^2$ \cite{BLYZ12adv}, or $K$ is a convex body in $\rn$
that is symmetric about coordinates hyperplanes  \cite{Saro1}.

\subsubsection{Logarithmic Brunn-Minkowski Inequality}

Given symmetric convex bodies $K, L$ and $0 < t < 1$, define the geometric mean
of $K$ and $L$ to be the body $K^{1-t}L^t$ to be the Wulff shape of the geometric mean
$h_K^{1-t}h_L^t$.

The logarithmic Minkowski inequality \eqref{lMi} is equivalent to the following conjecture \cite{BLYZ12adv}:

\vspace{3pt}

{\it
  If $K, L \subset \rn$ are symmetric convex bodies, then
\begin{equation}\label{lBMi}
V(K^{1-t}L^t) \ge V(K)^{1-t} V(L)^t, \quad 0<t <1.
\end{equation}
}

Since $(1-t)K + t L \supseteq K^{1-t}L^t$,
the logarithmic Brunn-Minkowski inequality \eqref{lBMi} implies the Brunn-Minkowski inequality
\eqref{BMI} for symmetric convex bodies.

Since the logarithmic Minkowski inequality is equivalent to the logarithmic Brunn-Minkowski inequality,
the results described in the previous section also imply the logarithmic Brunn-Minkowski inequality.

Colesanti-Livshyts-Marsiglietti \cite{CLM17} and Colesanti-Livshyts \cite{CLi21}
proved \eqref{lBMi} when $K$ and $L$ are in a small $C^2$ neighborhood of the
Euclidean unit ball $B^n$.
Kolesnikov-Milman \cite{KM18} showed \eqref{lBMi} when $K$ and $L$ are in a small $C^2$ neighborhood of
 the unit ball of the space $\ell_q^n$, $q\ge 2$ when the dimension $n$ is sufficiently large.

\subsection{Minkowski Problems of Dual Curvature Measures}

In \cite{HLYZ16acta}, Huang-LYZ posed the following problem for dual curvature measures.
\vs{3}

\subsubsection{Dual Minkowski Problem}
{\it Given a finite Borel measure $\mu$ on $\sn$ and $q\in \R$,
what are necessary and sufficient conditions for the existence of
a convex body $K\in \kon$ solving the measure equation,
\[
\wt C_q(K,\cdot) = \mu\,?
\]
}

Since the $n$th dual curvature measure $\wt C_n(K,\cdot)$ is the
cone-volume measure $V(K,\cdot)$, the case of $q=n$ of the dual Minkowski problem is
the logarithmic Minkowski problem for cone-volume measure.
Since the $0$th dual curvature measure $\wt C_0(K,\cdot)$ is
Aleksandrov's integral curvature (the $0$th curvature measure on $\sn$)
of the polar body $K^*$ (with the constant factor $1/n$),
the case of $q=0$ of the dual Minkowski problem is equivalent to the Aleksandrov problem.
The logarithmic Minkowski problem and
the Aleksandrov problem were thought to be two entirely different problems.
They are now unified as special cases of the dual Minkowski problem. Because of the
geometric meaning of dual quermassintegrals about sections of convex bodies, the other
integer cases, $q=1,\ldots, n-1$, of the dual Minkowski problem, are of special significance.

When the
measure $\mu$ has a density function $f:\sn \to \R$, the partial differential
equation of the dual Minkowski problem can be written as the following Monge-Amp\`ere type equation on $\sn$:
\begin{equation}\label{dualp}
 \det(\nabla_{ij} h+ h\delta_{ij}) = \frac{|\nabla h|^{n-q}}h f.
\end{equation}

To solve the dual Minkowski problem, the paper \cite{HLYZ16acta} used a new variational method and considered the
maximization problem,
\begin{equation}\label{dualMax}
\max_{g \in C_e^+(\sn)} \Big\{\frac1{|\mu|} \int_{\sn} \log g \, d\mu +\frac1q
\log\Big(\frac1{n\omega_n}\int_{\sn} h_{\sbla g \sbra}^{-q}\, dv\Big)\Big\},
\end{equation}
where $\bla g \bra$ is the convex hull of $g$, see \eqref{0.19},
while the classical variational method for the Minkowski problem
uses the Wulff shape. The use of convex hull for the dual Minkowski problem is quite natural because
the convex hull and the Wulff shape are dual concepts (see \cite{HLYZ16acta}).
The variational formulas \eqref{df1} and \eqref{df2} are the crucial tools to
show that a solution to the maximization problem gives a solution to the dual Minkowski problem.

Sufficient conditions, which are also necessary conditions for dual curvature measures,
for the existence of a solution to the maximization problem \eqref{dualMax} are crucial
for obtaining sharp estimates for
the two integrals in the maximizing functional --- the entropy-type integral and the dual quermassintegral.
The entropy-type integral with respect to the given measure $\mu$
was dealt with earlier in the paper \cite{BLYZ13jams}. The estimates obtained in that paper for the logarithmic
Minkowski problem are not enough for the dual Minkowski problem and more delicate estimates are needed,
but the technique of spherical partition can still be used.
Estimates for the dual quermassintegrals of degree $q > 0$
are much more difficult to obtain than the known case of $q = n$, where the dual quermassintegral is just volume.
Sharp estimates for the dual quermassintegral (the second integral in the maximizing functional) are obtained by
the application of generalized spherical coordinates. The sharp estimates for both the
entropy-type integral and the dual quermassintegral are then applied to analyzing the convergence
of a maximizing sequence, and this shows the existence of a solution to the maximization problem.

\vs{5}

For the symmetric case, the necessary and sufficient conditions for dual curvature measures are the conditions that show
the concentration of the measures on great subspheres.

\subsubsection{Subspace Mass Inequality}
For $0<q<n$, a finite Borel measure $\mu$ on $S^{n-1}$ is said to satisfy
\emph{the $q$th subspace mass inequality} if
\begin{equation}\label{smi2}
\frac{\mu(\xi\cap S^{n-1})}{\mu(S^{n-1})} <
\begin{cases}
i/q  &i < q \\
1 &i \ge q
\end{cases}
\end{equation}
for any proper subspace $\xi$ of dimension $i$ in $\rn$.
\vs{3}

When $0<q\leq 1$, the $q$th subspace mass inequality simply says that the measure
$\mu$ is not concentrated on a great subsphere.

\subsubsection{Existence of Solution to the Symmetric Dual Minkowski Problem}
{\it Let $0<q<n$ and $\mu$ be a non-zero even finite Borel measure on $S^{n-1}$.
Then there exists an origin-symmetric convex body $K$ in $\rn$ such that
$\widetilde{C}_q(K,\cdot)=\mu$ if and only if $\mu$ satisfies
the $q$th subspace mass inequality \eqref{smi2}.
}
\vs{3}

The case of $0<q\le 1$ was proved by Huang-LYZ \cite{HLYZ16acta}.
For the case of $1<q<n$, with a stronger sufficient condition than
\eqref{smi2}, the existence of solution was also proved in \cite{HLYZ16acta}.
That sufficient condition was refined as \eqref{smi2} by Zhao \cite{Zhao}
for integer $q=2, \ldots, n-1$, and by
B\"{o}r\"{o}czky,
Henk \& Pollehn \cite{BHP18} for any $1<q<n$, independently. That is,
the $q$th subspace mass inequality is a necessary condition for
the $q$th dual curvature measure to satisfy.
Zhao \cite{Zhao} also showed that
for integer $q=2, \ldots, n-1$,
the $q$th subspace mass inequality is a sufficient condition for a finite Borel measure
to be the $q$th dual curvature measure of a symmetric convex body.
B\"or\"oczky-LYZ-Zhao \cite{BLYZZ19adv} generalized Zhao's result to any real $1<q<n$.
The dual Minkowski problem for even data and $q=0$ is equivalent to
the Aleksandrov problem for even data, solved by Aleksandrov himself.
A variational proof appeared in \cite{HLYZ18jdg}.
When $q = n$, the dual Minkowski problem for even data is the logarithmic Minkowski problem for even data,
 solved in \cite{BLYZ13jams}.
Thus, the dual Minkowski problem is completely solved for the symmetric case when $0\le q\le n$.

When the given measure $\mu$ on $\sn$ is absolutely continuous with respect to the spherical
Lebesgue measure, the measure $\mu$ has no concentration on a great subsphere. Thus, the subspace mass inequality
is true. In this case, the maximization problem \eqref{dualMax} can be solved quite easily similar to solving
the classical Minkowski problem. In fact, if $\mu$ has density $f$ that is even and bounded below by a positive
constant, then for any $q>0$ the dual Minkowski problem has a symmetric convex body as a solution.
However, it is difficult to see how to solve the dual Minkowski problem for a general measure
by using approximation and the solution to the case of the measure with a density.

\vs{3}

The dual Minkowski problem for $q<0$, similar to the classical
Minkowski problem, does not require any non-trivial
measure concentration condition and was solved by Zhao \cite{Zhao2}.

\subsubsection{Existence of Solution to the Dual Minkowski Problem with $q<0$}
{\it Suppose $q < 0$ and $\mu$ is a non-zero finite Borel measure on $\sn$. Then
there exists a convex body $K$ in $\rn$ that contains the origin in its interior
such that $\wt C_q(K,\cdot)=\mu$ if and only if $\mu$ is not concentrated in any closed
hemisphere.
}

\vs{5}

\subsubsection{Uniqueness of the Dual Minkowski Problem}
The uniqueness problem of dual curvature measures was shown for $q\le 0$. The case when $q=0$ is classical and follows
from the uniqueness result of integral curvature shown by Aleksandrov. For convex bodies $K, L$ in $\rn$ that
contain the origin in their interiors, $\wt C_0(K,\cdot)=\wt C_0(L, \cdot)$
implies that $K=c L$ for $c>0$, and for $q<0$, $\wt C_q(K,\cdot) = \wt C_q(L,\cdot)$
implies that $K=L$. This was proved by Zhao \cite{Zhao2}.
The case when $q>0$  has not yet been settled.

The uniqueness together with
the variational solution to the dual Minkowski problem would imply
unproven Brunn-Minkowski type inequalities of dual quermassintegrals, for example,
\[
\wt V_q(K+L)^\frac1q \ge \wt V_q(K)^\frac1q + \wt V_q(L)^\frac1q,
\]
for origin-symmetric convex bodies $K, L$ in $\rn$ and $1<q\neq n$. This inequality
was conjectured by Lutwak in an unpublished paper shown to Yang and Zhang in 1990s.
The case of $q=n$ is the classical Brunn-Minkowski inequality and the case of $q\le 1$
can be easily shown by using the inclusion $K\wt+ L \subseteq K+L$.

\subsubsection{Regularity of Solution to the Dual Minkowski Problem and $q$-Gauss Curvature Flow}
Motivated by the dual Minkowski problem, Li-Sheng-Wang \cite{LSW20} introduced the $q$-Gauss curvature flow,
\begin{equation}\label{aGf1}
\frac{\partial X}{\partial t} (u,t) = -  |X(u,t)|^{n-q} \lkappa(u,t) f(\nu(u,t)) \nu(u,t), \quad q\in \R,
\end{equation}
where $X(\cdot,t):  \sn \to \rn$ is a parametrization of the boundaries of a family of
smooth convex bodies $K_t \in \kon$ with positive Gauss curvature $\lkappa(\cdot, t)$ and outer unit normal
 $\nu(\cdot, t)$ at $X(\cdot, t)$, and $f$ is a given positive smooth function $f$ on $\sn$.
When $q=n$, this flow equation is the anisotropic Gauss curvature flow \eqref{aGf}.

The self-similar solution of a normalized flow equation \eqref{aGf1} satisfies
the Monge-Amp\`ere type equation \eqref{dualp}
of the dual Minkowski problem.  Li-Sheng-Wang \cite{LSW20} gave solutions to the $q$-Gauss curvature flow
which yield regularity results to the dual Minkowski problem:
\vs{3}

{\it Suppose a Borel measure $\mu$ on $\sn$ has a density $f$ that is positive and $C^\infty$.
When $q<0$, the solution $K$ to the dual Minkowski problem is $C^\infty$ and uniformly convex
(i.e., of positive principal curvatures). When $q>0$ and $f$ is even,  the solution $K$
is $C^\infty$ and uniformly convex.}
\vs{3}

When $n=2$ and $q>0$, similar regularity results, but without assuming evenness, were shown by
Chen-Li \cite{CL18}. Again in the plane when $q\ge 2$ is an even integer, $C^2$ regularity
of the solution $K$ was proved by Huang-Jiang \cite{HJ19} under a weaker assumption on $f$.
When $f=1$ and the Gauss curvature is replaced by an intermediate curvature in \eqref{aGf1},
convergence results were proved by Li-Sheng-Wang \cite{LSW20jga}.

\subsubsection{Affine Dual Minkowski Problem} In the formula \eqref{dqxi} of dual quermassintergals
of a convex body $K$ in $\rn$, if the intersection volume is replaced by its $n$th power, one gets Lutwak's
{\it affine dual quermassintegrals} of $K$ (up to some normalization),
\[
\wt V_i^a(K) = \int_{G_{n,i}} V_i(K\cap\xi)^n\, d\xi, \ \ \ i=1,\ldots, n-1.
\]
It is an important question to ask if one can construct geometric measures as differentials of the affine
dual quermassintegrals.
 Recently, this was answered by Cai-Leng-Wu-Xi \cite{CLWX1} who proved the following differential formula,
\[
\frac{d\wt V_i^a(\blb h_t\brb)}{dt}\Big|_{t=0} = n\int_{\sn} f(v)\, d\wt C_i^a(K,v),
\]
where $\blb h_t \brb$ is the Wulff shape of the variations $h_t$ defined by $\log h_t = \log h_K + tf$.
The measure $\wt C_i^a(K,\cdot)$ is called the $i$th {\it affine dual curvature measure} of $K$.
The case of $i=n-1$ is closely related to the intersection body of $K$.
The construction requires to develop new tools of Radon transforms. In \cite{CLWX1},
the associated {\it affine dual Minkowski problem} was solved for origin-symmetric convex bodies under
the strict mass concentration inequality \eqref{ssc}.

\section{Chord Minkowski Problems from Integral Geometry}\

While quermassintegrals and dual quermassintegrals are the two families of
fundamental geometric invariants in
the Brunn-Minkowski theory and its dual theory,
 a third family of geometric invariants of convex bodies called  {\it chord integrals}
 arises from integral geometry and geometric probability. These are
moments of lengths of intersections of random lines
with convex bodies.  They also give the second moment of areas of intersections of
random planes with convex bodies \cite{Z99tams}.
Chord integrals are closely related to $X$-rays in geometric tomography \cite{G06book} and
Riesz potentials of characteristic functions of convex bodies.

In the recent paper \cite{LXYZ21}, the authors introduce the geometric measures derived from
chord integrals, called {\it chord measures}. They are the second family of translation invariant
geometric measures of convex bodies besides Aleksandrov-Fenchel-Jessen's area measures.
We describe Minkowski problems associated with these measures.

\subsection{Chord Integrals and Measures}

\subsubsection{Chord Integrals}
Let $\mathscr L^n$ be the affine Grassmannian of lines in $\rn$.
Lines in $\mathscr L^n$ are determined by their rotations (directions) and translations.
The Haar measure on $\mathscr{L}^n$ is normalized to be a probability measure when restricted to rotations of lines
and to be $(n-1)$-dimensional
Lebesgue measure when restricted to translations of lines.
The {\it chord integral} $I_q(K)$ of convex body $K$ in $\rn$ is defined by
\begin{equation}\label{ci}
I_q(K) = \int_{\mathscr{L}^n}
|K\cap \ell|^q\, d\ell, \ \ \ \ \ q\ge 0,
\end{equation}
where $|K\cap \ell|$ denotes the length of the chord $K\cap \ell$, $\ell \in \mathscr L^n$, and
$d\ell$ is the normalized Haar measure on $\mathscr{L}^n$.

Volume and surface area are chord integrals. Indeed, there are integral formulas,
\begin{equation*}
I_1(K) =  V(K), \ \
I_0(K)= \frac{\omega_{n-1}}{n\omega_n} S(K), \ \
I_{n+1}(K) = \frac{n+1} {\omega_n} V(K)^2,
\end{equation*}
which are Crofton's volume formula, Cauchy's integral formula for surface area, and
the Poincar\'e-Hadwiger formula, respectively (see \cite{Ren, San}).

The relationship between chord integrals and the Riesz potentials of
characteristic functions of convex bodies is given by
\begin{equation}\label{6.1}
I_q(K) = \frac{q(q-1)}{n\omega_n} \int_{\rn}\int_{\rn}
\frac{{\pmb 1}_K(x) {\pmb 1}_K(y)}{|x-y|^{n-q+1}}\,
dxdy, \ \ \ q>1,
\end{equation}
(see \cite{Ren, San}). There is also the formula,
\begin{equation}\label{6.2}
I_q(K) = \frac{q(1-q)}{n\omega_n} \int_{\rn}\int_{\rn}
\frac{{\pmb 1}_K(x) {\pmb 1}_{K^c}(y)}{|x-y|^{n-q+1}}\,
dxdy, \ \ \ 0<q<1,
\end{equation}
where $K^c$ is the complement of $K$. Formula \eqref{6.2} was shown by Qin \cite{Q22}, and implicitly by
Ludwig \cite{Lud14}.

The integral on the right-hand side of \eqref{6.2} is called
a fractional integral of $K$. Formula \eqref{6.2} and Cauchy's integral formula for $I_0(K)$ show that
the fractional integral tends to the surface area of $K$, up to a constant factor, as $q\to 0^+$,
which is a fact proved by Bourgain-Brezis-Mironescu \cite{BBM01} for any
bounded Borel set in $\rn$. A generalization of the fractional integral, called anisotropic fractional integral,
  was given by Ludwig \cite{Lud14}
when the Euclidean norm is replaced by a non-Euclidean norm. Results in \cite{BBM01} were extended in \cite{Lud14}.

The relationship between chord integrals and dual quermassintegrals is given by
\begin{equation*}
I_q(K) = \frac{q}{\omega_n} \int_K \wt V_{q-1}(K,z) \, dz, \ \  q>0,
\end{equation*}
(see \cite{GZ98, Z99tams}), where $\wt V_{q-1}(K,z)$ is the dual quermassintegral of $K$ with respect to
a point $z$ in $K$,
\[
\wt V_{q-1}(K, z) = \frac1n \int_{S_z^+} \rho_{K-z}^{q-1}(u)\, du, \ \ \  S_z^+=\{u\in\sn : \rho\lsub{K-z}(u)>0\}.
\]
As a function of $z\in K$, $\wt V_{q-1}(K, z)$ is continuous when $q>1$, but only lower semi-continuous when
$0<q\le 1$. When $z\in \partial K$ and $0<q<1$, $\wt V_{q-1}(K, z)$ may be infinite.

\subsubsection{Chord Measures}
In the paper \cite{LXYZ21}, differential formulas of chord integrals were derived. Then a new family of geometric measures
was defined.
For each $q>0$, the {\it $q$th chord measure} $F_q(K,\cdot)$
of a convex body $K$ in $\rn$ is defined by
\begin{equation}\label{cm}
F_q(K,\eta) = \frac{2q}{\omega_n} \int_{\nu_K^{-1}(\eta)} \wt V_{q-1}(K, z)\, dz,
\ \  \text{ for any Borel set $\eta\subset \sn$.}
\end{equation}

Each $F_q(K,\cdot)$ is a finite Borel measure on $\sn$ and is the differential of $I_q$ at $K$; i.e.,
\begin{equation}\label{cmd}
\frac{d}{dt}\Big|_{t=0^+} I_q(K+tL) = \int_{\sn} h_L(v)\, dF_q(K,v),
\end{equation}
for each convex body $L$. The case of $q=1$ corresponds to the classical formula \eqref{b2.2}
and the measure $F_1(K,\cdot)$ is the classical surface area measure $S(K,\cdot)$, i.e.,
$
F_1(K,\cdot) = S(K, \cdot).
$

The extremal case $q\to 0^+$ is very interesting. It exhibits a surprising connection between dual
quermassintegrals and mean curvature: For each sufficiently smooth convex body $K$
and each $z\in \partial K$,
\[
\lim_{q \to 0^+} q \wt V_{q-1}(K,z) = \frac{(n-1)\omega_{n-1}}{2n} H(z),
\]
where $H(z)$ is the mean curvature of $\partial K$ at $z$. This
 implies that for a sufficiently smooth convex body
the area measure $S_{n-2}(K,\cdot)$ is a limiting case of chord measures,
\begin{equation}\label{0.13}
\lim_{q\to 0^+} F_q(K, \eta) = \frac{(n-1)\omega_{n-1}}{n\omega_n} S_{n-2}(K, \eta),
\end{equation}
for each Borel set $\eta \subset \sn$.
Thus, one can define
\[
F_0(K,\cdot)= \frac{(n-1)\omega_{n-1}}{n\omega_n} S_{n-2}(K, \cdot).
\]
However,  it is not clear what the limit in \eqref{0.13} is for non-smooth convex bodies.

The chord measures $F_q(K,\cdot)$, $q\ge 0$, are translation invariant. This implies that
\begin{equation}\label{c14.1}
\int_{\sn} v\, dF_q(K,v) = 0.
\end{equation}

The relation between surface area measure and cone-volume measure can be extended to
chord measures and cone-chord measures. The {\it $q$th cone-chord measure $G_q(K,\cdot)$}  is defined by
\begin{equation}\label{0.10}
dG_q(K,\cdot) =\frac1{n+q-1} h_K\, dF_q(K,\cdot).
\end{equation}
When $q=1$, $G_1(K,\cdot)$ is  the cone volume measure $V(K,\cdot)$. In this case, \eqref{0.10} corresponds
to \eqref{5.1}.
The chord integral, the total mass of cone chord measure, and the chord measure are related by
\begin{equation*}
I_q(K)= G_q(K,\sn) = \frac1{n+q-1} \int_{\sn} h_K(v)\, dF_q(K,v).
\end{equation*}

\subsection{Minkowski Problems of Chord Measures}

The Minkowski problem of prescribing chord measures was posed in \cite{LXYZ21}:
\subsubsection{Chord Minkowski Problem}
{\it
Let $q\ge0$. If $\mu$ is a finite Borel measure on $\sn$,
what are necessary and sufficient conditions on $\mu$ to guarantee the existence of a convex body
$K\subset\rn$ that solves the equation,
\begin{equation}\label{0.14}
F_q(K,\cdot) = \mu\, ?
\end{equation}
}

When the given finite measure $\mu$ has a density and $q>0$,
equation \eqref{0.14} becomes a new Monge-Amp\`ere type equation on $\sn$,
\begin{equation}\label{0.15}
\det\big(\nabla_{ij}h(v) + h(v)\delta_{ij}) =\frac{f(v)}{\wt V_{q-1}(\wul{h},\nabla h(v))},
\end{equation}
where $\wt V_{q-1}(\wul{h}, \nabla h(v))$ denotes the dual quermassintegral of the
Wulff shape $\wul{h}$ with respect to the point $\nabla h(v) \in \rn$, which is a function of $v$.
Note that the dual quermassintegral is expressed by the radial function and depends implicitly on
$h$ because $h$ corresponds to the support function. Thus, $\wt V_{q-1}(\wul{h}, \nabla h(v))$
is a function of $h$ and $\nabla h(v)$
only in the implicit sense and is not locally defined. Also,
since $\nabla h(v)$ is a boundary point of $\wul{h}$, $\wt V_{q-1}(\wul{h}, \nabla h(v))$
may be infinite when $0<q<1$.
\vs{3}

This chord Minkowski problem was solved in \cite{LXYZ21} when $q>0$:

\subsubsection{Solution to the Chord Minkowski Problem}
{\it Let $q>0$.
If $\mu$ is a finite Borel measure on $\sn$, then there exists a convex body $K$ in $\rn$ satisfying
$
F_q(K,\cdot) = \mu
$
if and only if $\mu$ is not concentrated on a
closed hemisphere and
\[
\int_{\sn} v\, d\mu(v) =0.
\]}

This solution extends that of the classical Minkowski problem in \S \ref{4.1}. The proof
is variational, similar to the classical case. It involves solving
a maximization problem for the functional $\Phi(h)$ over $C^+(\sn)$,
\[
\Phi(h) = I_q([h])^\frac1{n+q-1} \Big( \int_{\sn} h\, d\mu\Big)^{-1}.
\]
This is shown by using the variational formula \eqref{cmd}.

The log-Minkowski problem of cone-volume measure in \S \ref{5.2} was extended to the following:

\subsubsection{The Chord Log-Minkowski Problem} {\it Let $q\ge0$.
If $\mu$ is a finite Borel measure on $\sn$,
what are necessary and sufficient conditions on $\mu$ to guarantee the existence of a convex body
$K$ in $\kon$ that solves the equation,
\begin{equation}\label{0.16}
G_q(K,\cdot) = \mu\, ?
\end{equation}
}

When $q=1$, the chord log-Minkowski problem is just the log-Minkowski problem,
When $q=0$, it is a log-Christoffel-Minkowski problem, which may well be even more difficult
to solve than the still open Christoffel-Minkowski problem.

The partial differential equation associated with \eqref{0.16} is the following Monge-Amp\`ere equation on $\sn$,
\begin{equation}\label{0.17}
\det\big(\nabla_{ij}h + h\delta_{ij})=\frac{f}{h\, \wt V_{q-1}(\wul{h}, \nabla h)}.
\end{equation}

Similar to the situation with the log-Minkowski problem versus the classical Minkowski problem,
the chord log-Minkowski problem is a harder problem than the chord Minkowski problem.
The following solution is obtained in \cite{LXYZ21} and \cite{Q22} for the symmetric case.

\subsubsection{Solution to the Symmetric Chord Log-Minkowski Problem}
{\it Let $\mu$ be an even finite Borel measure on $S^{n-1}$ and $0< q\le n+1$. Then
there exists an origin-symmetric convex body $K$ in $\rn$ that solves the equation,
\[
G_q(K,\cdot) = \mu,
\]
if $\mu$ satisfies the subspace mass inequality:
\begin{equation}
\frac{\mu(\xi_k\cap S^{n-1})}{|\mu|}  <  \frac{k+\min\{k,q-1\}}{n+q-1},
\end{equation}
for each $k$-dimensional subspace $\xi_k\subset\rn$ and each $k=1, \ldots, n-1$.}

The above solution of the chord log-Minkowski problem for the case where $q=1$ is the solution to the
log-Minkowski problem obtained in \cite{BLYZ13jams}. The case of $1<q\le n+1$ was proved in \cite{LXYZ21}.
By using a different approach via formula \eqref{6.2}, Qin \cite{Q22} proved the case of $0<q<1$.
For more recent studies on chord Minkowski problems, see \cite{GXZ24, HQ24, Li24, Q24}.

\subsubsection{Problems}

Some important open problems related to chord measures need to be studied. The uniqueness problem of chord measures
is a major problem:
\vs{2}

{\it If $F_q(K,\cdot) = F_q(L, \cdot)$ for a positive $q$, is $K$ a translation of $L$?}
\vs{2}

\noindent
The cases of $q=1, n+1$ are known, which are the uniqueness result of surface area measure.
An affirmative answer to this uniqueness question together with the variational solution to
the chord Minkowski problem
would imply Brunn-Minkowski type geometric inequalities for chord integrals,
\[
I_q(K+L)^\frac1{n+q-1} \ge I_q(K)^\frac1{n+q-1} + I_q(L)^\frac1{n+q-1}, \ \  q>0,
\]
with equality if and only if $K$ and $L$ are homothetic.
It is the classical Brunn-Minkowski inequality when  $q=1, n+1$. The case of $q=0$ is the classical
Brunn-Minkowski inequality of surface area which is a consequence of the Aleksandrov-Fenchel inequality
(see \cite{S14}). Other cases have not been studied.

Regularity in the chord Minkowski problems is also highly interesting. It is natural to solve
the chord Minkowski problems by using the continuity method or an anistropic Gauss curvature flow so that
smooth solutions are obtained. The nonlocal quantity $\wt V_{q-1}(\wul{h}, \nabla h)$ in the PDE \eqref{0.15}
may create some obstacles. Recently, Hu-Huang-Lu \cite{HHL24} obtained regularity results for Riesz potential
and applied them to the regularity of the chord log-Minkowski problem.

While chord integrals come from intersections of convex bodies with straight lines in $\rn$, it is natural
to consider intersections of convex bodies with flat planes.

\section{$L_p$ Minkowski Problems in the $L_p$ Brunn-Minkowski Theory}\

The basic algebraic operation in the Brunn-Minkowski theory is the Minkowski addition of convex bodies
(see \eqref{Madd}). By the motivation of interactions between convex geometry and functional analysis,
extending the Minkowski addition to the $L_p$ Minkowski addition leads to
the $L_p$ Brunn-Minkowski theory which generalizes the classical Brunn-Minkowski theory.
The  $L_p$ Brunn-Minkowski theory has found
many connections and applications in functional analysis and other subjects.
The $L_p$ Minkowski addition was first studied by Firey \cite{F62} for $p>1$ in early 1960s.
The $L_p$ Brunn-Minkowski theory gained its status only after the celebrated work of
Lutwak \cite{L93jdg, L96adv} in 1990s.
Subsequent studies on affine isoperimetric inequalities \cite{LYZ00jdg, LYZ02jdg, CG02, HS09jdg, HS09jfa},
valuations \cite{Lud03, Lud04, Lud06, Lud09, LR10}, affine surface areas \cite{MW00, SW04, W07, WY08adv, WY10mathann},
and Minkowski problems \cite{CW06, HLYZ05dcg, BLYZ13jams, LO95jdg, LYZ04tams, LYZ06imrn}
have made it a significant theory with rich and deep results.

\vs{3}
The critical concept introduced by Lutwak is the
$L_p$ surface area measure that extends the classical surface area measure of Aleksandrov, Fenchel and Jessen.
The problem of prescribing $L_p$ surface area measure is called the $L_p$ Minkowski problem.

\subsubsection{Basics of $L_p$ Brunn-Minkowski Theory}

For convex bodies $K$ and $L$ in $\mathcal K_o^n$ and $p\ge1$, $(h_K^p+h_L^p)^\frac1p$ is a positive convex
function of homogeneous of degree 1, and hence it is the support function of a convex body in $\mathcal K_o^n$,
denoted by $K\psum L$ and called the {\it $L_p$ sum} of $K$ and $L$. Thus,
\[
h_{K\spsum L} = (h_K^p+h_L^p)^\frac1p.
\]
The $L_1$ sum is the usual Minkowski sum.

The $L_p$ sum can be extended to the case of $p<1$ by using Wulff shape.
For $p, s, t\in \R$, suppose
$sh_K^p + t h_L^p$ is positive. The $L_p$ sum  $s\thin\cdot K \psum t\thin\cdot L$
 is defined as the Wulff shape generated by
the function $(sh_K^p + th_L^p)^\frac1p$ when $p\neq 0$ and by $h_K^s h_L^t$ when $p=0$, i.e.,
\[
s\thin\cdot K\psum t\thin\cdot L = \text{\bf \small [}(sh_K^p + th_L^p)^\frac1p\text{\bf\small ]}, \quad
s\thin\cdot K\qsum{0} t\thin\cdot L = \text{\bf \small [} h_K^s h_L^t \text{\bf \small ]}.
\]
Note that the geometric mean $K^{1-t}L^t$ is the $L_0$ sum $(1-t)\thin\cdot K\qsum{0} t\thin\cdot L$.

We remark that the $L_p$ sum $s\thin\cdot K\psum t\thin\cdot L$ is defined for $K, L\in\kon$, i.e.,
when the origin is inside the interiors of $K$ and $L$. But it is still well defined
even if the origin is on the boundary of $K$ or $L$ as long as $sh_K^p + t h_L^p$ is positive.

As demonstrated by Lutwak \cite{L93jdg}, area measures can be extended to
$L_p$ area measures. For a convex body $K$ in $\kon$ and $p\in \R$, its {\it $L_p$ $i$th area measure}
is the finite Borel measure on $\sn$ defined by
\begin{equation}\label{pam}
dS_{p,i}(K, \cdot) = h_K^{1-p} dS_i(K,\cdot), \ \ \ i=0, 1,\ldots, n-1.
\end{equation}
The case of $i=n-1$, $S_{p, n-1}(K,\cdot)$, is called the
{\it $L_p$ surface area measure} and is often denoted by $S_p(K,\cdot)$ when there is no confusion.
$L_p$ area measures arise naturally as the right-side differentials of quermassintegrals with respect to the $L_p$ sum,
\begin{equation}\label{pwd}
\frac d{dt} W_{n-i-1}(K\psum t\thin\cdot L)\Big|_{t=0^+} = \frac{i+1}{np} \int_{\sn} h_L^p(v)\, dS_{p,i}(K,v), \ \ \
i=0, 1,\ldots, n-1,
\end{equation}
for $p\neq 0$. There is a similar formula for the case of $p=0$.

A more geometric definition of the $L_p$ surface area measure is
\begin{equation}\label{p1.1}
S_p(K,\eta) = \int_{x\in\nu_K^{-1}(\eta)} (x\cdot \nuk(x))^{1-p}\, dx, \ \ \ \text{for any Borel set } \eta \subset \sn,
\end{equation}
where for $p>1$, $K$ contains the origin in its interior, for $p<1$, $K$ contains the origin (which may be on
the boundary of $K$), and for $p=1$, $K$ is any convex body.
The $L_1$ surface area measure is obviously the surface area measure $S(K,\cdot)$.
When $p=0$, the $L_0$ surface area measure is $n$ times the cone-volume measure $V(K,\cdot)$.

The $L_p$ surface area measure is not only the right-side differential of the volume functional at a convex body
with respect to the $L_p$ Minkowski sum, but also the left-side differential.
In fact, the differential formula \eqref{b2.2} implies that
\begin{equation}\label{p1.2}
\frac d{dt} V(K\psum t\thin\cdot L)\Big|_{t=0} = \frac1p \int_{\sn} h_L^p(v)\, dS_p(K,v).
\end{equation}
This formula is critical for solving the $L_p$ Minkowski problem.
The case of $p=0$ corresponds to the geometric mean and cone-volume measure, discussed in the previous section.

The integral above defines the $L_p$ mixed volume $V_p(K,L)$ of convex bodies $K$ and $L$,
\[
V_p(K,L) =\frac1n \int_{\sn} h_L^p(v) \, dS_p(K,v).
\]
When $K=L$, $V_p(K,K)=V(K)$. When $L$ is the unit ball $B^n$, $V_p(K, B^n)$, denoted by $S_p(K)$, is
called the $L_p$ surface area of $K$. There is the $L_p$ mixed volume inequality,
\begin{equation}\label{p1.3}
V_p(K,L)^n \ge V(K)^{n-p} V(L)^p,  \ \ \ p>1,
\end{equation}
with equality if and only if $K$ is a dilate of $L$.

The $L_p$ mixed volume inequality follows from Minkowski's first mixed volume inequality \eqref{b1.7}
and H\"older's inequality.
The equality conditions of \eqref{p1.3} implies the following uniqueness result:
\vs{3}

{\it Let $K$ and $L$ be convex bodies in $\rn$ that contain the origin in their interiors and $p>1$.
If $S_p(K,\cdot) = S_p(L,\cdot)$, then $K=L$ when $p\neq n$, and $K$ is a dilation of $L$ when $p=n$.}
\vs{3}

When $p<1$, $L_p$ surface area measure is not in general unique. It was conjectured that
for $0<p<1$ the $L_p$ mixed volume inequality still holds for origin-symmetric convex bodies. This would
be true if the logarithmic Minkowski inequality \eqref{lMi} is true, see B\"or\"oczky-LYZ \cite{BLYZ12adv}.
See Chen-Li-Zhu \cite{CLZ17, CLZ19}, Jian-Lu-Wang \cite{JLW15} and Li-Liu-Lu \cite{LLL21} for results about
the case of $p < 1$.
Recently, Milman \cite{M22jdg} gave a new strong non-uniqueness result for the even $L_p$-Minkowski
problem when $-n < p < 0$.

Suppose that $p$ is less than and close to 1.
Kolesnikov-Milman \cite{KM18} showed that $K=L$ if $K$ and $L$ have positive curvature and are close
in $C^2$ norm.  Chen-Huang-Li-Liu \cite{CHLL20} strengthened this local result as a global one by using
the Schauder estimates in PDEs to show that $K=L$ if $K$ and $L$ have positive curvature and are of class $C^{2,\alpha}$.

\vs{3}

Geometric measures in the Brunn-Minkowski theory and the dual Brunn-Minkowski theory have
$L_p$ analogues. Therefore, Minkowski problems discussed before have their $L_p$ versions.
Here, we describe the $L_p$ Minkowski problem, the $L_p$ Aleksandrov problem, the $L_p$
dual Minkowski problem, and the $L_p$ chord Minkowski problem. The $L_p$ Minkowski problem
has been studied most.

\subsubsection{$L_p$ Minkowski Problem}

The classical Minkowski problem of prescribing surface area measure was generalized to
the $L_p$ Minkowski problem of prescribing $L_p$ surface area measure by Lutwak \cite{L93jdg}:
\vs{3}

{\it Given a finite Borel measure $\mu$ on $\sn$ and $p\in \R$,
what are the necessary and sufficient conditions for the existence of
a convex body $K\in \kon$ satisfying
\begin{equation}\label{p1.4}
S_p(K,\cdot) = \mu\,?
\end{equation}
}

The associated partial differential equation of the $L_p$ Minkowski problem is the Monge-Amp\`ere equation
\begin{equation}\label{p1.6}
\det(\nabla_{ij} h + h \delta_{ij}) = h^{p-1} f.
\end{equation}

The problem contains three special cases which have important geometric meanings.
When $p=1$, it is the classical Minkowski problem discussed in Section 5.2.
When $p=0$, it is the logarithmic Minkowski problem described in Section 6.2. When $p=-n$, it is called the
{\it centro-affine Minkowski problem}. In this case, the absolutely continuous part of the $L_{-n}$
surface area measure with respect to the spherical Lebesgue measure is $h_{K}^{1+n}f_K$ whose $(n+1)$-th root is the affine
distance in affine differential geometry (see \cite{NS94}, pp. 62-63). Prescribing the affine distance is a long standing unsolved problem,
which is then a special case of the $L_p$ Minkowski problem and the partial differential equation is
\begin{equation}\label{p1.7}
\det(\nabla_{ij} h + h \delta_{ij}) = \frac{f}{h^{n+1}}.
\end{equation}

The $L_p$ Minkowski problem has been solved for $p\ge 1$ and also solved for $0\le p<1$ in the symmetric case.

\fr{Existence of Solutions to the Symmetric $L_p$ Minkowski Problem}
{\it Let $\mu$ be a finite even Borel measure that is not concentrated on a great subsphere. Then
there exists an origin-symmetric convex body $K$ such that $S_p(K,\cdot)=\mu$ when $0<p\neq n$
and $S_n(K,\cdot) = V(K) \mu$.}

\vs{5}

This was proved for $p>1$ by Lutwak \cite{L93jdg}  and for $0<p<1$ by Haberl-LYZ \cite{HLYZ10}.
A different proof and the volume normalized formulation were given by LYZ \cite{LYZ04tams}.

\fr{Existence of Solutions to the Asymmetric $L_p$ Minkowski Problem}
{\it Let $\mu$ be a finite Borel measure on $\sn$ that is not concentrated
in a closed hemisphere. Then, for each $p>1$, there exists a unique
convex body $K$ that contains the origin such that
\begin{equation}\label{p1.5}
dS(K, \cdot)=h_K^{p-1} d\mu, \ \ 1<p\neq n, \ \ \text{and } \ S_n(K,\cdot)=V(K)\mu.
\end{equation}
Moreover, $K$ contains the origin in its interior if either $p\ge n$ or $\mu$ is discrete.}
\vs{5}

This result was first proved by Chou-Wang \cite{CW06}, while the fact that the solution contains the origin in its
interior if either $p\ge n$ or $\mu$ is discrete
was shown by Hug-LYZ \cite{HLYZ05dcg}.
Equation \eqref{p1.5} is a volume normalized formulation of \eqref{p1.4},
first given by LYZ \cite{LYZ04tams}.
When $p\neq n$, the volume factor in \eqref{p1.5} can always be absorbed by a dilation.
When the origin is contained in the interior of $K$, equations \eqref{p1.4} and \eqref{p1.5}
are equivalent. However, there exist examples which show that the origin can be on the
boundary of $K$ when $1<p<n$, see \cite{HLYZ05dcg}. In this case, \eqref{p1.5} does not imply \eqref{p1.4}.

Chou-Wang \cite{CW06} proved the result for measures with smooth densities first by using the continuity method,
which also gives regularity of solutions. Then the general case of measures was shown by approximation.
Hug-LYZ \cite{HLYZ05dcg} used the variational method to prove the result for discrete measures first and followed
by an approximation argument. The optimization problem associated with the $L_p$ Minkowski problem is
\[
\inf_{h\in C^+(\sn)}\Big\{V(\wul{h})^{-\frac1n} \Big(\int_{\sn} h^p\, d\mu\Big)^\frac1p\Big\}.
\]
If this minimization problem has a positive solution, then its Wulff shape gives a solution to the
$L_p$ Minkowski problem. This can easily be seen by a proof similar to the classical case.
When $\mu$ is discrete or $p\ge n$, the minimization problem has a positive solution when $p>1$. However,
when $1<p<n$, the minimization problem, in general,  has only a nonnegative solution whose Wulff shape may not
contain the origin in its interior even if the given measure has a positive smooth density.
Thus, additional techniques are used to finish the proof.

\fr{Existence of Solutions When $0<p<1$}

{\it If $0<p<1$, and $\mu$ is a finite Borel measure
on $\sn$ not concentrated on a great subsphere, then there exists a convex body $K$ containing the origin
so that $S_p(K,\cdot)=\mu$.}

\vs{3}
The polytope case was proved by Zhu \cite{Zhu1} and then the general case was shown by
Chen-Li-Zhu \cite{CLZ17} by an approximation argument.

Note that for $p<1$ the $L_p$ surface area measure $S_p(K,\cdot)$ may concentrate on a great subsphere when
the origin is on the boundary of $K$. Indeed, if $K$ is a cone whose vertex is the origin, then $S_p(K,\cdot)$
is a single point mass. It was conjectured \cite{BBCY}:
\vs{3}

{\it Let $0<p<1$, and let $\mu$ be a non-trivial finite Borel measure on $\sn$. Then there exists a convex body
$K$ containing the origin so that $S_p(K,\cdot)=\mu$
if and only if $\mu$ does not concentrate on a pair of antipodal points.}
\vs{3}

The planar case $n = 2$ was proved independently by B\"or\"oczky-Trinh \cite{BT17} and
Chen-Li-Zhu \cite{CLZ17}. Partial results on higher dimensions were obtained in \cite{CLZ17, BBCY}.

\fr{Existence of Solutions When $p<0$} When $p<0$, the $L_p$ Minkowski problem has been studied
when the given Borel measure on $\sn$ is either discrete or absolutely continuous. For the discrete case,
there is the following solution:
\vs{3}

{\it If $p<0$ and $\mu$ is a discrete measure on $\sn$, then there exists a polytope $P$ in $\rn$ whose outer unit normals
are in general position so that $S_{p}(P, \cdot) =\mu$ if and only if the support of $\mu$ is in general position
and not concentrated on a closed hemisphere.}
\vs{3}

When $p=-n$, this result gives a solution to the centro-affine Minkowski problem for polytopes whose outer unit normals
are in general position, proved by Zhu \cite{Zhu3}. Then he \cite{Zhu4} generalized his result to any $p<0$
and to a larger class of polytopes.  Jian-Lu-Zhu \cite{JLZ16} gave a solution to the centro-affine Minkowski problem
for the class of convex bodies symmetric about coordinates hyperplanes.
Li \cite{LiQ19} showed that there may exist infinitely many solutions which are not affine equivalent.

\vs{3}

For the absolutely continuous case, there is the following result:
\vs{3}

{\it If $-n<p<0$, and $\mu$ is a non-trivial Borel measure on $\sn$ with a density function in $L_{\frac n{n+p}}(\sn)$,
then there exists a convex body $K$ containing the origin so that $S_p(K,\cdot) =\mu$.}
\vs{3}

The result was proved by Chow-Wang \cite{CW06} when the density function of $\mu$ is bounded and positive.
The planar case was proved by Chen \cite{C06} when the density function $\mu$ is continuous and even.
The general case was shown by Bianchi-B\"or\"oczky-Colesanti-Yang \cite{BBCY} by an approximation argument.
\vs{5}

For the smooth symmetric case, Bryan-Ivaki-Scheuer \cite{BIS19} gave a unified approach for all $p>-n$
by using curvature flow.

For $p\le 0$, the $L_p$ Minkowski problem has multiple solutions for a given measure.
This was proved by Jian-Lu-Wang \cite{JLW15} for $p \in (-n, 0)$, by He-Li-Wang \cite{HLW16}
for $p \in (-\infty, -n + 1)$, and by Chen-Li-Zhu \cite{CLZ19} for $p=0$.

For $p\le -2$, Do-Zhu \cite{DZ12} showed that the $L_p$ Minkowski problem in $\R^2$ has a solution when
the density function $\mu$ is continuous and even for each coordinate.
For $p<-n$, Guang-Li-Wang \cite{GLW23a} proved that the $L_p$ Minkowski problem in $\rn$
has a $C^{3,\alpha}$, $0<\alpha<1$, solution if $\mu$ has a positive $C^{1,1}$ density function.

\vs{3}
There are applications to Sobolev-type inequalities of functions,
in particular, the affine Sobolev inequalities \cite{LYZ02jdg, Z99jdg}.
Geometric measures are defined for functions, and convex bodies are constructed from such geometric measures
by using solutions of the $L_p$ Minkowski problem.
See \cite{CLYZ, HS09jfa, HSX12, LXZ11,LYZ02jdg, LYZ06imrn,X07,Z99jdg}.

\vs{5}

While the $L_p$ Minkowski problem has been extensively studied, the following
{\it $L_p$ Christoffel-Minkowski problem} needs more attention.
\vs{3}

{\it Given a finite Borel measure $\mu$ on $\sn$, a $p\in\R$,
and an $i\in \{1, \ldots, n-1\}$,
what are the necessary and sufficient conditions for the existence of
a convex body $K\in \kon$ satisfying
\begin{equation}\label{p1.8}
S_{p,i}(K,\cdot) = \mu\,?
\end{equation}
}

When $p\ge i+1$ and $\mu$ has a smooth density which satisfies some convexity conditions,
existence results were proved by Hu-Ma-Sheng \cite{HMS04}, and the even case of $1<p<i+1$ was studied
by Guan-Xia \cite{GX18}. The uniqueness of solutions was proved by Lutwak \cite{L93jdg} for $p>1$ and
for $1-i<p<1$ was partially proved by Chen \cite{CL20}.

\subsubsection{$L_p$ Aleksandrov Problem}

In \cite{HLYZ18jdg}, the $L_p$ integral curvature was defined and the $L_p$ Aleksandrov problem was studied.

\fr{$L_p$ Curvature Measures} For convex body $K$ in $\mathcal K^n_o$, the  {\it $L_p$ $i$th curvature measure}
$C_{p,i}(K,\cdot)$ is defined by
\[
dC_{p,i}(K,\cdot) = \rho_K^p dC_i(K,\cdot).
\]
The case of $i=0$, $C_{p,0}(K,\cdot)$ is called the {\it $L_p$ integral curvature} of $K$,
which can also be defined by
\[
C_{p,0}(K, \omega) = \int_{\alpha_K(\omega)} \rho_K^p\big(\alpha_K^{-1}(v)\big)\, dv.
\]

  It can be viewed
as the differential of the entropy of convex body with respect to $L_p$ sum,
\begin{equation}\label{pA1}
\frac d{dt} E(K^*\psum t\thin\cdot L^*)\Big|_{t=0} =\frac1p \int_{\sn} \rho\lsub{L}(u)^{-p}\, dC_{p,0}(K,u), \quad p\neq0,
\end{equation}
and for $p=0$,
\begin{equation}\label{pA2}
\frac d{dt} E(K^*\qsum{0} t\thin\cdot L^*)\Big|_{t=0} = - \int_{\sn} \log \rho\lsub{L}(u)\, dC_{0}(K,u), \quad p\neq0.
\end{equation}
These formulas are consequences of the differential formula \eqref{cm1} together with the formula
$\bla \rho\bra^* = \blb \rho^{-1}\brb$.

\fr{$L_p$ Aleksandrov-Minkowski Problem}
{\it Given a finite Borel measure $\mu$ on $\sn$, a $p\in\R$,
and an $i\in \{0, 1, \ldots, n-1\}$,
what are the necessary and sufficient conditions for the existence of
a convex body $K\in \kon$ satisfying
\begin{equation}\label{pA3}
C_{p,i}(K,\cdot) = \mu\,?
\end{equation}
}
The case of $i=0$ is called the
{\it $L_p$ Aleksandrov Problem}.

When the given measure $\mu$ has a density $f$, the partial differential equation associated with the
$L_p$ Aleksandrov problem can be written as the following Monge-Amp\`ere equation on $\sn$,
\begin{equation}\label{pAp}
\det(\nabla_{ij}h + h\delta_{ij}) = h^{p-1} |\nabla h|^n f.
\end{equation}

A solution of the equation would be the support function $h_{K^*}$ of the polar body of $K$.

The $L_0$ Aleksandrov problem is just the classical Aleksandrov problem described in Section \ref{AP}.
Huang-LYZ \cite{HLYZ18jdg} solved the $L_p$ Aleksandrov problem when $p>0$.

\fr{Existence of a Solution to the $L_p$ Aleksandrov Problem}
{\it Let $\mu$ be a finite Borel measure on $\sn$ and $p>0$. Then
there exists a convex body $K\in \kon$  so that $\mu$ is the $L_p$-integral curvature
of $K$ if and only if $\mu$ is not concentrated in any closed hemisphere.
}

\vs{5}

When $p<0$, the symmetric case was studied in \cite{HLYZ18jdg} and the following result was proved.

\fr{Existence of a Solution to the Symmetric $L_p$ Aleksandrov Problem}
{\it Let $\mu$ be a non-zero finite even measure on $\sn$.
Then there exists an origin-symmetric convex body $K$ in $\rn$ so that $\mu$ is the $L_p$-integral curvature
of $K$ if $\mu$ vanishes on great sub-spheres.
}

\vs{3}

These results are proved by using variational arguments.  The associated optimization problem is the following,
\[
\sup_{\rho \in C^+(\sn)} \Big\{\frac{E(\bla \rho\bra^*)}{n\omega_n} - \frac1p \log\int_{\sn} \rho^{-p}\, d\mu \Big\},
\]
where $E(\bla \rho\bra^*)$ is the entropy of convex body $\bla \rho\bra^*$, see \eqref{ent}.

When $-1<p<0$,  the symmetric $L_p$ Aleksandrov problem has the following solution:

\vs{3}
{\it Let $\mu$ be a finite even Borel measure on $\sn$ and $p\in (-1,0)$. Then there exists an origin-symmetric
convex body $K$ in $\rn$ such that $C_{p,0}(K,\cdot)=\mu$ if and only if $\mu$ is not concentrated on a great subsphere.
}

This result was first proved by Zhao \cite{Zhao3} for the discrete case and was recently proved by Mui \cite{Mui}
for the general case. Although the approaches are still variational, their proofs are different from that in \cite{HLYZ18jdg}.

For $p\le -1$, a sufficient measure concentration condition for the existence of solution to
the symmetric $L_p$ Aleksandrov problem was given by Mui \cite{Mui}.

\subsubsection{$L_p$ Dual Minkowski Problem}

In \cite{LYZ18adv}, LYZ constructed $L_p$ dual curvature measures which include
the $L_p$ surface area measure, the $L_p$ integral curvature, and dual curvature measures.
This shows that those geometric measures in the Brunn-Minkowski theory, in the dual Brunn-Minkowski theory,
and in the $L_p$  Brunn-Minkowski theory can be unified. Therefore, the Minkowski problems
associated with the geometric measures can also be unified.

\fr{$L_p$ Dual Curvature Measures}
Let $K$ be a convex body in $\kon$, $Q$ a star body in $\rn$, and $p, q \in \R$.
The {\it $(p,q)$th dual curvature measure of $K$ with respect to $Q$}, $\wt C_{p,q}(K,Q,\cdot)$,
is a finite Borel measure on $\sn$  defined by
\begin{equation}\label{pdcm}
\int_{\sn} g(v) \, \wt C_{p,q}(K,Q, v) = \frac1n \int_{\sn}
g(\alpha_K(u)) h_K(\alpha_K(u))^{-p}\rho_K(u)^q \rho_Q(u)^{n-q}
\, du,
\end{equation}
for each continuous $g:\sn \to \R$. When $Q$ is the unit ball $B$, denote $\wt C_{p,q}(K,B,\cdot)$
by $\wt C_{p,q}(K,\cdot)$.
The $L_p$ dual curvature measures have the following three important cases:
\begin{itemize}
\sbr{1} $L_p$ surface area measure, \[ \wt C_{p,n}(K,Q,\cdot) = \frac1n S_p(K,\cdot);\]
\sbr{2} dual curvature measure, \[ \wt C_{0,q}(K,B, \cdot) = \wt C_q(K,\cdot);\]
\sbr{3} $L_p$ integral curvature, \[ \wt C_{p,0}(K,B, \cdot) = \frac1n C_{p,0}(K^*, 0).\]
\end{itemize}

The $L_p$ dual curvature measures are differentials of dual quermassintegrals with respect to the $L_p$ Minkowski sum.
To state the differential formulas in a unified way, we define the {\it normalized power function},
$t^{\bar q}$, $t>0$, $q\in \R$,
\[
t^{\bar q} = \begin{cases} \frac1q t^q &q\neq 0 \\ \log t &q=0 \end{cases},
\]
and the {\it normalized dual mixed volume} $\wt V_{\bar q} (K,Q)$ is defined by,
\begin{equation}\label{n-dmv}
\wt V_{\bar q} (K,Q) = \frac1n \int_{\sn}
\Big(\frac{\rho_K}{\rho_Q}\Big)^{\bar q}(u) \rho_Q^n(u)\, du.
\end{equation}
Obviously, when $q\neq 0$, $q \wt V_{\bar q}(K, Q)$ is just the dual mixed volume $\wt V_{q}(K, Q)$; and when
$q=0$, $\wt V_{\bar 0}(K, B)$ is the dual entropy $\wt E(K)$. Then, we have the differential formula,
\[
\frac d{dt} \wt V_{\bar q}(K\psum t\thin\cdot L, Q) \Big|_{t=0} = \int_{\sn} h_L^{\bar p}(v) \, d\wt C_{p,q}(K,Q, v).
\]

\fr{$L_p$ Dual Minkowski Problem}
{\it Let $Q$ be a star body in $\rn$ and  $p, q \in\R$. If $\mu$ is a finite Borel measure
on $\sn$, what are the necessary and sufficient conditions for the existence of
a convex body $K\in \kon$ satisfying
\begin{equation}\label{pdMP}
\wt C_{p,q}(K,Q, \cdot) = \mu\,?
\end{equation}
}

When the given measure $\mu$ has a density $f$, the partial differential equation associated with the
$L_p$ dual Minkowski problem can be written as the following Monge-Amp\`ere equation on $\sn$,
\begin{equation}\label{pdMp}
\det(\nabla_{ij}h + h\delta_{ij}) = h^{p-1} \rho_Q^{q-n}(\nabla h) f.
\end{equation}

The partial differential equations \eqref{Mp1}, \eqref{APp}, \eqref{logp}, \eqref{dualp},
\eqref{p1.6}, \eqref{p1.7}, and \eqref{pAp} are all special cases of \eqref{pdMp}.

\vs{3}

When $p=q=0$, the $L_0$ dual Minkowski problem is the absolutely continuous case of the Gauss image problem
which will be described in the next section. It was solved by B\"or\"oczky-LYZ-Zhao \cite{BLYZZ20cpam}
(see the next section), and for the smooth case it was also solved by Li-Wang \cite{LW18} with a mass transport approach.

For $q>n-1$ and some negative values of $p$, Guang-Li-Wang \cite{GLW23} showed existence of solutions to
the $L_p$ dual Minkowski problems when $\mu$ has a positive smooth density and $Q$ has a smooth boundary.

When $Q$ is the unit ball, the $L_p$ dual Minkowski problem has been studied in a number of papers.
Huang and Zhao \cite{HZ18} proved that it has a solution
if and only if $\mu$ is not concentrated in any closed hemisphere of $\sn$
when $p>0$ and $q<0$. When $p\neq q$, $p, q>0$, and $\mu$ is even, they also showed
that there is a symmetric solution
 if and only if $\mu$ is not concentrated on any great subsphere of $\sn$.
The solution to the symmetric case was extended to the general case by B\"or\"oczky-Fodor \cite{BF19}
when $p>1$, $q>0$, and $p\neq q$.
These results were proved by Chen-Li \cite{CL21jfa} for the larger range of $p>0$ and $q\neq n$
by an expanding Gauss curvature flow,
\begin{equation*}
\frac{\partial X}{\partial t} (u,t) = |X|^\frac{n-q}{1-p} \lkappa^{\frac{1}{1-p}} f(\nu)^{-1} \nu.
\end{equation*}
Chen-Huang-Zhao \cite{CHZ19} obtained regularity results for symmetric solutions when $pq\ge 0$
by a (rescaled) contracting Gauss curvature flow,
\begin{equation*}
\frac{\partial X}{\partial t} (u,t) = - |X|^{n-q} (X\cdot \nu)^p \lkappa f(\nu) \nu,
\end{equation*}
which extends \eqref{aGf1}.

Huang-Zhao \cite{HZ18} used a variational method and B\"or\"oczky-Fodor \cite{BF19} dealt with
the discrete case first by a variational method followed by an approximation. The optimization problem
for the $L_p$ dual Minkowski problem is
\[
\inf_{h\in C^+(\sn)}\Big\{\wt V_q(\wul{h})^{-\frac1q} \Big(\int_{\sn} h^p\, d\mu\Big)^\frac1p\Big\}.
\]

The harder case of the $L_p$ dual Minkowski problem is when $p\le 0$ and $q>0$. Both the logarithmic
Minkowski problem ($p=0$ and $q=n$) and the centro-affine Minkowski problem ($p=-n$ and $q=n$) belong
to this case.

When the given measure $\mu$ is absolutely continuous, uniqueness results were proved by Huang-Zhao \cite{HZ18}
when $p>q$, and nonuniqueness results were proved by Li-Liu-Lu \cite{LLL21} when $p<0<q$.

\subsubsection{$L_p$ Chord Minkowski Problem}

Chord measures were extended to {\it $L_p$ chord measures} \cite{LXYZ21},
\begin{equation}\label{c11.2}
F_{p,q}(K,\eta) = \frac{2q}{\omega_n} \int_{\nu_K^{-1}(\eta)}
(z\cdot \nu_K(z))^{1-p}\wt V_{q-1}(K, z)\, dz,
\ \ \eta\subset \sn,
\end{equation}
for $q>0$ and any $p\in \R$ if $K$ contains the origin in its interior. There is the formula,
\[
dF_{p,q}(K,\cdot) = h_K^{1-p} dF_q(K,\cdot).
\]
It was shown in \cite{LXYZ21} that the differential of the chord integral $I_q$
with respect to $L_p$ Minkowski combinations leads to the $L_p$ chord measure: for $p\neq 0$,
\begin{equation*}
	\left.\frac{d}{dt}\right|_{t=0} I_q(K \psum t\thin\cdot L) = \frac{1}{p}\int_{\sn} h_L^p(v)dF_{p,q}(K,v).
\end{equation*}

The following   {\bf $L_p$ chord Minkowski problem} was studied in \cite{XYZZ23}:

\vskip 3pt

{\it
Let $\mu$ be a finite Borel measure on $\sn$, $p\in\R$ and $q\ge0$.
What are the necessary and sufficient conditions for the existence of a convex body
 $K\in \kon$ that solves the equation,
\begin{equation}\label{8.1}
F_{p,q}(K,\cdot) = \mu\, ?
\end{equation}
}

When $p=1$, this is the chord Minkowski problem, and the $q=1$ case is the $L_p$ Minkowski problem.

When the given measure $\mu$ has a density $f$  that is an integrable nonnegative function on $\sn$,
equation \eqref{8.1} becomes a new Monge-Amp\`ere type partial differential equation,
\begin{equation}\label{8.2}
\det\big(\nabla_{ij}h + h\delta_{ij}) = \frac{h^{p-1}f}{\wt V_{q-1}(\wul{h},\nabla h)}, \ \ \ \text{on $\sn$}.
\end{equation}

The symmetric case of the $L_p$ chord Minkowski problem was solved in \cite{XYZZ23} when $p, q>0$:
\vs{3}

{\it Let $p,q>0$. If $\mu$ is an even finite Borel measure on $\sn$ that is not concentrated on a great subsphere,
then there exists an origin-symmetric convex body $K$ in $\rn$ such that
\begin{align*}
F_{p,q}(K,\cdot)&= \mu, \ \  \text{when } p\neq n+q-1, \\
\frac{F_{p,q}(K,\cdot)}{V(K)} &= \mu, \ \ \text{when } p=n+q-1.
\end{align*}
}

When $q=1$, this result is a solution to the symmetric $L_p$ Minkowski problem, see \cite{L93jdg, LYZ04tams, HLYZ10}.
When $p>1$, the symmetry condition  can be dropped. Based on the techniques of \cite{HLYZ05dcg},
\cite{XYZZ23} gave the following solution:
\vs{3}

{\it Let $p>1,q>0$. If $\mu$ is a finite Borel measure on $\sn$ that is not concentrated
in any closed hemisphere, then
there exists a convex body $K$ with nonnegative support function $h_K\ge 0$ so that
\begin{align*}
dF_q(K,\cdot) &= h_{K}^{p-1} d \mu, \ \  \text{when }p\neq n+q-1, \\
\frac{dF_{q}(K,\cdot)}{V(K)} &= h_K^{p-1} d\mu, \ \ \text{when } p=n+q-1.
\end{align*}
Moreover, $h_K>0$ if $\mu$ is discrete or if $p\ge n$.
}

Again, when $q=1$, this result is a solution to the $L_p$ Minkowski problem,
see \cite{HLYZ05dcg, CW06}.

\subsubsection{Orlicz Minkowski Problems}

Extending the power function in the $L_p$ Brunn-Minkowski theory to a more general function is the natural next step
to further develop the theory. This has led to the so-called Orlicz Brunn-Minkowski theory, see the papers
\cite{GHW14, HLYZ10, Lud09, LYZ10jdg, LYZ10adv, XJL14, ZZW14}. The Orlicz Minkowski problem for surface area measure
can be stated as follows
\cite{CW06,HLYZ10}:

{\it If $\varphi:(0,\infty)\to(0,\infty)$ is a fixed continuous
function and $\mu$ is a finite Borel measure on $\sn$
which is not concentrated on a closed hemisphere of $\sn$, then does
there exist a convex body $K$ in $\R^n$ such that
$$
\varphi(h_K)\, dS(K, \cdot)=d\mu?
$$
}
This equation arises when the $L_p$ norm in the optimization problems associated with $L_p$ Minkowwski problems
is replaced by the Orlicz norm
\[
\|f\|_\phi = \inf\Big\{t>0 : \frac1{|\mu|} \int_{\sn}\phi\big(\frac ft\big)\, d\mu \le \phi(1)\Big\}, \ \
\text{ where } \phi(t) = \int_0^t \frac1{\varphi(s)}\, ds.
\]
Because of the nonhomogeneous nature of the equation, it is of interest to consider the equation
$$
c\, \varphi(h_K)\, dS(K, \cdot)=d\mu,
$$
for some positive constant $c$ that may be dependent on $K$.  A natural choice of the constant is
a power of volume, $c= V(K)^\alpha$, for some $\alpha \in \R$.
The Orlicz Minkowski problem has been studied in
\cite{BBC19, HH12, HLYZ10, JL19, Li14, LL20, Sun18, SL15, SZ19, WXL19, WXL20, Xie, Xie2, XY20}.

It is also of interest to generalize the Orlicz Minkowski problem to geometric measures other than surface
area measure. For example, the equation for the Orlicz Minkowski problem for dual curvature measure, called
the Orlicz dual Minkowski problem, can be written,
\[
c\,\varphi(h_K)\, d\wt C_q(K, \cdot)=d\mu,
\]
with $c=V(K)^\alpha$, for some $\alpha \in \R$.
See \cite{FH21, GHWXY, GHXY, ZXY18}. Gardner et.al. \cite{GHWXY, GHXY} gave a construction of general dual
curvature measures and posed their related Orlicz Minkowski problems.

\subsubsection{Anisotropic Gauss curvature flows}
The solvability of $L_p$ and Orlicz Minkowski problems, except the $L_p$ chord Minkowski problem,
described in this section amounts to solving the following  Monge-Amp\`ere type equation,
\begin{equation}\label{ME}
\varphi(h)G(\nabla h)\det(\nabla_{ij}h+h\delta_{ij})=f.
\end{equation}
This equation is associated with an anistropic Gauss curvature flow, see for example Liu-Lu \cite{LL20},
\begin{equation}
\frac{\partial X(x,t)}{\partial t}=- \frac{(X\cdot \nu)\lkappa }{\varphi(X\cdot \nu)G(X)}f(\nu)\nu+X,
\end{equation}
where $G$ is a positive continuous function in $\rn$.
 For different cases of $\varphi$ and $G$, these problems have been studied in \cite{CHZ19,CW06,LSW20,LSW20jga,LL20}.

\section{Gauss Image Problem and Spherical Mass Transport}\

Recall that the radial Gauss image $\balpha_K$ of a convex body $K$ in $\rn$ defined
by \eqref{radial-gauss-image}
maps a Borel set $\omega \subset \sn$ to a Lebesgue measurable set $\balpha_K(\omega)$.

\subsubsection{Gauss Image Measure}
   Suppose $\lambda$ is a measure defined on
Lebesgue measurable subsets of $\sn$, and $K\in\kon$. The {\it Gauss image measure} $\lambda (K,\cdot)$
associated with $\lambda$ and $K$ is defined by
\[
  \lambda(K,\omega) = \lambda(\balpha_K(\omega)), \ \ \ \ \text{Borel } \omega \subset\sn.
\]
In general, $\lambda(K,\cdot)$ is only an outer measure.
If $\lambda$ is absolutely continuous with respect to the Lebesgue measure on $\sn$ and $K\in \kon$,
then $\lambda(K,\cdot)$ is a Borel measure on $\sn$.

The concept of the Gauss image measure underlies the definitions of many geometric measures presented here.
For example, the integral curvature of $K$ is the Gauss image measure of spherical Lebesgue measure via $K$,
see \eqref{4.2}. The surface area measure $S(K,\cdot)$ is the Gauss image measure of the $(n-1)$th curvature
$C_{n-1}(K, \cdot)$ via $K^*$, see \eqref{4.3}.
Dual curvature measures are also  Gauss image measures, see \eqref{5.4}.

\subsubsection{Gauss Image Problem}

Motivated by this, the following problem was posed in \cite{BLYZZ20cpam}:

\vspace{3pt}

{\it
Suppose $\lambda$ is a measure defined on the
Lebesgue measurable subsets of $\sn$ and $\mu$ is a Borel measure on $\sn$.
What are the necessary and sufficient conditions, on $\lambda$ and $\mu$, so that
there exists a convex body $K\in\kon$ such that
\begin{equation}\label{gip}
\lambda(K,\cdot\,)=\mu\,?
\end{equation}
}

When $\lambda$ is spherical Lebesgue measure, the Gauss image problem is
the Aleksandrov problem. The essence of the Gauss image problem is
to understand how measures on $\sn$ are transformed by the Gauss image map of a convex body.
Thus, the Gauss image problem is essentially a problem of spherical mass transport.

When $\mu$ has a smooth density $f$ and $\lambda$ has a smooth density $g$,
the geometric equation \eqref{gip} can be written as a Monge-Amp\`ere type equation,
\begin{equation}\label{gp}
g\Big(\frac{\nabla h}{|\nabla h|}\Big)  \det\big(\nabla_{ij}h + h \delta_{ij}\big) = \frac{|\nabla h|^n} hf.
\end{equation}
Note that this equation is a special case of the equation \eqref{pdMp} when $g$ is the $n$th power of the
radial function of a star body and $p=q=0$. Thus, the Gauss image problem becomes the Aleksandrov problem.

\vspace{3pt}

\subsubsection{Aleksandrov Related Measures}

There is an obvious necessary condition required for the Gauss image problem to have a solution
for two given measures $\lambda$ and $\mu$. First, recall that if $\omega\subset \sn$ is contained
in a closed hemisphere, then its {\it polar set} $\omega^*$ is given by
\begin{equation*}
\omega^* =\{v\in \sn : \text{$u\cdot v\leq 0$, for all $u\in\omega$}\}=
\bigcap_{u\in\omega}\{v\in \sn : u\cdot v\leq 0 \}.
\end{equation*}
A subset $\omega \subset \sn$ is called spherically convex if the cone $\{tu : t\ge 0, u\in\omega\}$ is a
nonempty convex proper subset of $\rn$.
Given a convex body $K \subset \rn$ and a spherically convex $\omega \subset \sn$,
it is not hard to show that (see \cite[Lemma 3.2]{BLYZZ20cpam})
\[
  \balpha_K(\omega) \subset \sn\setminus \omega^*,
\]
and the set $(\sn\setminus \omega^*)\setminus \balpha_K(\omega)$ has interior points.
It follows that given any spherically convex set $\omega \subset \sn$,
\[
 \lambda(K,\omega) = \lambda(\balpha_K(\omega)) \le \lambda(\sn) - \lambda(\omega^*).
\]
Moreover, if $\lambda$ is assumed to be positive for every nonempty open subset of $\sn$,
then strict inequality must hold. This necessary condition motivates the following:

{\it Two Borel measures $\lambda$ and $\mu$ are said to be Aleksandrov related,
if for any spherically convex compact set $\omega \subset \sn$,
\[
  \lambda(\omega^*) + \mu(\omega) < \lambda(\sn)=\mu(\sn).
\]}
Since $\omega^{**} = \omega$, this is a symmetric relationship.
By the discussion above, if $\lambda$ is positive on open subsets of $\sn$,
then this is a necessary condition for $\mu$ to be the Gauss image measure of $\lambda$.
Below we will see that it is also a sufficient condition.

\subsubsection{Solutions to the Gauss Image Problem}

In \cite{BLYZZ20cpam} the following theorem was proved:
\vs{3}

{\it Let $\mu$ be a finite Borel measure on $\sn$ and
$\lambda$ be a measure on $\sn$ that is absolutely continuous
with respect to the spherical Lebesgue measure.
If $\mu$ and $\lambda$ are Aleksandrov related, then
there exists a convex body $K$ in $\kon$ such that $\mu=\lambda(K,\cdot\,)$.
}

\vs{3}

When both $\mu$ and $\lambda$ have positive densities, the result was proved also by
Li-Wang \cite{LW18} by the approach of spherical mass transport.
When $\mu$ is discrete and $\lambda$ is absolutely continuous, the result was strengthened
by Semenov \cite{S24} who introduced the so-called {\it weakly Aleksandrov related} condition.

It is important to study the Gauss image problem when
 $\lambda$ is a discrete measure,
due to connections with discrete mass transport and combinatorics.
Significant progress to this end was made by Semenov \cite{S22}.

\section{{\hspace{-11pt}} Capacitary Minkowski Problems from Harmonic Analysis}\

In his remarkable work \cite{Je96}, Jerison introduced the electrostatic capacitary measure of a convex body
as the differential of electrostatic (Newtonian, or $L_2$) capacity and solved its Minkowski problem completely;
see also \cite{Je96adv} for a direct proof.
This establishes amazing connections between convex geometry and harmonic analysis related to partial differential
equations.
Important connections of convex geometry to harmonic analysis, functional analysis and probability
theory related to geometric and analytic inequalities were well described by Ball \cite{B-ICM}.
In the paper \cite{CNSXYZ}, Colesanti-Nystr\"om-Salani-Xiao-Yang-Zhang
extended the concept of electrostatic capacitary measure
to $L_p$ capacitary measure which is the differential of the $L_p$ capacity. They solved the $L_p$ capacitary
Minkowski problem in $\rn$ when $1<p<2$. Then, Akman-Gong-Hineman-Lewis-Vogel
\cite{AGHLV} solved the $L_p$ capacitary Minkowski problem when $2<p<n$.

\subsubsection{$L_p$ Capacitary Measure}
Let $K$ be a convex body in $\rn$ and $1<p<n$.  The $L_p$ capacity Cap$_p(K)$ of $K$ is defined by
\begin{equation*}
   \text{Cap}_p(K)=\inf\Big\{\int_{\rn}|\nabla u|^pdx:\ u\in C_c^\infty(\rn)\  \text{and $u\ge 1$ on $K$} \Big\}.
 \end{equation*}
The $L_2$ case Cap$_2(K)$ is the electrostatic capacity.

The $p$-equilibrium potential $U_K$ of $K$ is the unique solution to the boundary value problem
\begin{equation*}
\left\{
\begin{array}{lll}
\Delta_p U=0\quad\mbox{in $\rn \setminus K$,}\\
\\
U=1\quad\quad\mbox{on $\partial K$ \ \ \  and \ \ \ } \ds \lim_{|x|\to\infty}U(x)=0\,,
\end{array}
\right.
\end{equation*}
where $\Delta_p$ is the $p$-Laplace operator,
\[
\Delta_pU= \nabla \cdot ( |\nabla  U|^{p-2} \, \nabla  U ).
\]
A  proof of the existence and uniqueness of $U_K$ can be found in \cite{Le77}.
When approaching the boundary of $\rn\setminus K$ non-tangentially, the limit of $\nabla U_K$
exists $\Hn$-almost everywhere on $\partial K$ and its norm $|\nabla U_K|$ is $L_p$ integrable
over $\partial K$, i.e.,
\[
\int_{\partial K} |\nabla U_K|^p\, dx < \infty.
\]
This is a well-known result of Dahlberg \cite{D} when $p=2$, and the general case
was shown by Lewis and Nystr\"om  \cite{LN}.
Therefore, a measure $\mu_p(K,\cdot)$ of $K$, called
the {\it $p$-capacitary measure of $K$}, can be defined by
\begin{equation*}
\mu_p(K, \eta)=\int_{\nu_K^{-1}(\eta)}|\nabla U_K|^p\,dx, \quad \text{ for any Borel set } \eta \subset \sn.
\end{equation*}

An important fact is that the $p$-capacitary measure is the differential of the $p$-capacity, that is,
\begin{equation}\label{capv}
\frac{d}{dt} \text{Cap}_p(K+tL)\Big|_{t=0} = (p-1) \int_{\sn} h_L \, d\mu_p(K,\cdot).
\end{equation}
The proof of this differential formula is not easy. The case of $p=2$ was proved by Jerison \cite{Je96}
and the general case $1<p<n$ was proved in \cite{CNSXYZ}. It is crucial for solving the following capacitary
Minkowski problem.

\subsubsection{The $p$-Capacitary Minkowski Problem}
{\it Let $\mu$ be a finite Borel measure on $\sn$ and $1<p<n$. Under what necessary and sufficient
conditions does there exist a convex body $K$ in $\rn$ such that
$\mu_p(K, \cdot) = \mu\,?$}

If the measure $\mu$ has a density $f$, the equation associated with this Minkowski
problem is a Monge-Amp\`ere type equation on $\sn$
\[
\det\big(\nabla_{ij} h + h\delta_{ij}\big) = \frac f{|\nabla U_{\swul{h}}|^p},
\]
where $h$ is the unknown positive function on $\sn$ to be found, $U_{\swul{h}}$ is the $p$-equilibrium potential
of the Wulff shape $\wul{h}$, and $|\nabla U_{\swul{h}}|$ is viewed as a function on $\sn$ under the inverse Gauss map
$\nu_{\swul{h}}^{-1}$.

\subsubsection{Solution to the $p$-Capacitary Minkowski Problem}
{\it Let $\mu$ be a finite Borel measure on $\sn$ and $1<p<n$.
Then there exists a convex body $K$ in $\rn$ such that
$\mu_p(K, \cdot) = \mu$
if and only if $\mu$ is not concentrated on a closed hemisphere and
\[
\int_{\sn} v\, d\mu(v) =0.
\]}

This was first proved by Jerison \cite{Je96} when $p=2$, then Colesanti, {et al.} \cite{CNSXYZ}
proved the case of $1<p<2$ under an extra technical condition, and finally Akman, {et al.} \cite{AGHLV}
proved the general case.

\subsubsection{Uniqueness Problem of the $p$-Capacitary Measure}
{\it
Let $K_0, K_1$ be convex bodies in $\rn$ and $1<p<n$. If
\[
\mu_p(K_0,\cdot) =\mu_p(K_1,\cdot),
\]
then $K_0$ is a translate of $K_1$ when $p\neq n-1$, and $K_0, K_1$ are homothetic when $p= n-1$.
}

\vspace{5pt}

This was first proved by Caffarelli-Jerison-Lieb in \cite{CJL96} when $p=2$,
and the general case was proved by Colesanti et al. \cite{CNSXYZ} by using the differential formula \eqref{capv}
and the capacitary Brunn-Minkowski inequality of Colesanti-Salani \cite{CS03}.

\vspace{5pt}

\subsubsection{$(p,q)$-Capacitary Minkowski Problem}

In \cite{ZX20}, Zou and Xiong defined the {\it $(p,q)$-capacitary measure} of a convex body $K$,
\[
d\mu_{p,q}(K, \cdot) = h_K^{1-q} d\mu_p(K,\cdot), \ \ \ 1<p<n, \ q\in \R,
\]
which is the differential of the $p$-capacity with respect to $L_q$ Minkowski addition,
\[
\frac{d\text{Cap}_p(K\qsum{q} t \thin\cdot L)}{dt}\Big|_{t=0} = \frac{p-1}q \int_{\sn} h_L^q(v) \, d\mu_{p,q}(K, v),
\ \ \ 1<p<n, \ q \ge 1.
\]
They solved the associated $(p,q)$-Capacitary Minkowski problem when $1<p<n$ and $q>1$. In particular,
for the discrete case, the following result was proved:

{\it If $\mu$ is a discrete measure that is not concentrated on any closed hemisphere, then
for $1<p<n$ and $q>1$ with $p+q\neq n$, there exists a unique convex polytope $K$ containing
the origin in its interior such that
\[
\mu_{p,q} (K, \cdot) = \mu.
\]}

Existence proofs for capacitary Minkowski problems take a variational
approach, relying on \eqref{capv} and its analogues.
It is of interest to use continuity methods or geometric flows to find smooth solutions
to the capacitary Minkowski problems.

Generalizations of the capacitary Minkowski problems and other Minkowski-type problems
for various capacities have been studied in
\cite{ALSV, ALV, C20, CD20, CF10, HYZ18, LuoYeZhu, X15, X16, X17, X18, X20, XXL19, ZX20}.

\section{{\hspace{-11pt} Gaussian Minkowski Problems from Probability Theory}}\

Geometric measures of convex bodies that are related to probability measures are also of great interest.
Huang-Xi-Zhao \cite{HXZ20} constructed the {\it Gaussian surface measure} of convex bodies from the Gaussian probability measure
and studied its {\it Gaussian Minkowski problem}. A general construction of a geometric measure of a convex body
 from a probability measure was given by Livshyts \cite{Liv19}. This provides new connections among convex geometry,
 probability theory, and partial differential equations.

\subsubsection{Gaussian Surface Area Measure}
Let $K\in \kon$. The Gaussian measure $\gamma_n(K)$ of $K$ is defined by
\[
 \gamma_n(K) = \frac1{(2\pi)^{n/2}} \int_{K} e^{-\frac{|x|^2}2}\, dx.
\]
The {\it Gaussian surface area measure $S_{\gamma_n}(K,\cdot)$ of $K$} is a Borel measure on $\sn$
defined by
\begin{equation}\label{Gsm}
 S_{\gamma_n}(K,\eta) = \frac1{(2\pi)^{n/2}} \int_{\nu_K^{-1}(\eta)} e^{-\frac{|x|^2}2}\, dx, \quad
 \text{ for any Borel set } \eta \subset \sn.
 \end{equation}

 The measure $S_{\gamma_n}(K,\cdot)$ is the differential of the Gaussian measure of $K$:
 \begin{equation}\label{G1}
 \frac{d}{dt} \gamma_n(K+tL)\Big|_{t=0} = \int_{\sn} h_L\, dS_{\gamma_n}(K,\cdot).
 \end{equation}

The total measure $S_{\gamma_n}(K,\sn)$ is a known invariant, called the
{\it Gaussian surface area (or Gaussian perimeter) of $K$}.
It was first studied by Borell \cite{Bo75} for proving his Gaussian isoperimetric inequality, see also Lata\l{}a \cite{La02}.
Then as a different invariant, it was studied by
Ball \cite{B93} for a reverse Gaussian isoperimetric inequality, see also Nazarov \cite{N03}. The formula \eqref{G1}
proved by Huang-Xi-Zhao \cite{HXZ20} clarifies the connections among these studies and shows a striking similarity
between the Brunn-Minkowski theory in the Euclidean space and its analog in the Gaussian space.

\subsubsection{Uniqueness Problem of Gaussian Surface Area Measure}
Let $K, L$ be convex bodies in $\kon$ with $\gamma_n(K), \gamma_n(L) \ge \frac12$. It is shown in \cite{HXZ20}
that
\[
S_{\gamma_n}(K,\cdot) =  S_{\gamma_n}(L,\cdot) \ \ \Longrightarrow \ \  K=L.
\]

Huang-Xi-Zhao \cite{HXZ20} showed this uniqueness result by using the Ehrhard inequality, which is a Brunn-Minkowski type
inequality for the Gauss measure, see \cite{E83, E86}.

\subsubsection{Gaussian Minkowski problem} {\it Given a finite Borel measure $\mu$ on $\sn$, what are
the necessary and sufficient conditions on $\mu$ so that there exists a convex body $K\in \kon$ satisfying
\[
S_{\gamma_n}(K,\cdot) =\mu\, ?
\]}

When the given measure has a density, the partial differential equation associated with the Gaussian Minkowski problem
can be written as the Monge-Amp\`ere type equation on $\sn$,
\begin{equation}\label{GMp}
\det\big(\nabla_{ij}h + h\delta_{ij}\big) = e^{\frac{|\nabla h|^2}2} f.
\end{equation}

\subsubsection{Existence of Solutions to the Symmetric Gaussian Minkowski Problem}
{\it Let $\mu$ be a finite even Borel measure on $\sn$ that is not concentrated on any great subsphere.
Then for any $0<\alpha<\frac1n$, there exists an origin-symmetric convex body $K$ such that
\[
\mu = \frac1{\gamma_n(K)^{1-\alpha}} S_{\gamma_n}(K,\cdot).
\]}

A variational proof of this result is given in \cite{HXZ20}, which is similar to the proof
given in \cite{HLYZ10}. By using the degree theory for second-order nonlinear elliptic operators developed in
Li \cite{Li89} for the smooth case,  followed by an approximation argument, Huang-Xi-Zhao \cite{HXZ20}
also proved the following non-normalized solution:

{\it Let $\mu$ be an even measure on $\sn$ that is not concentrated
on any great subsphere and $|\mu|< \frac1{\sqrt{2\pi}}$. Then there exists a unique origin-symmetric
$K$ with $\gamma_n(K) > \frac12$ such that $S_{\gamma_n}(K,\cdot)=\mu$.}

A similar solution to the non-symmetric case was obtained by Feng-Liu-Xu \cite{FLX23}.

Related studies on the Gaussian Minkowski problem and its $L_p$ variant
can be found in the papers \cite{Liu22, CHLZ23, FHX22, SX22}.

A more general construction is to replace the Gaussian
density by any density function $g$ of a Borel measure in $\rn$,
\[
 S_{g}(K,\eta) = \int_{\nu_K^{-1}(\eta)} g(x)\, dx, \quad
 \text{ for any Borel set } \eta \subset \sn.
\]
This was defined by Livshyts \cite{Liv19}. When $g$ is even, homogeneous, and of certain concavity,
the corresponding Minkowski problem was solved in \cite{Liv19}, extending the solution of the classical Minkowski problem
for symmetric convex bodies. The results in \cite{Liv19} were
generalized by Kryvonos-Langharst \cite{KL23}.

\vspace{3pt}

While convex geometry primarily focuses on convex bodies in $\rn$, one can often replace convex bodies
by functions in $\rn$, usually with certain convexity properties. It is interesting to construct geometric measures
from functions and solve the corresponding Minkowski problems. In \cite{LYZ06imrn}, the $L_p$ surface area measure
$S_p(f, \cdot)$, $p\ge 1$, a Borel measure on $\sn$, was defined for $f: \rn \to \R$ with $L_p$ weak derivative,
\[
\int_{\sn}  \varphi(u)^p dS_p(f, u) = \int_{\rn} \varphi(-\nabla f(x))^p \, dx,
\]
where $\varphi$ is a nonnegative continuous homogeneous function of degree 1 in $\rn$.
A similar concept in $\rn$, called the moment measure $\mu(f, \cdot)$, was defined in \cite{CK15} for a
nonnegative integrable log-concave function $f$ in $\rn$,
\[
\int_{\rn} \varphi(x)\, d\mu(f, x) = \int_{\rn} \varphi(-\nabla \log f(x)) f(x)\, dx
\]
where $\varphi$ is a nonnegative Borel function in $\rn$.

Another direction of study is to replace convex bodies by pseudo-cones. They are unbounded closed convex sets not containing the
origin and viewed as the counterpart to convex bodies containing the origin. Schneider \cite{S18} defined
the surface area measure for pseudo-cones and studied its Minkowski problem. Further studies were carried out in
\cite{AYY, LiYZ, S23, XLL23, YYZ23}.

\vspace{5pt}

\section*{Acknowledgements}
The authors are very grateful to
K\'aroly B\"or\"oczky,  Richard Gardner, Pengfei Guan, Jinrong Hu,  Monika Ludwig,  Erwin Lutwak, Rolf Schneider,
Xu-Jia Wang,  Dongmeng Xi, and Yiming Zhao
for helpful comments during the writing of this paper.


\begin{thebibliography}{999}

\bibitem{AGHLV}
M. Akman, J. Gong, J. Hineman, J. Lewis, A. Vogel,
\emph{The Brunn-Minkowski inequality and a Minkowski problem for nonlinear capacity},
Mem. Amer. Math. Soc.  275 (2022),  no. 1348.

\bibitem{ALSV}
M. Akman, J. Lewis, O. Saari, A. Vogel,
\emph{The Brunn-Minkowski inequality and a Minkowski problem for A-harmonic Green's function},
Adv. Calc. Var.  14 (2021), 247-302.

\bibitem{ALV}
M. Akman, J. Lewis, A. Vogel,
\emph{Note on an eigenvalue problem with applications to a Minkowski type regularity problem in $\rn$},
Calc. Var. PDE 59 (2020), Paper No. 47.

\bibitem{AYY}
W. Ai, Y. Yang, D. Ye,
\emph{The $L_p$ dual Minkowski problem for unbounded closed convex sets},
arXiv:2404.09804v1.


\bibitem{Al2} A.D. Aleksandrov,
{\it On the theory of mixed volumes. III. Extension of two theorems of Minkowski on convex polyhedra
to arbitrary convex bodies},
Mat. Sbornik N.S.  3 (1938), 27-46.

\bibitem{Al3} A.D. Aleksandrov,
{\it On the surface area measure of convex bodies},
Mat. Sbornik N.S.  6 (1939), 167-174.

\bibitem{A42} A.D. Aleksandrov,
\emph{Existence and uniqueness of a convex surface with a given integral curvature},
C.R. (Doklady) Acad. Sci. URSS 35 (1942), 131-134.

\bibitem{A56} A.D. Aleksandrov,
\emph{Uniqueness theorems for surfaces in the large I},
Vestnik Leningrad. Univ. 11 (1956), 5-17.

\bibitem{And99}
B. Andrews,
\emph{Gauss curvature flow: the fate of the rolling stones},
Invent. Math. 138 (1999), 151-161.

\bibitem{And00}
B. Andrews,
\emph{Motion of hypersurfaces by Gauss curvature},
Pacific J. Math. 195 (2000), no. 1, 1-34.



\bibitem{AC12}
B. Andrews, X. Chen,
\emph{Surfaces moving by powers of Gauss curvature},
Pure Appl. Math. Q. 8 (2012), no. 4, 825-834.

\bibitem{AGN16}
B. Andrews, P. Guan, L. Ni,
\emph{Flow by powers of the Gauss curvature},
Adv. Math. 299 (2016), 174-201.

\bibitem{B93}
K. Ball,
\emph{The reverse isoperimetric problem for Gaussian measure},
Discrete Comput. Geom. 10 (1993), 411-420.

\bibitem{B-ICM}
K. Ball,
\emph{Convex geometry and its connections to harmonic analysis, functional analysis and probability theory},
Proc. Int. Cong. Math. 2022, Vol. 4, pp. 3104-3139.

\bibitem{B69}
C. Berg, \emph{Corps convexes et potentiels sphe\'riques}, Danske Vid. Selsk., Mat.-Fys. Medd. 37(b)
(1969), 1-64.

\bibitem{B16}
J. Bertrand,
\emph{Prescription of Gauss curvature using optimal mass transport},
Geom. Dedicata {183} (2016), no. 1, 81-99.


\bibitem{BBC19}
G. Bianchi, K.J. B\"or\"oczky, A. Colesanti,
\emph{The Orlicz version of the $L_p$ Minkowski problem for $-n < p < 0$},
Adv. Appl. Math. 111 (2019), Paper No. 101937.

\bibitem{BBCY}
G. Bianchi, K.J. B\"or\"oczky, A. Colesanti, D. Yang,
\emph{The $L_p$-Minkowski problem for $-n < p < 1$},
Adv. Math. 341 (2019), 493-535.


\bibitem{Bo75}
C. Borell,
\emph{The Brunn-Minkowski inequality in Gauss space},
Inventiones Math. 30 (1975), 205-216.


\bibitem{BF19}
K.J. B\"or\"oczky,  F. Fodor,
\emph{The $L_p$ dual Minkowski problem for $p > 1$ and $q > 0$},
J. Differential Equ. 266 (2019), 7980-8033.

\bibitem{BHZ16}
K.J. B\"or\"oczky, P. Heged\H{u}s, G. Zhu,
\emph{On the discrete logarithmic Minkowski problem},
Int. Math. Res. Not. 6 (2016), 1807-1838.

\bibitem{BH16}
K.J. B\"or\"oczky, M. Henk,
\emph{Cone-volume measure of general centered convex bodies},
Adv. Math. 286 (2016), 703-721.

\bibitem{BHP18}
K.J. B\"or\"oczky, M. Henk, H. Pollehn,
\emph{Subspace concentration of dual curvature measures of symmetric convex bodies},
J. Differ. Geom. 109 (2018), 411-429.



\bibitem{BLYZ12adv}
K.J. B\"or\"oczky, E. Lutwak, D. Yang, G. Zhang,
\emph{The log-Brunn-Minkowski inequality},
Adv. Math.  231 (2012), 1974-1997.

\bibitem{BLYZ13jams}
K.J. B\"or\"oczky, E. Lutwak, D. Yang, G. Zhang,
\emph{The logarithmic Minkowski problem},
J. Amer. Math. Soc.  26 (2013), 831-852.

\bibitem{BLYZZ19adv}
K.J. B\"or\"oczky, E. Lutwak, D. Yang, G. Zhang, Y. Zhao,
\emph{The dual Minkowski problem for symmetric convex bodies},
Adv. Math.   356 (2019), Paper No. 106805.


\bibitem{BLYZZ20cpam}
K.J. B\"or\"oczky, E. Lutwak, D. Yang, G. Zhang, Y. Zhao,
\emph{The Gauss image problem},
Comm. Pure Appl. Math.  73 (2020), 1406-1452.

\bibitem{BT17}
K.J. B\"or\"oczky, H.T. Trinh,
\emph{The planar $L_p$-Minkowski problem for $0 < p < 1$},
Adv. Appl. Math. 87 (2017), 58-81.



\bibitem{BBM01}
J. Bourgain, H. Brezis, P. Mironescu,
\emph{Another look at Sobolev spaces},
In (J.L. Menaldi, E. Rofman and A. Sulem, eds.), Optimal Control and Partial Differential Equations,
a volume in honor of A. Bensoussans's 60th birthday, Amsterdam:
IOS Press; Tokyo: Ohmsha, 2001.




\bibitem{BCD17}
S.\  Brendle, K.\ Choi, P.\ Daskalopoulos,
\emph{Asymptotic behavior of flows by powers of the Gaussian curvature},
Acta Math. 219 (2017), 1-16.

\bibitem{BIS19}
P. Bryan, M. N. Ivaki, J. Scheuer,
\emph{A unified flow approach to smooth, even $L_p$-Minkowski problems},
Anal. PDE  12 (2019), 259-280.



\bibitem{BIS20}
P. Bryan, M. N. Ivaki, J. Scheuer,
\emph{Christoffel-Minkowski flows},
Trans. Amer. Math. Soc. 376 (2023), no. 4, 2373-2393.



\bibitem{Caf90Ann}
L. Caffarelli,
\emph{Interior $W^{2,p}$-estimates for solutions of the Monge-Amp\`ere
equation},
Ann. of Math. 131 (1990), 135-150.


\bibitem{CJL96}
L.A. Caffarelli, D. Jerison, E.~H. Lieb,
\emph{On the case of equality in the Brunn-Minkowski inequality for capacity},
Adv. Math.  117 (1996), 193-207.

\bibitem{CLWX1}
X. Cai, G. Leng, Y. Wu, D. Xi,
\emph{Affine dual Minkowski problems},
Preprint.


\bibitem{CG02}
S. Campi, P. Gronchi,
\emph{The $L_p$-Buseman-Petty centroid inequality},
Adv. Math. 167 (2002), 128-141.

\bibitem{CHZ19}
C. Chen, Y. Huang, Y. Zhao,
\emph{Smooth solutions to the $L_p$ dual Minkowski problem},
Math. Ann. 373 (2019), 953-976.

\bibitem{CL21jfa}
H. Chen, Q.-R. Li,
\emph{The $L_p$ dual Minkowski problem and related parabolic flows},
J. Funct. Anal. 281 (2021), Paper No. 109139.



\bibitem{CL20}
L. Chen,
\emph{Uniqueness of solutions to $L_p$-Christoffel-Minkowski problem for $p<1$},
J. Funct. Anal. 279 (2020), Paper No. 108692.

\bibitem{CHLZ23}
S. Chen, S. Hu, W. Liu, Y. Zhao,
\emph{On the planar Gaussian-Minkowski problem},
Adv. Math. 435 (2023), Paper No. 109351, 32 pp.


\bibitem{CHLL20}
S. Chen, Y. Huang, Q.-R. Li, J. Liu,
\emph{$L_p$-Brunn-Minkowski inequality for $p \in \big(1 - \frac c{n^{3/2}}, 1\big)$},
Adv. Math. 368 (2020), Paper No. 107166.

\bibitem{CL18}
S. Chen, Q.-R. Li,
\emph{On the planar dual Minkowski problem}, Adv. Math. 333 (2018), 87-117.

\bibitem{CLZ17}
S. Chen, Q.-R. Li, G. Zhu,
\emph{On the $L_p$ Monge-Amp\`ere equation},
J. Differential Equ. 263 (2017), 4997-5011.

\bibitem{CLZ19}
S. Chen, Q.-R. Li, G. Zhu,
\emph{The logarithmic Minkowski problem for non-symmetric measures},
Trans. Amer. Math. Soc. 371 (2019), 2623-2641.


\bibitem{CDS15I}
X. Chen, S. Donaldson, S. Sun,
\emph{{K\"ahler-Einstein metrics on Fano manifolds. I: Approximation of metrics with cone singularities}},
J. Amer. Math. Soc.  28 (2015), 183-197.

\bibitem{CDS15II}
X. Chen, S. Donaldson, S. Sun,
\emph{{K\"ahler-Einstein metrics on Fano manifolds. II: Limits with cone angle less than $2\pi$}},
J. Amer. Math. Soc. 28 (2015), 199-234.

\bibitem{CDS15III}
X. Chen, S. Donaldson, S. Sun,
\emph{{K\"ahler-Einstein metrics on Fano manifolds. III:
Limits as cone angle approaches $2\pi$ and completion of the main proof}},
J. Amer. Math. Soc. 28 (2015), 235-278.

\bibitem{C06}
W. Chen,
\emph{$L_p$ Minkowski problem with not necessarily positive data},
Adv. Math. 201 (2006), 77-89.

\bibitem{C20}
Z. Chen,
\emph{The $L_p$ Minkowski problem for $q$-capacity},
Proc. Royal Soc. Edin.,  Section A Math. 151 (2021), 1247-1277.

\bibitem{CD20}
Z. Chen, Q. Dai,
\emph{The $L_p$ Minkowski problem for torsion},
J. Math. Anal. Appl. 488 (2020), Paper No. 124060.

\bibitem{ChengYau}
S.-Y.~Cheng, S.-T.~Yau,
\emph{On the regularity of the solution of the $n$-dimensional Minkowski problem},
Comm. Pure Appl. Math.
29 (1976), 495-516.


\bibitem{Chern1} S.-S. Chern,
\emph{Differential geometry and integral geometry},
Proc. Int. Cong. Math. 1958, 441-449, Cambridge Univ. Press, New York, 1960.

\bibitem{Chern2} S.-S. Chern,
\emph{Integral formulas for hypersurfaces in Euclidean space and their applications to
uniqueness theorems}, J. Math. Mech. 8 (1959), 947-955.


\bibitem{Choi17}
K. Choi, \emph{The Gauss Curvature Flow: Regularity and Asymptotic Behavior}, Ph.D. Thesis,
Columbia University, New York, NY, 2017.

\bibitem{CW00}
K.-S. Chou,  X.-J. Wang,
\emph{A logarithmic Gauss curvature flow and the Minkowski problem},
Ann. Inst. H. Poincar\'{e} Anal. Non Lin\'{e}aire  17 (2000), 733-751.

\bibitem{CW06}
K.-S.~Chou, X.-J.~Wang,
\emph{The {$L\sb p$}-Minkowski problem and the Minkowski problem
in centroaffine geometry},
Adv. Math. 205 (2006), 33-83.

\bibitem{CZ01}
K.-S. Chou, X.-P. Zhu,
\emph{The curve shortening problem},
Chapman \& Hall/CRC, Boca Raton, FL, 2001.

\bibitem{Chow87}
B. Chow,
\emph{Deforming convex hypersurfaces by the square root of the scalar curvature},
Invent. Math. 87 (1987), no. 1, 63-82.


\bibitem{C1865}
E.B. Christoffel,
\emph{\"Uber die Bestimmung der Gestalt einer krummen Oberfl\"ache durch lokale Messungen auf derselben},
J. Reine Angew. Math. 64 (1865), 193-209.


\bibitem{CLYZ}
A.~Cianchi, E.~Lutwak, D.~Yang, G.~Zhang,
\emph{Affine Moser-Trudinger and Morrey-Sobolev inequalities},
Calc. Var. PDE 36 (2009), 419-436.

\bibitem{CF10}
 A. Colesanti, M. Fimiani,
\emph{The Minkowski problem for the torsional rigidity},
Indiana Univ. Math. J. 59 (2010), 1013-1039.

\bibitem{CLi21}
A. Colesanti, G.V. Livshyts,
\emph{A note on the quantitative local version of the log-Brunn-Minkowski inequality},
Adv. Anal. Geom. 2 (2020), 85-98.

\bibitem{CLM17}
A. Colesanti, G.V. Livshyts, A. Marsiglietti,
\emph{On the stability of Brunn-Minkowski type inequalities},
J. Funct. Anal. 273 (2017), 1120-1139.

\bibitem{CNSXYZ}
A. Colesanti, K. Nystr\"om, P. Salani, J. Xiao, D. Yang, G. Zhang,
\emph{The Hadamard variational formula and the Minkowski problem for $p$-capacity},
Adv. Math. 285 (2015), 1511-1588.


\bibitem{CS03}
A. Colesanti, P. Salani,
\emph{The Brunn-Minkowski inequality for $p$-capacity of convex bodies},
Math. Ann.  327 (2003), 459-479.

\bibitem{CK15}
D. Cordero-Erausquin, B. Klartag,
\emph{Moment measures},
J. Funct. Anal. 268 (2015), 3834-3866.

\bibitem{D}
B.E.J. Dahlberg,
\emph{Estimates of harmonic measure},
Arch. Rational Mech. Anal.  65 (1977), 275-288.

\bibitem{DZ12}
J. Dou, M. Zhu,
\emph{The two dimensional $L_p$ Minkowski problem and nonlinear equations with negative exponents},
Adv. Math. 230 (2012), 1209-1221.


\bibitem{E83}
A. Ehrhard,
\emph{Sym\'etrisation dans l'espace de Gauss},
Math. Scand.  53(1983), 281-301.


\bibitem{E86}
A. Ehrhard,
\emph{\'El\'ements extr\'emaux pour les in\'egalit\'es de Brunn-Minkowski gaussiennes},
Ann. Inst. H. Poincar\'e Probab. Statist.  22 (1986), 149-168.


\bibitem{F59}
H. Federer,
\emph{Curvature measures},
Trans. Amer. Math. Soc. 93 (1959), 418-491.


\bibitem{FJ}
W.~Fenchel, B.~Jessen,
\emph{Mengenfunktionen und konvexe K\"orper},
Danske Vid. Selskab. Mat.-fys. Medd. 16 (1938), 1-31.

\bibitem{FH21}
Y. Feng, B. He,
\emph{The Orlicz Aleksandrov problem for Orlicz integral curvature},
Int. Math. Res. Not.  2021 (2021), 5492-5519.

\bibitem{FHX22}
Y. Feng, S. Hu, L. Xu,
\emph{On the $L_p$ Gaussian Minkowski problem},
 arXiv:2211.10956.




\bibitem{FLX23}
Y. Feng, W. Liu, L. Xu,
\emph{Existence of non-symmetric solutions to the Gaussian Minkowski problem},
J. Geom. Anal. 33 (2023), Paper No. 89.




\bibitem{F62}
W. Firey,
\emph{$p$-means of convex bodies},
Math. Scand. 10 (1962), 17-24.

\bibitem{F68}
W. Firey, \emph{Christoffel's problem for general convex bodies}, Mathematika 15 (1968), 7-21.

\bibitem{F70}
W. Firey, \emph{Convex bodies with constant outer $p$-measure}, Mathematika 17 (1970), 21-27.

\bibitem{F74}
W. Firey,
\emph{Shapes of worn stones},
Mathematika 21 (1974), 1-11.

\bibitem{F70.1}
W. Firey, \emph{Intermediate Christoffel-Minkowski problems for figures of revolution},
Israel J. Math. 8 (1970), 384-390.

\bibitem{Ga93}
M. Gage,
\emph{Evolving plane curves by curvature in relative geometries},
Duke Math. J. 72 (1993), 441-466.




\bibitem{GL94}
  M. Gage, Y. Li,
  \emph{Evolving plane curves by curvature in relative geometries, II},
  Duke Math. J. 75 (1994), 79-98.

\bibitem{G06book}
R.J. Gardner,
\emph{Geometric Tomography},
Second Edition, Encyclopedia of Mathematics and its Applications,
Cambridge Univ. Press, New York, 2006.

\bibitem{G94annals}
R.J. Gardner,
\emph{A positive answer to the Busemann-Petty problem in three dimensions},
Ann. Math. 140 (1994), 435-447.

\bibitem{GKS99}
R.J. Gardner, A. Koldobsky, T. Schlumprecht,
\emph{An analytic solution to the Busemann-Petty problem on sections of convex bodies},
Ann.  Math. 149 (1999), 691-703.

\bibitem{G02bull}
R.J. Gardner, \emph{The Brunn-Minkowski inequality},
Bull. Amer. Math. Soc. 39 (2002), 355-405.

\bibitem{GHW14}
R.J. Gardner, D. Hug, W. Weil,
\emph{The Orlicz-Brunn-Minkowski theory: a general framework, additions, and inequalities},
J. Differential Geom. 97 (2014), 427-476.

\bibitem{GHWXY}
R.J. Gardner, D. Hug, W. Weil, S. Xing, D. Ye,
\emph{General volumes in the Orlicz-Brunn-Minkowski
theory and a related Minkowski problem I},
Calc. Var. PDE 58 (2019),  No. 12, 35 pp.

\bibitem{GHXY}
R.J. Gardner, D. Hug, S. Xing, D. Ye,
\emph{General volumes in the Orlicz-Brunn-Minkowski theory
and a related Minkowski problem II},  Calc. Var. PDE, 59 (2020), No. 15, 33 pp.

\bibitem{GZ98}
R.J. Gardner, G. Zhang,
\emph{Affine inequalities and radial mean bodies}, Amer. J. Math.
 120 (1998), 505-528.

\bibitem{Gl1}
H. Gluck,
\emph{The generalized Minkowski problem in differential geometry in the large},
Ann. Math.  96 (1972), 245-276.

\bibitem{Gl2}
H. Gluck,
\emph{Manifolds with preassigned curvature - a survey},
Bull. Amer. Math. Soc.  81 (1975), 313-329.


\bibitem{GW92}
P. Goodey, W. Weil, \emph{Centrally symmetric convex bodies and the spherical Radon
transform}, J. Differential Geom. 35 (1992), 675-688.


\bibitem{GYY11}
P. Goodey, V. Yaskin, M. Yaskina,
\emph{A Fourier transform approach to Christoffel's problem},
Trans. Amer. Math. Soc. 363 (2011), 6351-6384.



\bibitem{GZ99}
E. Grinberg, G. Zhang,
\emph{Convolutions, transforms, and convex bodies},
Proc. London Math. Soc. 78 (1999), 77-115.



\bibitem{GG02}
B. Guan, P. Guan,
\emph{Convex hypersurfaces of prescribed curvatures},
Ann.  Math. 156 (2002), 655-673.



\bibitem{GL97com}
P. Guan, Y. Li,
\emph{$C^{1,1}$ estimates for solutions of a problem of Alexandrov},
Comm. Pure Appl. Math. 50 (1997), 189-811.

\bibitem{GLM09}
P. Guan, C.S. Lin, X. Ma,
\emph{The Existence of Convex Body with Prescribed Curvature Measures},
Int. Math. Res. Not.  (2009), 1947-1975.


\bibitem{GLL12duke}
P. Guan, J. Li, Y.Y. Li,
\emph{Hypersurfaces of prescribed curvature measure},
Duke Math. J. 161 (2012), 1927-1942.



\bibitem{GM03}
P. Guan, X. Ma,
\emph{The Christoffel-Minkowski problem I: Convexity of solutions of a Hessian equation},
Invent. Math. 151 (2003), 553-577.



\bibitem{GN17}
P. Guan, L. Ni,
\emph{Entropy and a convergence theorem for Gauss curvature flow in high dimension},
J. Eur. Math. Soc.
19 (2017), 3735-3761.

\bibitem{GX18}
P. Guan, C. Xia,
 \emph{$L_p$ Christoffel-Minkowski problem: the case $1<p<k+1$},
 Calc. Var. PDE 57 (2018), no. 2, Paper No. 69.


\bibitem{GLW23}
Q. Guang, Q.-R. Li, X.-J. Wang,
\emph{Flow by Gauss curvature to the $L_p$ dual Minkowski problem},
Math. Eng. 5 (2023), no. 3, Paper No. 049.

\bibitem{GLW23a}
Q. Guang, Q.-R. Li, X.-J. Wang,
\emph{The $L_p$-Minkowski problem with super-critical exponents},
arXiv:2203.05099.

\bibitem{GXZ24}
L. Guo, D. Xi, Y. Zhao,
\emph{The $L_p$ chord Minkowski problem in a critical interval},
Math. Ann. 389 (2024), 3123-3162.



\bibitem{HS09jdg}
C. Haberl, F. Schuster,
\emph{General $L_p$ affine isoperimetric inequalities},
J. Differential Geom.  83 (2009), 1-26.


\bibitem{HS09jfa}
C. Haberl, F. Schuster,
\emph{Asymmetric affine $L_p$ Sobolev inequalities},
J. Funct. Anal. 257 (2009), 641-658.



\bibitem{HLYZ10}
C. Haberl, E. Lutwak, D. Yang, G. Zhang,
\emph{The even Orlicz Minkowski problem},
Adv. Math.  224 (2010), 2485-2510.

\bibitem{HSX12}
C. Haberl, F. Schuster, J. Xiao,
\emph{An asymmetric affine P\'olya-Szeg\"o principle},
Math. Ann. 352 (2012), 517-542.

\bibitem{HLL}
B. He, G. Leng, K. Li,
\emph{Projection problems for symmetric polytopes},
Adv. Math. 207 (2006), 73-90.

\bibitem{HLW16}
Y. He, Q.R. Li, X.J. Wang,
\emph{Multiple solutions of the $L_p$-Minkowski problem},
Calc. Var. PDE 55 (2016), Art. 117.

\bibitem{HSW}
M. Henk, A. Sch\"urmann, J.M. Wills,
\emph{Ehrhart polynomials and successive minima},
Mathematika 52 (2005), 1-16.

\bibitem{HL14}
M. Henk, E. Linke,
\emph{Cone-volume measures of polytopes},
Adv. Math. 253 (2014), 50-62.



\bibitem{HNS19}
M. Hering, B. Nill, H. S\"uss,
\emph{Stability of tangent bundles on smooth toric Picard-rank-2 varieties and surfaces},
Facets of algebraic geometry. Vol. II, 1-25, London Math. Soc. Lecture Note Ser., 473,
Cambridge Univ. Press, Cambridge, 2022.


\bibitem{HYZ18}
H. Hong, D. Ye, N. Zhang,
\emph{The $p$-capacitary Orlicz-Hadamard variational formula and Orlicz-Minkowski problems},
Calc. Var. PDE 57 (2018), https://doi.org/10.1007/s00526-017-1278-6.


\bibitem{HMS04}
C. Hu, X. Ma, C. Shen,
\emph{On the Christoffel-Minkowski problem of Firey's $p$-sum},
Calc. Var. PDE 21 (2004), 137-155.


\bibitem{HHL24}
J. Hu, Y. Huang, J. Lu,
\emph{Boundary regularity of Riesz potential, smooth solution to the chord log-Minkowski problem},
 (2024).




\bibitem{HH12}
Q. Huang, B. He,
\emph{On the Orlicz Minkowski problem for polytopes},
Discrete Comput. Geom. 48 (2012), 281-297.


\bibitem{HLYZ16acta}
Y. Huang, E. Lutwak, D. Yang,  G. Zhang,
\emph{Geometric measures in the dual Brunn-Minkowski theory and their
associated Minkowski problems}, Acta Math. 216 (2016), 325-388.

\bibitem{HLYZ18jdg}
Y. Huang, E. Lutwak, D. Yang,  G. Zhang,
\emph{The $L_p$-Aleksandrov problem for $L_p$-integral curvature},
J. Differential Geom. 110 (2018), 1-29.

\bibitem{HJ19}
Y. Huang, Y. Jiang,
\emph{Variational characterization for the planar dual Minkowski problem},
J. Funct. Anal. 277 (2019), 2209-2236.

\bibitem{HLX15}
Y. Huang, J. Liu, L. Xu,
\emph{On the uniqueness of $L_p$-Minkowski problems: the constant $p$-curvature case in $\R^3$},
Adv. Math. 281 (2015), 906-927.

\bibitem{HQ24}
Y. Huang, L. Qin,
\emph{The Gaussian chord Minkowski problem},
Discrete Contin. Dyn. Syst. Ser. S 17 (2024), no. 2, 930-944.


\bibitem{HXZ20}
Y. Huang, D. Xi, Y. Zhao,
\emph{The Minkowski problem in Gaussian probability space},
Adv. Math. 385 (2021), Paper No. 107769.

\bibitem{HZ18}
Y. Huang, Y. Zhao,
\emph{On the $L_p$-dual Minkowski problem},
Adv. Math. 322 (2018), 57-84.




\bibitem{HLYZ05dcg}
D.~Hug, E.~Lutwak, D.~Yang,  G.~Zhang,
\emph{On the $L\sb p$ Minkowski problem for polytopes},
Discrete Comput. Geom. 33 (2005), 699-715.

\bibitem{Je96}
D. Jerison,
\emph{A Minkowski problem for electrostatic capacity}, Acta Math.  176 (1996), 1-47.

\bibitem{Je96adv}
D. Jerison,
\emph{The direct method in the calculus of variations for convex bodies},
Adv. Math. 122 (1996), 262-279.

\bibitem{JL19}
H. Jian, J. Lu,
\emph{Existence of solutions to the Orlicz-Minkowski problem},
Adv. Math. 344 (2019), 262-288.

\bibitem{JLW15}
H. Jian, J. Lu, X.-J. Wang,
\emph{Nonuniqueness of solutions to the $L_p$-Minkowski problem},
Adv. Math. 281 (2015), 845-856.

\bibitem{JLZ16}
H. Jian, J. Lu, G. Zhu,
\emph{Mirror symmetric solutions to the centro-affine Minkowski problem},
Calc. Var. PDE 55 (2016), Art. 41, 22 pp.



\bibitem{Klain}
D. Klain,
\emph{The Minkowski problem for polytopes},
Adv. Math. 185 (2004), 270-288.

\bibitem{K03}
A. Koldobsky,
\emph{Intersection bodies, positive definite distributions, and the Busemann-Petty problem},
Amer. J. Math. 120 (1998), 827-840.


\bibitem{K05}
A. Koldobsky,
\emph{Fourier Analysis in Convex Geometry}, Mathematical Surveys and Monographs,
 116, Amer. Math. Soc., Providence, RI, 2005.


\bibitem{KM18}
A.V. Kolesnikov, E. Milman,
\emph{Local $L_p$-Brunn-Minkowski inequalities for $p < 1$},
Mem. Amer. Math. Soc. 277 (2022), no. 1360.

\bibitem{KL23}
L. Kryvonos, D. Langharst,
\emph{Weighted Minkowski existence theorem and projection bodies},
Trans. Amer. Math. Soc. 376 (2023), no. 12, 8447-8493.

\bibitem{La02}
R. Lata\l{}a,
\emph{On some inequalities for Gaussian measures},
Proc. Int. Cong. Math. 2002, Vol. II, pp. 813-822.

\bibitem{Le77}
J.L. Lewis,
\emph{Capacitary functions in convex rings},
Arch. Rational Mech. Anal.  66 (1977), 201-224.

\bibitem{LN}
J.L. Lewis, K. Nystr{\"o}m,
\emph{Boundary behaviour for $p$-harmonic functions in Lipschitz and starlike Lipschitz ring domains},
Ann. Sci. \' {E}cole Norm. Sup.  40 (2007), 765-813.



\bibitem{Lew38}
H. Lewy,
\emph{On differential geometry in the large I (Minkowski problem)},
Trans. Amer. Math. Soc. 43 (1938), 258-270.

\bibitem{Li14}
A.-J. Li,
\emph{The generalization of Minkowski problems for polytopes},
Geom. Dedicata 168 (2014), 245-264.

\bibitem{LiYZ}
N. Li, D. Ye, B. Zhu,
\emph{The dual Minkowski problem for unbounded closed convex sets},
Math. Ann. 388 (2024), no. 2, 2001-2039.




\bibitem{LiQ19}
Q.-R. Li,
\emph{Infinitely many solutions for centro-affine Minkowski problem},
Int. Math. Res. Not. (2019), no. 18, 5577-5596.

\bibitem{LLL21}
Q.-R. Li, J. Liu, J. Lu,
\emph{Nonuniqueness of solutions to the $L_p$ dual Minkowski problem},
Int. Math. Res. Not. (2022), no. 12, 9114-9150.

\bibitem{LSW20}
Q.-R. Li, W. Sheng, X.-J. Wang,
\emph{Flow by Gauss curvature to the Aleksandrov and dual Minkowski problems},
J. Eur. Math. Soc. 22 (2020), 893-923.

\bibitem{LSW20jga}
Q.-R. Li, W. Sheng, X.-J. Wang,
\emph{Asymptotic convergence for a class of fully nonlinear curvature flows},
J. Geom. Anal. 30 (2020), 834-860.

\bibitem{LWW20}
Q.-R. Li, D. Wan, X.-J. Wang,
\emph{The Christoffel problem by the fundamental solution of the Laplace equation},
Science China in Mathematics, Special Issue on Differential Geometry, 2020.

\bibitem{LW18}
Q.-R. Li, X.-J. Wang,
\emph{A class of optimal transportation problems on the sphere}, (Chinese),
Scientia Sinica Mathematica, 48 (2018), 181-200.


\bibitem{Li24}
Y. Li,
\emph{Nonuniqueness of solutions to the $L_p$ chord Minkowski problem},
Calc. Var. PDE 63 (2024), no. 4, Paper No. 83, 22 pp.

\bibitem{Li89}
Y.-Y. Li,
\emph{Degree theory for second order nonlinear elliptic operators and its applications},
Comm. Partial Differential Equations 14 (1989), 1541-1578.


\bibitem{Liu22}
J. Liu,
\emph{The $L_p$-Gaussian Minkowski problem},
Calc. Var. PDE 61 (2022), no. 1, Paper No. 28.


\bibitem{LL20}
Y. Liu, J. Lu,
\emph{A flow method for the dual Orlicz-Minkowski problem},
Trans. Amer. Math. Soc. 373 (2020), 5833-5853.

\bibitem{LLSX24}
Y. Liu, X. Lu, Q. Sun, G. Xiong,
\emph{The logaritumic Minkowski problem in ${\mathbb R}^2$},
Pure Appl. Math. Q. 20 (2024), no. 2, 869-902.



\bibitem{Liv19}
G. Livshyts,
\emph{An extension of Minkowski's theorem and its applications to questions about projections for measures},
Adv. Math. 356 (2019), Paper No. 106803.





\bibitem{Lud03}
M. Ludwig,
\emph{Ellipsoids and matrix-valued valuations},
Duke Math. J. 119 (2003), 159-188.

\bibitem{Lud04}
M. Ludwig,
\emph{Intersection bodies and valuations},
Amer. J. Math. 128 (2006), 1409-1428.

\bibitem{Lud05}
M. Ludwig,
\emph{Valuations on Sobolev spaces},
Amer. J. Math. 134 (2012), 824-842.

\bibitem{Lud09}
M. Ludwig,
\emph{General affine surface areas},
Adv. Math. 224 (2009), 2346-2360.

\bibitem{Lud06}
M. Ludwig,
\emph{Minkowski areas and valuations},
J. Differential Geom.  86 (2010), 133-161.

\bibitem{Lud14}
M. Ludwig,
\emph{Anisotropic fractional perimeters},
J. Differential Geom. 96 (2014), 77-93.


\bibitem{LR10}
M. Ludwig, M. Reitzner,
\emph{A classification of $SL(n)$ invariant valuations},
Ann. of Math.  172 (2010), 1219-1267.

\bibitem{LXZ11}
M. Ludwig, J. Xiao, G. Zhang,
\emph{Sharp convex Lorentz-Sobolev inequalities},
Math. Ann. 350 (2011), 169-197.


\bibitem{L75paci}
E. Lutwak,
\emph{Dual mixed volumes},
Pacific J. Math. 58 (1975), 531-538.

\bibitem{L88adv}
E. Lutwak,
\emph{Intersection bodies and dual mixed volumes},
Adv. Math. 71 (1988), 232-261.

\bibitem{L93hcg}
E.~Lutwak,
\emph{Selected affine isoperimetric inequalities},
Handbook of Convex Geometry (P.M. Gruber, J.M. Wills, eds), vol. A, pp. 161-176, North-Holland, 1993.

\bibitem{L93jdg}
E.~Lutwak,
\emph{The Brunn-Minkowski-Firey theory. I. Mixed volumes and the Minkowski problem},
J. Differential Geom. 38 (1993), 131-150.

\bibitem{L96adv}
E. Lutwak,
\emph{The Brunn-Minkowski-Firey theory. II. Affine and geominimal surface areas},
Adv. Math. 118 (1996), 244-294.

\bibitem{LO95jdg}
E.~Lutwak, V.~Oliker,
\emph{On the regularity of solutions to a generalization of the Minkowski problem},
J. Differential Geom. 41 (1995), 227-246.


\bibitem{LXYZ21}
E.~Lutwak, D. Xi, D.~Yang, G.~Zhang,
\emph{Chord measures in integral geometry  and their Minkowski problems},
Comm. Pure Appl. Math.  77 (2024), no. 7,  3277-3330.


\bibitem{LYZ00jdg}
E.~Lutwak, D.~Yang, G.~Zhang,
\emph{${L}\sb p$ affine isoperimetric inequalities},
J. Differential Geom. 56 (2000), 111-132.



\bibitem{LYZ02jdg}
E.~Lutwak, D.~Yang, G.~Zhang,
\emph{Sharp affine {$L\sb p$} Sobolev inequalities},
J. Differential Geom. 62 (2002), 17-38.

\bibitem{LYZ04tams}
E.~Lutwak, D.~Yang, G.~Zhang,
\emph{On the $L\sb p$-Minkowski problem},
Trans. Amer. Math. Soc. 356 (2004), 4359-4370.





\bibitem{LYZ06imrn}
E.~Lutwak, D.~Yang, G.~Zhang,
\emph{Optimal Sobolev norms and the $L\sp p$ Minkowski problem},
Int. Math. Res. Not.  (2006), Art. ID 62987, 21 pp.

\bibitem{LYZ10jdg}
 E. Lutwak, D. Yang, G. Zhang,
 \emph{Orlicz centroid bodies},
 J. Differential Geom. 84 (2010), 365-387.

\bibitem{LYZ10adv}
E. Lutwak, D. Yang, G. Zhang,
\emph{Orlicz projection bodies},
Adv.  Math. 223 (2010), 220-242.


\bibitem{LYZ18adv}
E. Lutwak, D. Yang, G. Zhang,
\emph{$L_p$ dual curvature measures},
Adv.  Math. 329 (2018), 85-132.



\bibitem{LuoYeZhu}
X. Luo, D. Ye, B. Zhu,
\emph{On the Polar Orlicz-Minkowski Problems and the p-capacitary Orlicz-Petty Bodies},
Indiana Univ. Math. J. 69 (2020), 385-420.


\bibitem{Ma15}
L. Ma,
\emph{A new proof of the log-Brunn-Minkowski inequality},
Geom. Dedicata  177 (2015), 75-82.

\bibitem{MW00}
M. Meyer, E. Werner,
\emph{On the $p$-affine surface area},
Adv. Math. 152 (2000), 288-313.

\bibitem{M22jdg}
E. Milman,
\emph{A sharp centro-affine isospectral inequality of Szeg\"o-Weinberger type and
the $L_p$-Minkowski problem},
J. Differential Geom. 127 (2024), 373-408.




\bibitem{M1897}
H. Minkowski,
\emph{Allgemeine Lehrs\"atze \"uber die convexen Polyeder},
 Nachr. Ges. Wiss. G\"ottingen, 198-219.
 Gesammelte Abhandlungen, Vol. II, Teubner, Leipzig, 1911, pp. 103-121.

\bibitem{M1903}
H. Minkowski,
\emph{Volumen und Oberfl\"ache}, Math. Ann.  57 (1903),  447-495.

\bibitem{Mui}
S. Mui,
\emph{On the $L^p$ Aleksandrov problem for negative $p$},
Adv. Math. 408 (2022), Part A, Paper No. 108573.



\bibitem{N03}
F. Nazarov,
\emph{On the maximal perimeter of a convex set in $\rn$ with respect to a Gaussian measure},
In Geometric Aspects of Functional Analysis,  Lecture Notes in Math. Vol. 1807, 169-187,
Springer, Berlin, 2003.


\bibitem{N53}
L. Nirenberg,
\emph{The Weyl and Minkowski problems in differential geometry in the large},
Comm. Pure Appl. Math.  6 (1953), 337-394.


\bibitem{NS94}
K. Nomizu, T. Sasaki,
\emph{Affine differential geometry},
Cambridge Univ. Press, New York, 1994.


\bibitem{Ol}
V. Oliker,
\emph{Embedding $\sn$ into $\R^{n+1}$ with given integral Gauss curvature
and optimal mass transport on $\sn$},
Adv. Math. 213 (2007), 600-620.

\bibitem{Ol2}
V. Oliker,
\emph{Existence and uniqueness of convex hypersurfaces with prescribed Gaussian curvature
in spaces of constant curvature},
Sem. Inst. Matem. Appl. Giovanni Sansone (1983), 1-64.

\bibitem{P69}
A.V. Pogorelov,
\emph{Extrinsic geometry of convex surfaces}, ``Nauka", Moscow, 1969; English
transl., Transl. Math. Mono., Vol. 35, Amer. Math. Soc., Providence, R.I., 1973.


\bibitem{P78}
A.V. Pogorelov,
\emph{The Minkowski Multidimensional Problem},
V.H. Winston \& Sons, Washington, D.C., 1978.

\bibitem{Q22}
L. Qin,
\emph{Nonlocal energies of convex body and their log-Minkowski problem},
 Adv. Math. 427 (2023), Paper No. 109132.

\bibitem{Q24}
L. Qin,
\emph{The chord log-Minkowski problem for $0<q<1$},
Proc. Amer. Math. Soc. 152 (2024), 2647-2655.


\bibitem{Ren}
D. Ren,
\emph{Topics in Integral Geometry},
World Scientific, Singapore, 1994.

\bibitem{San}
 L.A. Santalo,
 \emph{Integral Geometry and Geometric Probability},
 Addison-Wesley, Reading, MA, 1976.


\bibitem{Saro1}
C. Saroglou,
\emph{Remarks on the conjectured log-Brunn-Minkowski inequality},
Geom. Dedicata 177 (2015), 353-365.


\bibitem{S14}
R. Schneider,
\emph{Convex Bodies: The Brunn-Minkowski Theory},
 Second Edition, Encyclopedia of Mathematics and its Applications, Cambridge
 Univ. Press, New York, 2014.

\bibitem{S77}
R. Schneider,
\emph{Das Christoffel-Problem f\"ur Polytope},
Geom. Dedicata 6 (1977), 81-85.

\bibitem{S78}
R. Schneider,
\emph{Curvature measures of convex bodies},
Ann. Mat. Pura Appl. 116 (1978), 101-134.

\bibitem{S18}
R. Schneider,
\emph{A Brunn-Minkowski theory for coconvex sets of finite volume},
Adv. Math. 332 (2018), 199-234.

\bibitem{S23}
R. Schneider,
\emph{Pseudo-cones},
Adv. Appl. Math. 155 (2024), Paper No. 102657.

\bibitem{SW04}
C. Sch\"utt, E. Werner,
\emph{Surface bodies and $p$-affine surface area},
Adv. Math. 187 (2004), 98-145.

\bibitem{S22}
V. Semenov,
\emph{The Discrete Gauss Image Problem},
(2022), arXiv:2210.16974.

\bibitem{S24}
V. Semenov,
\emph{The Gauss image problem with weak Aleksandrov condition},
(2023), arXiv:2210.16778v2.



\bibitem{STW}
W. Sheng, N. Trudinger, X.-J. Wang,
\emph{Convex hypersurfaces of prescribed Weingarten curvatures},
Comm. Anal. Geom.  12 (2004), 213-232.

\bibitem{SX22}
W. Sheng, K. Xue,
\emph{Flow by Gauss curvature to the $L_p$-Gaussian Minkowski problem},
 arXiv:2212.01822.

\bibitem{St02}
A. Stancu,
\emph{The discrete planar $L_0$-Minkowski problem},
Adv. Math. 167 (2002), 160-174.

\bibitem{St03}
A. Stancu,
\emph{On the number of solutions to the discrete
two-dimensional $L_0$-Minkowski problem},
Adv. Math.  180 (2003), 290-323.

\bibitem{Sun18}
Y. Sun,
\emph{Existence and uniqueness of solutions to Orlicz Minkowski problems involving $0 < p < 1$},
Adv. Appl. Math. 101 (2018), 184-214.

\bibitem{SL15}
Y. Sun, Y. Long,
\emph{The planar Orlicz Minkowski problem in the $L_1$-sense},
Adv. Math. 281 (2015), 1364-1383.

\bibitem{SZ19}
Y. Sun, D. Zhang,
\emph{The planar Orlicz Minkowski problem for $p = 0$ without even assumptions},
J. Geom. Anal. 29 (2019), 3384-3404.


\bibitem{T85}
K. Tso,
\emph{Deforming a hypersurface by its Gauss-Kronecker curvature},
Comm. Pure Appl. Math. 38 (1985), 867-882.



\bibitem{W07}
E. Werner,
\emph{On $L_p$-affine surface areas},
Indiana Univ. Math. J. 56 (2007), 2305-2323.

\bibitem{WY08adv}
E. Werner, D. Ye,
\emph{New $L_p$ affine isoperimetric inequalities}, Adv. Math. 218 (2008), 762-780.

\bibitem{WY10mathann}
E. Werner, D. Ye,
\emph{Inequalities for mixed $p$-affine surface area},
Math. Ann. 347 (2010), 703-737.

\bibitem{WXL19}
Y. Wu, D. Xi, G. Leng,
\emph{On the discrete Orlicz Minkowski problem},
Trans. Amer. Math. Soc. 371 (2019), 1795-1814.

\bibitem{WXL20}
Y. Wu, D. Xi, G. Leng,
\emph{On the discrete Orlicz Minkowski problem II},
Geom. Dedicata 205 (2020), 177-190.

\bibitem{XL16}
D. Xi, G. Leng, \emph{Dar's conjecture and the log-Brunn-Minkowski inequality},
J. Differential Geom. 103 (2016), 145-189.

\bibitem{XJL14}
D. Xi, H. Jin, G. Leng,
\emph{The Orlicz Brunn-Minkowski inequality},
Adv. Math. 260 (2014), 350-374.

\bibitem{XYZZ23}
D. Xi, D. Yang, G. Zhang, Y. Zhao,
\emph{The $L_p$ chord Minkowski problem},
Advanced Nonlinear Studies, 23 (2023), no. 1, pp. 20220041.



\bibitem{X07}
J. Xiao,
\emph{The sharp Sobolev and isoperimetric inequalities split twice},
Adv. Math. 211 (2007), 417-435.

\bibitem{X15}
J. Xiao,
\emph{On the variational $p$-capacity problem in the plane},
Comm. Pure Appl. Anal. 14 (2015), 959-968.

\bibitem{X16}
J. Xiao,
\emph{The $p$-affine capacity},
J. Geom. Anal. 26 (2016), 947-966.

\bibitem{X17}
J. Xiao,
\emph{The $p$-affine capacity redux},
J. Geom. Anal. 27 (2017), 2872-2888.

\bibitem{X18}
J. Xiao,
\emph{Exploiting log-capacity in convex geometry},
Asian J. Math. 22 (2018), 955-980.

\bibitem{X20}
J. Xiao,
\emph{Prescribing capacitary curvature measures on planar convex domains},
J. Geom. Anal. 30 (2020), 2225-2240.

\bibitem{Xie}
F. Xie,
\emph{The Orlicz Minkowski problem for general measures},
Proc. Amer. Math. Soc. 150 (2022), no. 10, 4433-4445.

\bibitem{Xie2}
F. Xie,
\emph{The Orlicz Minkowski problem for cone-volume measure},
Adv. Appl. Math. 149 (2023), Paper No. 102523.





\bibitem{XY20}
S. Xing, D. Ye,
\emph{The general dual Orlicz-Minkowski problem},
Indiana Univ. Math. J. 69 (2020), 621-655.


\bibitem{X1}
G. Xiong,
\emph{Extremum problems for the cone-volume functional of convex polytopes},
Adv. Math. 225 (2010), 3214-3228.

\bibitem{XXL19}
G. Xiong, J. Xiong, L. Xu,
\emph{The $L_p$ capacitary Minkowski problem for polytopes},
J. Funct. Anal. 277 (2019), 3131-3155.

\bibitem{XLL23}
Y. Xu, J. Li, G. Leng,
\emph{Dualities and endomorphisms of pseudo-cones},
Adv. Appl. Math. 142 (2023), Paper No. 102434, 31 pp.

\bibitem{YYZ23}
J. Yang, D. Ye, B. Zhu,
\emph{On the $L_p$ Brunn-Minkowski theory and the $L_p$ Minkowski problem for $C$-coconvex sets},
Int. Math. Res. Not. (2023), no. 7, 6252-6290.

\bibitem{Z94}
G. Zhang,
\emph{Centered bodies and dual mixed volumes},
Trans. Amer. Math. Soc.  345 (1994), 777-801.

\bibitem{Z99tams}
G. Zhang,
\emph{Dual kinematic formulas},
Trans. Amer. Math. Soc. 351 (1999), 985-995.

\bibitem{Z99annals}
G. Zhang,
\emph{A positive solution to the Busemann-Petty problem in $\R^4$},
Ann. Math. 149 (1999), 535-543.


\bibitem{Z99jdg}
G. Zhang,
\emph{The affine Sobolev inequality},
J. Differential Geom. 53 (1999), 183-202.



\bibitem{Zhao}
Y. Zhao,
\emph{Existence of solutions to the even dual Minkowski problem},
J. Differential Geom. 110 (2018), 543-572.

\bibitem{Zhao2} Y. Zhao,
\emph{The dual Minkowski problem for negative indices},
Calc. Var. PDE 56 (2017).

\bibitem{Zhao3} Y. Zhao,
\emph{The $L_p$ Aleksandrov problem for origin-symmetric polytopes},
Proc. Amer. Math. Soc.  147 (2019), 4477-4492.

\bibitem{ZXY18}
B. Zhu, S. Xing, D. Ye,
\emph{The Dual Orlicz-Minkowski Problem},
J. Geom. Anal. 28 (2018), 3829-3855.

\bibitem{ZZW14}
B. Zhu, J. Zhou, W. Xu,
\emph{Dual Orlicz-Brunn-Minkowski theory},
Adv. Math. 264 (2014), 700-725.

\bibitem{Zhu1} G. Zhu,
\emph{The logarithmic Minkowski problem for polytopes},
Adv. Math. 262 (2014), 909-931.


\bibitem{Zhu3} G. Zhu,
\emph{The centro-affine Minkowski problem for polytopes},
J. Differential Geom. 101 (2015), 159-174.

\bibitem{Zhu4} G. Zhu,
\emph{The $L_p$ Minkowski problem for polytopes for $p<0$},
Indiana Univ. Math. J. 66 (2017), 1333-1350.

\bibitem{ZX20}
D. Zou, G. Xiong,
\emph{The $L_p$ Minkowski problem for the electrostatic $\mathfrak p$-capacity},
J. Differential Geom. 116 (2020), 555-596.



\end{thebibliography}
\end{document}